\documentclass{article}
   \usepackage{amssymb} 
   \usepackage[all,cmtip]{xy}
   \usepackage[width=440pt, height=650pt, hcentering, vcentering]{geometry}
\newcommand{\g}{\mbox{$\bf g$}}

\newcommand{\h}{\mbox{\textbf{h}}}
\newcommand{\n}{\mbox{\textbf{n}}}
\newcommand{\np}{\mbox{$\textbf{n}^+$}}
\newcommand{\nm}{\mbox{$\textbf{n}^-$}}
\newcommand{\al}{\alpha}
\newcommand{\eps}{\epsilon}
\newcommand{\la}{\lambda}
\newcommand{\La}{\Lambda}

\newcommand{\Th}{\Theta}

\newcommand{\mb}{\mbox}
\newcommand{\Mkl}[1]{\left\{ #1 \right\}}
\newcommand{\Mklz}[2]{\left\{\left. #1 \right| #2 \right\}}
\newcommand{\W}{\mbox{$\Delta$}}

\newcommand{\rW}{\mbox{$\Delta_{re}$}}

\newcommand{\iW}{\mbox{$\Delta_{im}$}}

\newcommand{\F}{\mathbb{F}}

\newcommand{\C}{\mathbb{C}}
\newcommand{\N}{\mathbb{N}}
\newcommand{\Nn}{\mathbb{N}_0}
\newcommand{\Q}{\mathbb{Q}}

\newcommand{\R}{\mathbb{R}}

\newcommand{\Z}{\mathbb{Z}}
\newcommand{\We}{\mbox{$\mathcal W$}}

\newcommand{\RkX}{\mbox{$\mb{Fa}(X)$}}
\newcommand{\iB}[2]{\left(#1\mid#2\right)}
\newcommand{\iBl}{\left( \;\ \mid \;\ \right)}

\newcommand{\ND}{\mbox{$\widehat{N}$}}
\newcommand{\GD}{\mbox{$\widehat{G}$}}

\newcommand{\WeD}{\mbox{$\widehat{\We}$}}
\newcommand{\PD}{\mbox{$\widehat{P}$}}

\newcommand{\ve}[1]{\mbox{$\varepsilon (#1)$}} 
\newcommand{\FK}[1]{\mbox{$\F\,[#1]$}}

\newcommand{\PCAF}{\mbox{Proj\,}(CA)(\mathbb{F})}
\newcommand{\PCA}{\mbox{Proj\,}(CA)}
\newcommand{\PA}{\mbox{Proj\,}(A)}
\newcommand{\PAF}{\mbox{Proj\,}(A)(\mathbb{F})}
\newcommand{\PB}{\mbox{Proj\,}(B)}
\newcommand{\PBF}{\mbox{Proj\,}(B)(\mathbb{F})}

\newcommand{\prd}{\mbox{\boldmath$\,\cdot\,$}}
\newcommand{\prdca}{\mbox{\boldmath$\,\cdot\,$}}
\newcommand{\ins}{\,\cdot\,}
\newcommand{\ti}{\widetilde}
\newcommand{\res}[1]{\!\mid_{#1}}

\newcommand{\red}[1]{red\left(#1\right)}

\newcommand{\Proof}{\mbox{\bf Proof: }}
\newcommand{\qed}{\mb{}\hfill\mb{$\square$}\\}

\newtheorem{thm}{Theorem}[section]

\newtheorem{prop}[thm]{Proposition}
\newtheorem{cor}[thm]{Corollary}
\newtheorem{rem}[thm]{Remark}

\newtheorem{lem}[thm]{Lemma}

\begin{document}

\title{An algebraic geometric model of an action of the face monoid associated to a Kac-Moody group on its building}
\author{Claus Mokler} 
\date{}
\maketitle
\begin{abstract}\noindent 
The face monoid $\GD$ described in \cite{M1}  acts on the integrable irreducible highest weight modules of a symmetrizable Kac-Moody algebra. It has similar structural properties as a reductive algebraic monoid whose unit group is a symmetrizable Kac-Moody group $G$. We found in \cite{M5} two natural extensions of the action of the Kac-Moody group $G$ on its building $\Omega$ to actions of the face monoid $\GD$ on the building $\Omega$. Now we give an algebraic geometric model of one of these actions of the face monoid $\GD$ on $\Omega$, where the building $\Omega$ is obtained as a part of the $\F$-valued points of the spectrum of all homogeneous prime ideals of the Cartan algebra $CA$ of the Kac-Moody group $G$. We describe the the spectrum of all homogeneous prime ideals of the Cartan algebra $CA$ and determine its $\F$-valued points.
\end{abstract}
\section*{Introduction}
The face monoid $\GD$ is an infinite-dimensional algebraic monoid. It has been obtained in \cite{M1} by a Tannaka reconstruction from categories determined by the integrable irreducible highest weight representations of a symmetrizable Kac-Moody algebra. By its construction and by the involved categories it is a very natural object. Its Zariski open dense unit group $G$ coincides, up to a slightly extended maximal torus, with the Kac-Moody group defined representation theoretically in \cite{KP1}.  

The face monoid is a purely infinite-dimensional phenomenon, quite unexpected. In the classical case, i.e., if one takes a split semisimple Lie algebra for the symmetrizable Kac-Moody algebra, it coincides with a split semisimple simply connected algebraic group. Put in another way, there seem to exist fundamentally different infinite-dimensional generalizations of a split semisimple simply connected algebraic group.

The results obtained in \cite{M1}, \cite{M2}, \cite{M3}, and \cite{M4} show that the face monoid $\GD$ has similar structural and algebraic geometric properties as a reductive algebraic monoid, e.g., the monoid of ($n\times n$)-matrices. The face monoid $\GD$ is the first example of an infinite-dimensional reductive algebraic monoid. Actually, it is particular. The investigation of the conjugacy classes in \cite{M6} will show, that the relation between the face monoid $\GD$ and its unit group $G$, the Kac-Moody group, is much closer than for a general reductive algebraic monoid.

Obviously, there is the following question: Does the face monoid $\GD$ fit in some way to the building theory of the Kac-Moody group $G$? In \cite{M5} we investigated how to extend the natural action of the Kac-Moody group $G$ on its building $\Omega$ to actions of the face monoid $\GD$ on $\Omega$. To explain the results obtained in \cite{M5} note the following facts: 

For the face monoid $\GD$ an infinite Renner monoid $\WeD$ plays the same role as the Weyl group $\We$ does for the Kac-Moody group $G$. For example, there are Bruhat and Birkhoff decompositions of $\GD$, similar as for $G$, but the Weyl group $\We$ replaced by the monoid $\WeD$. The monoid $\WeD$ can be constructed from the Weyl group $\We$ and the face lattice of the Tits cone $X$, where the term "face" means a face of the convex cone $X$ in the sense of convex geometry. The Weyl group $\We$ is the unit group of $\WeD$.

The building $\Omega$ is covered by certain subcomplexes, the apartments, which are isomorphic to the Coxeter complex $\mathcal C$ associated to the Weyl group $\We$. This connects the action of $G$ on $\Omega$ and the action of $\We$ on $\mathcal C$.

Now let $\mathcal A$ be the standard apartment associated to a fixed BN-pair $(B,\,N)$ of the Kac-Moody group $G$. The group $N$ may be obtained as normalizer $N=N_G(T)$ of the standard torus $T=B\cap N$, and the Weyl group $\We$ identifies with $N/T$. If we define similarly $\ND:=N_{\widehat{G}}(T)$, then the monoid $\WeD$ identifies with $\ND/T$. Consider the following diagram:
\begin{eqnarray*}
\xymatrix{ \Omega  \ar@{}[d] |{\mb{$ \circlearrowleft$}} \ar@{}[r] |{\mb{$\supset$}} & {\mathcal A}   \ar@{}[d] |{\mb{$ \circlearrowleft$}}\ar@{}[r] | {\mb{$\cong$} } & {\mathcal C}  \ar@{}[d] |{\mb{$ \circlearrowleft$}} \\
G  \ar@{}[d] |{\mb{$ \cap$}}\ar@{}[r] |{\mb{$\supset$}} & N  \ar@{}[d] |{\mb{$ \cap$}}\ar@{->>}[r] & \We \ar@{}[d] |{\mb{$ \cap$}}\\ 
\GD\,?\ar@{}[r] |{\mb{$\supset$}} & \ND\ar@{->>}[r] & \WeD
}
\end{eqnarray*}
Suppose it is given an action of $\WeD$ on $\mathcal C$, extending the natural action of $\We$ on $\mathcal C$. Then this action induces uniquely an action of $\ND$ on $\mathcal A$, extending the natural action of $N$ on $\mathcal A$. 
There are the following questions:
\begin{enumerate}
\item[$\bullet$] Does there exist an action of $\GD$ on $\Omega$, extending the natural action of $G$ on $\Omega$, such that the diagram commutes? 
\item[$\bullet$]  If such an action exists, is it uniquely determined?
\item[$\bullet$]  Which actions can be obtained in this way?
\end{enumerate}
In Theorem 45 and Corollary 48 of \cite{M5} we obtained the following beautiful result: There is a bijective correspondence between:
\begin{itemize}
\item[(i)] The actions of $\WeD$ on $\mathcal C$, extending the natural action of $\We$ on $\mathcal C$, satisfying certain conditions which we do not state here.
\item[(ii)] The actions of $\GD$ on $\Omega$, extending the natural action of $G$ on $\Omega$. 
\end{itemize}
Furthermore, the action of $\GD$ on $\Omega$ is obtained by an explicit formula from the action of $\WeD$ on $\mathcal C$.

There exist quite general actions of $\GD$ on $\Omega$, extending the natural action of $G$ on $\Omega$. Compare for example the action of Remark 46 in \cite{M5}. There is the following question:
\begin{enumerate}
\item[$\bullet$]  Is it possible to single out some actions of $\GD$ on $\Omega$, which extend the natural action of $G$ on $\Omega$ and also keep some of its properties?
\end{enumerate}
The natural action of $G$ on the building $\Omega$ satisfies:
\begin{itemize}
\item[(a)] It preserves the natural order on $\Omega$.
\item[(b)] $\mb{Stab}_G(Q)= Q$ for all standard parabolic subgroups $Q\in\Omega$.
\end{itemize}
Now every standard parabolic subgroup $Q$ of $G$ extends naturally to a standard parabolic submonoid $\widehat{Q}$ of $\GD$. An action of $\GD$ on $\Omega$, which extends the natural action of $G$ on $\Omega$, is called good if it satisfies:
\begin{itemize}
\item[(a)] It preserves the natural order on $\Omega$.
\item[(b)] $\mb{Stab}_{\widehat{G}}(Q)= \widehat{Q}$ for all standard parabolic subgroups $Q\in\Omega$.
\end{itemize}
It is open to determine all good actions. In \cite{M5} we found two good actions of $\GD$, good action 1 and good action 2. Both have been obtained by the correspondence mentioned above from actions of $\WeD$ on $\mathcal C$, which in turn have been found by an action of $\WeD$ on the facets of the Tits cone and by an action of $\WeD$ on Looijenga's modified Tits cone. However, the formulas for these actions in Corollary 50 and Corollary 51 of \cite{M5} are unintuitive. It is desirable to find natural models for these actions, if at all possible.

In this paper we give a natural algebraic geometric model of good action 1 of $\GD$. It is obtained as follows: Let $P^+$ be the set of dominant weights of the weight lattice $P$. Let $L(\La)$ be the irreducible highest weight module of highest weight $\La$, let $L(\La)^{(*)}$ be its restricted dual. The Cartan algebra
\begin{eqnarray*}
  CA=\bigoplus_{\La\in P^+}L(\La)^{(*)}
\end{eqnarray*}
is a $P^+$-graded algebra over the field $\F$ of characteristic zero, whose multiplication on the graded parts
\begin{eqnarray*}
   L(\La_1)^{(*)}\otimes L(\La_2)^{(*)}\to L(\La_1+\La_2)^{(*)}
\end{eqnarray*}
is obtained dually to $G$-equivariant embeddings $L(\La_1+\La_2)\to L(\La_1)\otimes L(\La_2)$, $\La_1$, $\La_2\in P^+$.

The spectrum $\PCA$ of all $P^+$-graded prime ideals, no irrelevant ideals excluded, is a locally ringed space. It is not a scheme, but it is stratified by schemes. 
The building $\Omega$, which is as a set the union 
\begin{eqnarray*}
  \Omega =\dot{\bigcup_{Q\; \textrm{\footnotesize standard parabolic}   }} G/Q
\end{eqnarray*}
of all thin flag spaces, identifies with a dense part of the $\F$-valued points $\PCAF$ of this spectrum, i.e., there exists a natural injection
\begin{eqnarray*}
  \Omega \hookrightarrow\PCAF
\end{eqnarray*}
with dense image. A natural completion $\Omega_{fn}$ of the building $\Omega$, which is as a set the union of all thick flag spaces, identifies with the $\F$-valued points $\PCAF$ itself.

Now the face monoid $\GD$ acts on every module $L(\La)$, $\La\in P^+$. Dually we get actions of the opposite monoid $\widehat{G}^{\,op}$ on the restricted duals $L(\La)^{(*)}$, $\La\in P^+$, which fit together to an action of $\widehat{G}^{\,op}$ on the Cartan algebra $CA$ by morphisms of graded algebras. This in turn induces actions of $\GD$ on $\PCA$ and $\PCAF$ by morphisms. The image of $\Omega$ in $\PCAF$ is invariant under the action of $\GD$. By pulling back the action of the face monoid $\GD$ to the building $\Omega$ we obtain good action 1 of \cite{M5}, Corollary 50.\\

This algebraic geometric model of an action of the face monoid $\GD$ on the building $\Omega$ has been found as follows. V. Kac reviews in Section 2.6 of \cite{K2} the construction of the flag space $G/B$ from \cite{KP1}. It uses the Kostant cone ${\mathcal V}_\La:=G (L(\La)_\La)$, the $G$-orbit of the highest weight space $L(\La)_\La$ of $L(\La)$, where $\La$ is chosen arbitrarily from 
the set of strongly dominant weights $P^{++}$. V. Kac suggests on page 199, line 3 and 4, to use $\PCA$ for a definition of the flag space $G/B$, independent of the choice of $\La\in P^{++}$. 
Now the author had obtained different algebraic actions of the face monoid $\GD$ on the building $\Omega$, and was looking for everything linked to algebraic geometric actions of the face monoid on flag spaces. As shown in \cite{M1}, Proposition 5.4, the face monoid acts naturally on the Kostant cones ${\mathcal V}_\La$, $\La\in P^+$. These actions are related to an action of the face monoid $\GD$ on the building $\Omega$, which does not preserve its natural order, the bad action of \cite{M5}. The author tried $\PCA$, which led to this article.\\

A Cartan algebra can be obtained for every normal reductive algebraic monoid. Its spectrum of homogeneous prime ideals and the action of the reductive algebraic monoid on this spectrum will be investigated in a subsequent article.
%
%
%
%
%
\section{Preliminaries\label{Preli}} In this section we collect some facts about Kac-Moody algebras, minimal Kac-Moody groups, its associated face monoids, and formal Kac-Moody groups. We assume, as in the whole article, that the corresponding generalized Cartan matrices are symmetrizable. In particular, we recall the representation theoretic constructions of these groups and monoids. This section is also a reference for the notation used in the article.

The facts on Kac-Moody algebras, which we state here, can be found in \cite{K}, \cite{Ku}, and \cite{MoPi}, most results of \cite{K} and \cite{Ku} also valid for a field of characteristic zero with the same proofs. The facts on minimal Kac-Moody groups can be found in \cite{KP1}, \cite{KP2}, \cite{KP3}, \cite{Re}, and \cite{MoPi}. The facts on formal Kac-Moody groups can be found in \cite{Ku}, \cite{Re}, and \cite{Sl}, most results of \cite{Sl} and \cite{Ku} also valid for a field of characteristic zero with the same proofs. The facts on face monoids can be found in \cite{M1} and \cite{M5}.

We denote by $\N=\Z^+$, $\Q^+$, resp. $\R^+$ the sets of strictly positive numbers of $\Z$, $\Q$, resp. $\R\,$, and the sets $\N_0=\Z^+_0$, $\Q^+_0$, $\R^+_0$ contain, in addition, the zero.
$\F$ is a field of characteristic 0 and $\F^\times$ its group of units.\vspace*{1ex}

{\bf The starting data:} The starting data for the construction of the Kac-Moody algebra, the minimal and formal Kac-Moody groups and its face monoids are the following:
\begin{itemize}
\item A symmetrizable generalized Cartan matrix $A=(a_{ij})_{i,j\in I}$ with finite index set $I:=\{1,\,2,\,\ldots,\, n\}$. We denote by $l$ the rank of $A$. 
\item A {\it simply connected minimal free root base} for $A$ consisting of:
\begin{itemize} 
\item Dual free  $\mathbb{Z}$-modules $H,P$ of rank $2n-l$.
\item Linear independent sets $\{h_1,\,\ldots,\, h_n\}\subseteq H $, 
$\{\al_1,\,\ldots,\,\al_n\}\subseteq P$, such that $\al_i(h_j)=a_{ji}\,$, $i,j=1,\,\dots,\, n$. 
\item Furthermore, there exist $\La_1,\,\ldots,\,\La_n\in P$ such that $\La_i(h_j)=\delta_{ij}\,$, $i,j=1,\,\dots,\, n$.
\end{itemize}
\end{itemize}

$P$ is called the {\it weight lattice}, $P^+:=\Mklz{\La\in P}{\La(h_i)\geq 0 \mb{ for all } i\in I}$ the set of {\it dominant weights} of $P$. $Q:=\Z\mb{-span}\Mklz{\al_i}{i\in I}$ is called the {\it root lattice}. Set $Q^\pm_0:=\Z^\pm_0\mb{-span}\Mklz{\al_i\,}{\,i\in I}$ and $Q^\pm:=Q^\pm_0\setminus\{0\}$. The height function $\mb{ht}:Q\to \Z$ is defined by $\mb{ht}(\sum_{i\in I}n_i\al_i):=\sum_{i\in I}n_i$.

$\Mklz{\al_i}{i\in I}$ is called the set of {\it simple roots}. Note: To abbreviate many formulas we often identify this set with the index set $I$. 

We fix a system of elements $\La_1,\,\ldots,\,\La_n\in P$ such that $\La_i(h_j)=\delta_{ji}$, $i,j=1,\,\dots,\, n$. We extend to dual bases $h_1,\,\ldots,\,\ldots,\, h_{2n-l}\in H$ and $\La_1,\,\ldots,\, \La_{2n-l}\in P$. We call $\La_1,\,\ldots,\, \La_{2n-l}$ a system of {\it fundamental dominant weights} of $P$.

Let $J\subseteq I$. Set $A_J:=(a_{ij})_{i,j\in J}$, which is a symmetrizable generalized Cartan matrix if $J$ is nonempty. Set $P_J:=\Z\mb{-span}\Mklz{\La_i}{i\in J}$ and $P_J^+:= \Z_0^+\mb{-span}\Mklz{\La_i}{i\in J}$. 
Set $Q_J:=\Z\mb{-span}\Mklz{\al_i}{i\in J}$ and $(Q_J)^\pm_0:= \Z_0^\pm\mb{-span}\Mklz{\al_i}{i\in J}$ and $Q_J^\pm:=(Q_J)^\pm_0\setminus\{0\}$. Furthermore, set $P_{rest}:=\Z\mb{-span}\{\,\La_i \,| \,i = {n+1},\,\ldots,\,{2n-l}\,\}$. As always, a span of the empty set is defined to be $\{0\}$.\vspace*{1ex}

{\bf The linear spaces $\h$ and $\h^*$:} The lattices $H$ and $P$ are identified with the corresponding sublattices of the $\F$-linear spaces
\begin{eqnarray*}
   \h  :=  \h_\F   :=    H \otimes_{\mathbb Z} \F   \quad\mb{ and }\quad 
   \h^* = \h^*_\F  :=  P \otimes_{\mathbb Z} \F.
\end{eqnarray*}
We consider $\h^*$ as the dual of $\h$. We order the elements of $\h^*$ by $\la\leq\la'$ if and only 
if $\la'-\la\in Q_0^+$. 

For $J\subseteq I$ set $\h_J:=\mb{span}\Mklz{h_i}{i\in J}$. Set $\h_{rest}:=\mb{span}\Mklz{h_i}{i=n+1,\,\ldots,\,2n-l}$.
In the same way as in \S 2.1 of \cite{K} we equip $\h$ with a non-degenerate symmetric bilinear form $\iB{\;}{\;}$ adapted to the decomposition $\h=\h_I\oplus\h_{rest}$. It induces a linear isomorphism $\nu:\h\to \h^*$. It induces a non-degenerate symmetric bilinear form on $\h^*$, which is also denoted by $\iB{\;}{\;}$.\vspace*{1ex}

{\bf The Weyl group $\We$:} The {\it Weyl group} $\We=\We(A)$ is the Coxeter group with generators $\sigma_i$, $i\in I$, and
relations
\begin{eqnarray*}
       \sigma_i^2 \;=\; 1 \qquad (i\in I)    \;\:&\;,\;&\;\:
       {(\sigma_i\sigma_j)}^{m_{ij}} \;=\; 1 \qquad (i,j\in I,\,i\ne j),
\end{eqnarray*}
where the $m_{ij}$ are given by:
    $\quad \begin{tabular}{c|ccccc}
      $a_{ij}a_{ji}$ & 0 & 1 & 2 & 3 &  $\geq$ 4  \\[0.5ex] \hline 
         $m_{ij}$    & 2 & 3 & 4 & 6 &  no relation between $\sigma_i$ and $\sigma_j$ 
     \end{tabular}$\vspace*{1ex}\\

We denote by $l:\We\to\Nn$ the {\it length function}. We denote by $\We_J\cong\We(A_J)$, $J\subseteq I$, the {\it standard parabolic subgroups} of $\We$. For $J$, $K\subseteq I$ we denote by $\We^J$ the set of minimal coset representatives of $\We/\We_J$, by $\mb{}^K\We$ the set of minimal coset representatives of $\We_K\backslash \We$, and by $\mb{}^K\We^J$ the set of minimal double coset representatives of $\We_K\backslash \We/\We_J$.\vspace*{1ex}

{\bf The Tits cone $X$:} The Weyl group $\We$ acts faithfully $\h^*$ by 
\begin{eqnarray*}
  \sigma_i \la   :=  \la - \la(h_i) \al_i & \qquad i\in I,\quad \la\in \h^*.
\end{eqnarray*}
The lattices $Q$ and $P$ are left invariant by this action.

The {\it Tits cone} $X$ is a convex $\We$-invariant cone in $\h^*_\R$ obtained by
\begin{eqnarray*}
    X:=\We \overline{C} \quad\mb{ with }\quad \overline{C}:=\Mklz{\la\in\h^*_\R}{\la(h_i)\geq 0 \mb{ for all } i\in I }. 
\end{eqnarray*}
Every $\We$-orbit of $X$ contains exactly one point of the (closed) {\it standard fundamental chamber} $\overline{C}$.
Now we recall shortly:
\begin{itemize}
\item[1.] A $\We$-invariant partition of $X$ given by the interiors of polyhedral cones, which we call the facets of $X$.
\item[2.] The set $\RkX$ of faces of the convex cone $X$ in the sense of convex geometry. It is a $\We$-space, partially ordered by the inclusion of sets, even a complete lattice.
\end{itemize}

To 1. For $J\subseteq I$ set
\begin{eqnarray*}
          F_J    &:=&  \Mklz{\la\in\h^*_\R}{\la(h_i)=0 \;\mb{ for }\; i\in J,\;\; \la(h_i) > 0 \;\mb{ for }\;i\in I\setminus J},\\
  \overline{F_J} &:=&  \Mklz{\la\in\h^*_\R}{\la(h_i)=0\; \mb{ for }\;i\in J ,\;\;\la(h_i)\geq 0\; \mb{ for }\;i\in I\setminus J}=\dot{\bigcup_{K\supseteq J}}F_K. 
\end{eqnarray*}  
In particular, $\overline{C} = \overline{F_\emptyset}$. Here $\overline{F_J}$ is a polyhedral cone with relative interior $F_J$. Its faces in the sense of convex geometry are $\overline{F_K}$, $K\supseteq J$. We call $F_J$ resp. $\overline{F_J}$ an open resp. closed {\it standard facet} of type $J$.

The parabolic subgroup $\We_J$ of $\We$ is the stabilizer of any element $\la\in F_J$. It is the stabilizer of $F_J$ as well as of $\overline{F_J}$ as a whole. 

The set $\Mklz{\sigma F_J}{\sigma\in\We\,,\,J\subseteq I}$ is a $\We$-invariant partition of $X$. Partially ordered by the reverse closure relation it gives the Coxeter complex of $\We$. We call $\sigma F_J$ resp. $\sigma \overline{F_J}$ an open resp. closed {\it facet} of type $J$.\vspace*{1ex}

To 2. For $\emptyset\neq\Th \subseteq I$ we denote by $\Th^0$, resp. $\Th^\infty$ the set of indices which correspond to the sum of the components of the generalized Cartan submatrix $A_\Th$ of finite, resp. non-finite type. We set $\emptyset^0:=\emptyset^\infty:=\emptyset$. 
For $J\subseteq I$ set $J^\bot:=\Mklz{i\in I}{a_{ij}=0\;\mb{ for all }j\in J}$. 

A set $\Th\subseteq I$ such that $\Th=\Th^\infty$ is called {\it special}. If $\Th$ is a special set then 
\begin{eqnarray*}
   R(\Th):=\We_{\Th^\bot}\overline{F_\Th} = \Mklz{\la\in X}{\la(h_i)=0\;\mb{ for all }\;i\in\Th}
\end{eqnarray*}
is a face of $X$ with relative interior  
\begin{eqnarray*}
 ri \left( R(\Th) \right) = \We_{\Th^\bot} \bigcup_{J\subseteq \Th^{\bot}\,,\, J=J^0} F_{\Th\cup J}.
\end{eqnarray*}
The $\We$-stabilizers of $R(\Th)$, pointwise and as a whole, are
\begin{eqnarray*}
      Z_{\mathcal W}(R(\Th))  = \We_{\Th} \quad\mb{ and }\quad N_{\mathcal W}(R(\Th)) = \We_{\Th\cup \Th^\bot}.
\end{eqnarray*}

It is $\Mklz{R(\Th)}{\Th\;\textrm{special}\;}=\Mklz{R\in \RkX}{ri(R)\cap\overline{C}\neq\emptyset}$. Furthermore, 
\begin{eqnarray*}
  \RkX=\dot{\bigcup_{\Th\; \textrm{\footnotesize special}}}\Mklz{\sigma R(\Th)}{\sigma\in\We}.
\end{eqnarray*}
The special set $\Th$ is called the {\it type} of the face $\sigma R(\Th)$, $\sigma\in\We$.

Let $\sigma,\sigma'\in\We$, $\Th,\Th'$ be special. Then
\begin{eqnarray*}
   \sigma'R(\Th') \subseteq \sigma R(\Th)\quad \iff \quad \Th' \supseteq \Th \mb{ and } \sigma^{-1}\sigma' \in \We_{\Th^\bot}\We_{\Th'}.  
\end{eqnarray*}
Different faces of $X$ of the same type are not comparable by $\subseteq$.

A well defined function $red:\We\to I$ is obtained as follows: Set $\red{1}:=\emptyset$. If $\sigma\in\We\setminus\{1\}$ and $\sigma_{i_1}\cdots \sigma_{i_k}$ is a reduced expression for $\sigma$ set $\red{\sigma}:= \{i_1,\,\ldots,\,i_k\}\subseteq I$.

The lattice intersection, which coincides with the set theoretical intersection, and the lattice join of two arbitrary faces can be reduced easily by the formulas for the stabilizers to the following formulas: Let $\tau\in\mb{}^{\Th_1\cup\Th_1^\bot}\We^{\Th_2\cup\Th_2^\bot}$, $\Th_1,\Th_2$ be special. Then
\begin{eqnarray*}
   R(\Th_1)\cap \tau R(\Th_2)   &=&  R(\Th_1\cup\Th_2\cup \red{\tau}), \\
   R(\Th_1)\vee\tau R(\Th_2)  & =&  R((\Th_1\cap \tau\Th_2)^\infty). 
\end{eqnarray*}
%
%

%
%
%
%
%

{\bf The Kac-Moody algebra $\g$:} The {\it Kac-Moody algebra} $\g=\g(A)$ associated to the symmetrizable generalized Cartan matrix $A$ over $\F$ is the Lie algebra over $\F$ generated by the abelian Lie algebra $\h$, and
the elements $e_i$, $f_i$, ($i\in I$), with the following relations, which hold for any $i,j \in I$, $h \in \h$: 
  \begin{eqnarray*}
    \left[ e_i,f_j \right] \,=\,   \delta_{ij} h_i   \;\;,\;\; 
    \left[ h,e_i \right]   \,=\,   \al_i(h) e_i    \;\;,\;\;  
    \left[ h,f_i \right]   \,=\,  -\al_i(h) f_i\;\;,     \\
    \left(ad\,e_i\right)^{1-a_{ij}}e_j  \,=\, \left(ad\,f_i\right)^{1-a_{ij}}f_j  \,=0  \qquad (i\neq j)\;\;.\qquad
  \end{eqnarray*}
The abelian Lie algebra $\h$ and the elements $e_i$, $f_i$, ($i\in I$), identify with their images in $\g$.

There is the {\it root space decomposition}
\begin{eqnarray*}
  \g=\bigoplus_{\al \in \h^*}\g_{\al} \quad \mb{where} \quad \g_\al := \Mklz{x\in \g}{[h,x]=\al(h)\,x\;\mb{ for all }\;h\in \h}.            
\end{eqnarray*}
In particular, $\g_0 = \h$, $\g_{\al_i}=\F e_i$, and $\g_{-\al_i}=\F f_i$, $i\in I$. The set of {\it roots} $\W:=\Mklz{\al\in\h^*\setminus\{0\}}{\g_\al\ne 0}$ 
satisfies $\W\subseteq Q$, $\W=-\W$, and it is invariant under the Weyl group. $\Delta_{re}:=\We\Mklz{\al_i}{i\in I}$ is the set of {\it real roots}. $\iW:=\W\setminus\rW$ is the set of {\it imaginary roots}. Every root space $\g_\al$, $\al\in\Delta$, is finite dimensional.

$\W$, $\rW$, and $\iW$ decompose into the disjoint union of the sets of {\it positive} and {\it negative roots} $\W^\pm:=\W\cap Q^\pm$, 
$\rW^\pm:=\rW\cap Q^\pm$, $\iW^\pm:=\iW\cap Q^\pm$. There is the corresponding {\it triangular decomposition}
\begin{eqnarray*}
    \g=\nm \oplus \h \oplus \np\quad \mb{ where }\quad \n^\pm :=\bigoplus_{\al\in {\Delta}^\pm} \g_\al\,.
\end{eqnarray*} 
For every $\al\in \rW$ the subalgebra $\g_\al\oplus[\g_\al,\g_{-\al}]\oplus\g_{-\al}$ of $\g$ is isomorphic to $sl(2,\F)$.

The non-degenerate symmetric bilinear form ( $|$ ) on $\h$ extends uniquely to a non-degenerate symmetric invariant bilinear form ( $|$ ) on $\g$. For $\al,\,\beta\in\Delta\cup\{0\}$ such that $\al+\beta\neq 0$ it is $(\g_\al|\g_\beta)=0$.
The restriction $(\; | \;):\g_\al\times\g_{-\al}\to\F$ is nondegenerate and $[x,y]=(x|y)\nu^{-1}(\al)$ for all $x\in\g_\al$, $y\in\g_{-\al}$, $\al\in\Delta$.

Let $J\subseteq I$. Set 
\begin{eqnarray*}
   \g_J:=\h_J\oplus\bigoplus_{\al\in\Delta_J} \g_\al=\n_J^-\oplus\h_J\oplus \n_J^+ \quad\mb{ where }\quad \n_J^\pm := \bigoplus_{\al\in \Delta_J^\pm}\g_\al 
\end{eqnarray*}
and $\Delta_J:=\Delta\cap Q_J$,  $\Delta_J^\pm:=\Delta^\pm\cap Q_J$. Then $\g_J\cong \g(A_J)'$ for $J\neq \emptyset$. In particular, $\g_I$ is the derived Lie algebra of $\g$.

An ideal of $\n^\pm$ is obtained by
\begin{eqnarray*}
  (\n^\pm)^J := \bigoplus_{\al\in(\Delta^\pm)^J}\g_\al \quad\mb{ with }\quad (\Delta^\pm)^J:=\Delta^\pm\setminus\Delta_J^\pm.
\end{eqnarray*}

{\bf The (minimal) Kac-Moody group $G$:}
A $\g$-module $V$ is called {\it integrable} if $V$ is $\h$-diagonalizable with set of weights $P(V)\subseteq P$ and the elements of $\g_\al$ act locally nilpotent on $V$ for all $\al\in\rW$.
Examples of integrable representations are the adjoint representation $\g$, the irreducible highest 
weight representations $L(\La)$ with highest weights $\La\in P^+$, and the irreducible lowest weight representations $L(\La)^{(*)}$ with lowest weights $-\La\in -P^+$.

Let $\mathcal I$ be the category whose objects are the integrable $\g$-modules and whose morphisms are the morphisms of $\g$-modules. This category generalizes the category of locally finite dimensional representations of a semisimple Lie algebra.

Let $Nat({\mathcal I})$ be the set of natural transformations of the forgetful functor from the category $\mathcal I$ into the category of $\F$-linear spaces. Explicitly, $Nat({\mathcal I})$ consists of the families of linear maps
\begin{eqnarray*}
   m &=& \left(\,m_V\in End(V)\,\right)_{\,V \:\textrm{\footnotesize obj.}\;\textrm{\footnotesize of } {\mathcal I}}\;,
\end{eqnarray*} 
such that for all objects $V$, $W$, and all morphisms $\phi:V\to W$ of $\mathcal I$ the diagram 
\begin{eqnarray*}
\xymatrix{ V \ar[d]_{\mbox{\footnotesize  $\phi$}}  \ar[r]^{\mb{\footnotesize $m_V$} }  & V \ar[d]^{\mbox{\footnotesize  $\phi$}} \\
               W \ar[r]^{\mb{\footnotesize  $m_W$}}  & W}
\end{eqnarray*}
commutes. $Nat({\mathcal I})$ gets in the obvious way the structure of an associative $\F$-algebra with unit. In particular, there are the following elements of $Nat({\mathcal I})$:

(1) For every $h\in H$, $s\in\F^\times$ there exists $t_h(s)\in Nat({\mathcal I})$, such that for every object $V$ of $\mathcal I$
it holds
\begin{eqnarray*}
  t_h(s)v_\la = s^{\la(h)}v_\la \quad,\quad v_\la\in V_\la\;\,,\;\,\la\in P(V).
\end{eqnarray*}

(2) For every $x\in\g_\al$, $\al\in\rW$, there exists  $exp(x)\in Nat({\mathcal I})$, such that for 
every object $V$ of $\mathcal I$ it holds
\begin{eqnarray*}
  exp(x)v = \sum_{k\in\Nn} \frac{x^k v}{k!}\quad,\quad v\in V.
\end{eqnarray*}

The {\it minimal Kac-Moody group} $G$, which we call {\it Kac-Moody group} for short, is the submonoid of $Nat({\mathcal I})$ generated by the elements of (1) and (2). Its elements are compatible with $\oplus$, $\otimes$, act as identity on the trivial representation, and induce dual maps on the integrable duals, i.e., on the biggest integrable submodules of the full duals. The Kac-Moody group $G$ has the following important structural properties:

A root group data system $(T,\, (U_\al)_{\al\in \Delta_{re}})$ of $G$ is obtained as follows: The elements of (1) induce an embedding of the torus $H\otimes_\Z\F^\times = Hom((P,+),\,(\F^\times,\cdot))$ into $G$, whose image is $T$. For every $\al\in\rW$ the elements of (2) induce an embedding of $(\g_\al,+)$ into $G$, whose image is the {\it root group} $U_\al$. 

The normalizer $N:=N_G(T)$ is generated by $T$ and the elements $n_i:=exp(e_i)exp(-f_i)exp(e_i)$, $i\in I$. An isomorphism $\kappa:\,\We\to N/T$ is induced by $\kappa(\sigma_i):=n_i T$, $i\in I$.

We denote by $\widetilde{\quad}:\We\to N$ the cross section of the canonical map from $N$ to $\We$ defined by $\widetilde{1}=1$, $\widetilde{\sigma_i}=n_i$, $i\in I$, and $\widetilde{\sigma\tau}=\widetilde{\sigma}\widetilde{\tau}$ if  $l(\sigma\tau)=l(\sigma)+l(\tau)$, $\sigma,\tau\in\We$.

We denote an arbitrary element $n\in N$ with $\kappa^{-1}(nT)=\sigma\in\We$ by $n_\sigma$. The set of weights $P(V)$ of an
integrable $\g$-module $V$ is $\We$-invariant and $n_\sigma V_\la =V_{\sigma \la}$, $\la\in P(V)$, $\sigma\in\We$.

Let $U^\pm$ be the subgroups of $G$ generated by $U_\al$, $\al\in\Delta_{re}^\pm $. Then $U^\pm$ are normalized by $T$. The pairs ($B^\pm :=  T\ltimes U^\pm $, $N$) are twin BN-pairs of $G$ with the property $B^+\cap B^-= B^\pm\cap N = T$. There are the {\it Bruhat} and {\it Birkhoff decompositions}
\begin{eqnarray*}
  G  =  \dot{ \bigcup_{\sigma\in N }} U^\epsilon n U^\delta = \dot{ \bigcup_{\sigma\in {\mathcal W} }} B^\epsilon\sigma B^\delta \qquad,\qquad \epsilon,\delta\;\in\;\{\,+\,,\,-\,\}.
\end{eqnarray*}

There are also {\it Levi decompositions} of the standard parabolic subgroups $P_J^\pm$, $J\subseteq I$, and of $U^\pm$:
\begin{eqnarray*}
   P_J^\pm = L^\pm_J \ltimes (U^\pm)^J \quad \mb{ and }\quad U^\pm =  U^\pm_J \ltimes (U^\pm)^J  .
\end{eqnarray*}
Here $U_J^\pm$ is the group generated by $U_\al$, $\al\in (\W_J^\pm)_{re}:=\W_J^\pm\cap \Delta_{re}$. $(U^\pm)^J$ is the smallest normal subgroup of $U^\pm$ containing $U_\al$, $\al\in (\W^\pm)^J_{re}:=(\W^\pm)^J\cap \Delta_{re}$. This group equals $\bigcap_{\sigma\in{\We}_J}\sigma U^\pm\sigma^{-1}$. Furthermore, $L_J^\pm:=U_J^\pm (T\We_J) U_J^\pm$.

For $J\subseteq I$ let $G_J$ be the subgroup of $G$ generated by $U_\al$, $\al\in (\Delta_J)_{re}:=\Delta_J\cap\Delta_{re}$.
Then $G_\emptyset =\{1\}$ and $G_J\cong G(A_J)'$ for $J\neq\emptyset$. The torus $T_J:=G_J\cap T$ is generated by the subtori $t_{h_i}(\F^\times)$, $i\in J$. The group $N_J:=G_J\cap N$ is generated by $T_J$ and the elements $n_i$, $i\in J$. 

In particular, $G_I$ is the derived group of $G$, which identifies with the Kac-Moody group as defined in \cite{KP1}. It is $G = G_I\rtimes T_{rest}$, where $T_{rest}$ is the subtorus of $T$ generated by the subtori $t_{h_i}(\F^\times)$, $i=n+1,\,\ldots,\,2n-l$.

The adjoint action $ad$ of the Kac-Moody algebra $\g$ on itself is integrable. Therefore, the Kac-Moody group $G$ acts on $\g$. This action is called the {\it adjoint action} $Ad$ of $G$.\vspace*{1ex}

{\bf Some stabilizers:} Let $J\subseteq I$ and $\La\in P^+\cap F_J$. Then
\begin{eqnarray*}
  N_{\bf g}(L(\La)_\La):=\Mklz{x\in \g}{x L(\La)_\La\subseteq L(\La)_\La} = \n_J^-\oplus \h\oplus\n\,.
\end{eqnarray*}
Furthermore, let $T^J$ be the subtorus of $T$ generated by $t_{h_i}(\F^\times)$, $i\in \{1,\,2,\,\ldots,\,2n-l\}\setminus J$. Set
\begin{eqnarray*}
   Y_\La := \Big\{\prod_{i=1,\,\ldots,\,2n-l,\,i\notin J} t_i(s_i)\in T \,\Big|\, \prod_{i=1,\,\ldots,\,2n-l,\,i\notin J}(s_i)^{\La(h_i)}=1 \Big\}\subseteq T^J\,.
\end{eqnarray*}
Then
\begin{eqnarray*}
   && N_G(L(\La)_\La):=\Mklz{g\in G}{g L(\La)_\La=L(\La)_\La} = P_J = U^J G_J T^J  = U N_J T^J U\,,\\
   && Z_G(L(\La)_\La):=\Mklz{g\in G}{g v=v \mb{ for all } v\in L(\La)_\La } = U^J G_J Y_\La  = U N_J Y_\La U\,.
\end{eqnarray*}

{\bf The face monoid $\GD$ associated to the Kac-Moody group $G$:} 
The category ${\mathcal O}^p$ is defined as follows: Its objects are the $\g$-modules $V$, which have the properties:
\begin{itemize}
\item[1.] $V$ is $\h$-diagonalizable with finite dimensional weight spaces.
\item[2.] There exist finitely many elements $\la_1,\,\ldots,\,\la_m\in\h^*$, such that the set of weights $P(V)$ of $V$ is 
contained in the union $\bigcup_{1=1}^m \Mklz{\la\in\h^*}{\la\leq \la_i}$.
\end{itemize}
The morphisms of ${\mathcal O}^p$ are the morphisms of $\g$-modules.

Let ${\mathcal O}_{int}^p$ the full subcategory of the category ${\mathcal O}^p$, whose objects are integrable modules. This category generalizes the category of finite dimensional representations of a semisimple Lie algebra, keeping the complete reducibility theorem: Every object of ${\mathcal O}_{int}^p$ is isomorphic to a direct sum of the integrable irreducible highest weight modules $L(\La)$, $\La\in P^+$.
Denote the set of weights of $L(\La)$ by $P(\La)$, $\La\in P^+$. The set of weights of every object of ${\mathcal O}_{int}^p$ is contained in $X\cap P$ because of $\bigcup_{\La\in P^+} P(\La) = X\cap P$.

Let $V$ be an object of ${\mathcal O}_{int}^p$ and $V=\bigoplus_{\la\in P(V)}V_\la$ its weight space decomposition. We may associate two $\F$-linear duals, the restricted dual $V^{(*)}$ and the full dual $V^*$ as follows 
\begin{eqnarray*}
   V^{(*)}:=\bigoplus_{\la\in P(V)}V_\la^*  \;&\subseteq& \; V^* =\prod_{\la\in P(V)}V_\la^*.
\end{eqnarray*}

Let $Nat({\mathcal O}_{int}^p)$ be the set of natural transformations of the forgetful functor from the category ${\mathcal O}_{int}^p$ into the category of $\F$-linear spaces. In particular, there are the following elements of the $\F$-algebra $Nat({\mathcal O}_{int}^p)$:

(1) For every $h\in H$, $s\in\F^\times$ there exists $t_h(s)\in Nat({\mathcal O}_{int}^p)$, such that for every object $V$ of ${\mathcal O}_{int}^p$
it holds
\begin{eqnarray*}
  t_h(s)v_\la = s^{\la(h)}v_\la \quad,\quad v_\la\in V_\la\;\,,\;\,\la\in P(V).
\end{eqnarray*}

(2) For every $x\in\g_\al$, $\al\in\rW$, there exists  $exp(x)\in Nat({\mathcal O}_{int}^p)$, such that for 
every object $V$ of ${\mathcal O}_{int}^p$ it holds
\begin{eqnarray*}
  exp(x)v = \sum_{k\in\Nn} \frac{x^k v}{k!}\quad,\quad v\in V.
\end{eqnarray*}

(3) For every $R\in \mb{Fa}(X)$ there exists $e(R)\in Nat({\mathcal O}_{int}^p)$, such that for every object $V$ of ${\mathcal O}_{int}^p$ it holds
\begin{eqnarray*}
  e(R)v_\la = \left\{ \begin{array}{ccl}
    v_\la &\mb{ if}&\la\in R\\
      0   &\mb{ if}& \la\in X\setminus R
  \end{array}\right. \quad,\quad \;v_\la\in V_\la\,,\; \;\la\in P(V).
\end{eqnarray*}

The {\it face monoid} $\GD$ is the submonoid of $Nat({\mathcal O}_{int}^p)$ generated by the elements of (1), (2), and (3). It is a Tannaka monoid: It is the biggest part of $Nat({\mathcal O}_{int}^p)$, whose elements are compatible with $\oplus$, $\otimes$, act as identity on the trivial representation, and induce dual maps on the restricted duals.

The unit group of $\GD$ is generated by the elements of (1) and (2). It is isomorphic to the Kac-Moody group $G$ in the obvious way, and for simplicity of notation we identify it with $G$.

Let $\Th$ be special. There are the following decompositions of certain parabolic subgroups of $G$: 
\begin{eqnarray*}
    P_{\Th\cup\Th^\bot}  &=& \left(G_{\Th^\bot}\rtimes T^{\Th\cup\Th^\bot}\right)\ltimes \left(G_\Th\ltimes U^{\Th\cup\Th^\bot}\right).\\
    P_{\Th\cup\Th^\bot}^- &=& \left((U^-)^{\Th\cup\Th^\bot}\rtimes G_\Th\right) \rtimes\left(G_{\Th^\bot}\rtimes T^{\Th\cup\Th^\bot}\right).
\end{eqnarray*}
The projections belonging to these semidirect products are denoted by $p_\Th^\pm :\, P_{\Th\cup\Th^\bot}^\pm  \to  G_{\Th^\bot}\rtimes T^{\Th\cup\Th^\bot}$. Now the monoid $\GD$ may be described algebraically as follows: It is
\begin{eqnarray*}
  \GD = \dot{\bigcup_{\Th\;   \textrm{\footnotesize special}}} G e(R(\Th)) G  .
\end{eqnarray*}
Let $g_1,g_2, h_1, h_2\in G$, let $\Th_1$, $\Th_2$ be special. It holds $g_1 e(R(\Th_1)) h_1= g_2 e(R(\Th_2)) h_2 $ if and only if $\Th_2 =\Th_1$ and 
\begin{eqnarray*} 
  (g_2)^{-1}g_1\in P_{\Th_1\cup\Th_1^\bot},\;\; h_2(h_1)^{-1}\in P_{\Th_1\cup\Th_1^\bot}^- \;\mb{ with }\; p_{\Th_1}((g_2)^{-1}g_1) = p_{\Th_1}^-(h_2 (h_1)^{-1}).
\end{eqnarray*}
To describe the multiplication of $g_1 e(R(\Th_1)) h_1$ and $g_2 e(R(\Th_2)) h_2 $ write $h_1 g_2= a_1 \widetilde{\sigma} a_2$ with $a_1\in P_{\Th_1\cup\Th_1^\bot}^-$, $a_2\in  P_{\Th_2\cup\Th_2^\bot}$, and $\sigma\in \mb{}^{\Th_1\cup\Th_1^\bot}\We^{\Th_2\cup\Th_2^\bot}$. Then
\begin{eqnarray*}
 g_1 e(R(\Th_1)) h_1 \, g_2 e(R(\Th_2)) h_2 = g_1 p_{\Th_1}^-(a_1)\, e(\,R(\Th_1\cup\Th_2\cup\red{\sigma})\,)\,p_{\Th_2}(a_2)h_2.
\end{eqnarray*}

The monoid $\ND:=N_{\widehat{G}}(T):=\{\,\hat{g}\in\GD\mid\hat{g}T=T\hat{g}\,\}$ is generated by $N$ and $\Mklz{e(R)}{R\in\RkX}$.
There are the Bruhat and Birkhoff decompositions 
\begin{eqnarray*}
   \GD  = \dot{\bigcup_{\widehat{n}\,\in\,\widehat{N}}} U^\eps\widehat{n}U^\delta 
        = \dot{\bigcup_{\widehat{w}\,\in\,\widehat{N}/T}} B^\eps\widehat{w}B^\delta
\qquad,\qquad \epsilon,\delta\;\in\;\{\,+\,,\,-\,\}.      
\end{eqnarray*}

The Weyl group $\We$ acts on the monoid (\,$\RkX\,,\,\cap$\,). The semidirect product $\We\ltimes \RkX $ consists of the set 
$\We\times \RkX$ on which the structure of a monoid is given by 
\begin{eqnarray*}
  (\sigma, R)\cdot (\tau, S) := (\sigma \tau, \tau^{-1} R\cap S) .
\end{eqnarray*} 
A congruence relation on $\We\ltimes \RkX $ is obtained by
\begin{eqnarray*}
   (\sigma, R) \sim (\sigma', R') &:\iff & R\:=\:R'\quad\mb{and}\quad \sigma^{-1}\sigma'\in Z_{\mathcal W}(R)\;.
\end{eqnarray*}
We denote the congruence class of $(\sigma, R)$ by $\sigma \ve{R}$.
Now the monoid $\ND/T$ is isomorphic to the monoid $\WeD:=(\We\ltimes \RkX)/\sim$, which we call the {\it face monoid} associated to $\We$, by 
\begin{eqnarray*}
\begin{array}{rccc}
  \kappa : &   \WeD    & \to     &   \ND/T       \\
           & \sigma \ve{R}        & \mapsto & \widetilde{\sigma} e(R)T
\end{array}.
\end{eqnarray*}

The Weyl group $\We$ is the unit group of $\WeD$. There are standard parabolic submonoids $\WeD_J\cong\WeD(A_J)$, $J\subseteq I$, of $\WeD$. These give the standard parabolic submonoids $\PD_J := B \,\WeD_J B$, $J\subseteq I$, of $\GD$. \vspace*{1ex}

{\bf The formal Kac-Moody group $G_{fn}$:} The Lie bracket of $\g$ extends in the obvious way to a Lie bracket of 
\begin{eqnarray*}
   \g_{fn}:=\n^-_f\oplus\h\oplus\n   \quad\mb{ with }\quad \n_f^-:=\prod_{\al\in \Delta^-}\g_\al\,.
\end{eqnarray*} 
This Lie algebra is called the {\it formal Kac-Moody algebra}. The Lie subalgebra $\n_f^-$ is a pronilpotent Lie algebra with (non-complete) defining system of ideals $\n_f^-(k):=\prod_{\al\in\Delta^-,\, -ht(\al)\geq k+1}\g_\al$, $k\in\Nn$. 

Define the category ${\mathcal O}^n$ similarly as the category ${\mathcal O}^p$ , but replace $\leq$ by $\geq$ in the second condition.
Let ${\mathcal O}_{int}^n$ the full subcategory of the category ${\mathcal O}^n$, whose objects are integrable modules. This category generalizes the category of finite dimensional representations of a semisimple Lie algebra, keeping the complete reducibility theorem: Every object of ${\mathcal O}_{int}^n$ is isomorphic to a direct sum of the integrable irreducible lowest weight modules $L(\La)^{(*)}$ of lowest weights $-\La\in -P^+$.

For every object $V$ of ${\mathcal O}_{int}^n$ the action of the Kac-Moody algebra $\g$ on $V$ extends in the obvious way to an action of the formal Kac-Moody algebra $\g_{fn}$ on $V$. 

Let $Nat({\mathcal O}_{int}^n)$ be the set of natural transformations of the forgetful functor from the category ${\mathcal O}_{int}^n$ into the category of $\F$-linear spaces. In particular, there are the following elements of the $\F$-algebra $Nat({\mathcal O}_{int}^n)$:

(1) For every $h\in H$, $s\in\F^\times$ there exists $t_h(s)\in Nat({\mathcal O}_{int}^n)$, such that for every object $V$ of ${\mathcal O}_{int}^n$
it holds
\begin{eqnarray*}
  t_h(s)v_\la = s^{\la(h)}v_\la \quad,\quad v_\la\in V_\la\;\,,\;\,\la\in P(V).
\end{eqnarray*}

(2) For every $x\in \n_f^-\cup\bigcup_{\al\in\Delta^+_{re}}\g_\al$ there exists $exp(x)\in Nat({\mathcal O}_{int}^n)$, such that for every object $V$ of ${\mathcal O}_{int}^n$ it holds
\begin{eqnarray*}
  exp(x)v = \sum_{k\in\Nn} \frac{x^k v}{k!}\quad,\quad v\in V.
\end{eqnarray*}

The {\it formal Kac-Moody group} $G_{fn}$ is the submonoid of $Nat({\mathcal O}_{int}^n)$ generated by the elements of (1) and (2). Its elements are compatible with $\oplus$, $\otimes$, act as identity on the trivial representation, and induce dual maps now on the full duals. The formal Kac-Moody group $G_{fn}$  has the following important structural properties:

$U_f^-:= exp(\n_f^-)$ is a subgroup of $G_{fn}$. As a group $U_f^-$ is the prounipotent group with pronilpotent Lie algebra $\n_f^-$ and exponential map $exp:\n_f^-\to U_f^-$. 

Now certain subsets of $\Delta^-$ induce prounipotent subgroups and prounipotent normal subgroups of $U_f^-$. Certain decompositions induce multiplicative decompositions of $U_f^-$ and of prounipotent subgoups of $U_f^-$. In this article we use many times: For $J\subseteq I$ the decomposition $\Delta^-= \Delta^-_J\,\dot{\cup}\,(\Delta^-)^J$ induces the semidirect decomposition
\begin{eqnarray*}
  \n^-_f = (\n_f^-)_J\ltimes (\n_f^-)^J \quad \mb{ with }\quad 
(\n_f^-)_J:=\prod_{\al\in\Delta_J^-}\g_\al\; \mb{ and }\; (\n_f^-)^J:=\prod_{\al\in(\Delta^-)^J}\g_\al \,,
\end{eqnarray*}
which in turn induces the semidirect decomposition 
\begin{eqnarray*}
   U_f^-= (U_f^-)_J \ltimes (U_f^-)^J \quad\mb{ with }\quad (U_f^-)_J:= exp((\n_f^-)_J)  \;\mb{ and }\; (U_f^-)^J:= exp((\n_f^-)^J)\,.
\end{eqnarray*}

The formal Kac-Moody group $G_{fn}$ contains in the obvious way a subgroup isomorphic the Kac-Moody group $G$, and for simplicity of notation we identify it with $G$. 

It is $G_{fn}=U_f^- G=G U_f^-$. Furthermore, $G_{fn}=G$ if and only if $U_f^-=U^-$ if and only if the generalized Cartan matrix $A$ contains only components of finite type.

$U_f^-$ is normalized by $T$ and ($B_f^-:=T U_f^-$, $N$) is a BN-pair of $G_{fn}$ with $B_f^-\cap N=T$. More generally, $G_{fn}$, $U_f^-$, $U$, $N$, $T$, $\Mklz{\sigma_i\in \We\cong N/T}{i\in I}$ is a refined Tits system. (To adapt the notation used for refined Tits systems in \cite{KP3}, Section 2, or in \cite{Ku}, Section 5.2 to our situation, the groups $U^+$ and $U^-$ in the axioms have to be interchanged.) In particular, there are the {\it Bruhat} and {\it Birkhoff decompositions}
\begin{eqnarray*}
  G_{fn} =  \dot{ \bigcup_{\sigma\in {\mathcal W} }}  B_f^- \sigma B_f^- =  \dot{ \bigcup_{\sigma\in {\mathcal W} }} B^\epsilon\sigma B_f^- =  \dot{ \bigcup_{\sigma\in {\mathcal W} }}  B_f^- \sigma B^\epsilon \qquad,\qquad \epsilon\;\in\;\{\,+\,,\,-\,\}.
\end{eqnarray*}

There are {\it Levi decompositions} of the standard parabolic subgroups 
\begin{eqnarray*}
    (P_{fn}^-)_J = (L_{fn})_J \ltimes (U_f^-)^J 
\end{eqnarray*}
with $(L_{fn})_J=(U_f^-)_J(\We_J T)(U_f^-)_J=(U_f^-)_J L_J=L_J (U_f^-)_J$, $J\subseteq I$.

For $J\subseteq I$ let $(G_{fn})_J$ be the subgroup of $G_{fn}$ generated by  $(U_f^-)_J$ and $U_\al$, $\al\in (\Delta^+_J)_{re}$. Then $(G_{fn})_\emptyset=\{1\}$ and $(G_{fn})_J\cong G_{fn}(A_J)'$ for $J\neq\emptyset$. 

The formal Kac-Moody group $G_{fn}$ acts on the formal Kac-Moody algebra $\g_{fn}$ by the {\it adjoint action} $Ad$. On the generators (1) and (2) of $G_{fn}$ it is given as follows. For $h\in H$ and $s\in\F^\times$ it is
\begin{eqnarray*}
   Ad(t_h(s)) y = \sum_{\al\in \Delta\cup\{0\}}s^{\al(h)} y_\al \quad \mb{ where }\;y=\sum_{\al\in\Delta\cup\{0\}}y_\al\in\g_{fn}.
\end{eqnarray*} 
For $x\in\n_f^-\cup\bigcup_{\al\in\Delta^+_{re}}\g_\al$ it is 
\begin{eqnarray*}
 Ad(exp(x)) y = exp (ad(x)) y = \sum_{k\in\Nn} \frac{1}{k!}(ad(x))^k y \quad \mb{ where }\;y\in\g_{fn}.
\end{eqnarray*} 
This action extends the adjoint action $Ad$ of $G$ on $\g$.

It is important to note the following: If $V$ is a $\g$-module contained in ${\mathcal O}_{int}^p$, then the restricted dual $\g$-module $V^{(*)}$ is contained in ${\mathcal O}_{int}^n$. The Kac-Moody group $G$ acts on $V$ and $V^{(*)}$ dually. The formal Kac-Moody group $G_{fn}$ only acts on  
\begin{eqnarray*}
   V_f:= \prod_{\la\in P(V)} V_\la = \prod_{\la\in P(V)} (V_\la^*)^* = (V^{(*)})^*  
\end{eqnarray*}
and $V^{(*)}$ dually, where $V_\la$ has been identified canonically with $(V_\la^*)^*$, $\la\in P(V)$. The $G_{fn}$-module $V_f$ extends the $G$-module $V$. The linear subspace $V$ of $V_f$ is $G$-invariant, but in general not $G_{fn}$-invariant.

Note also: In the literature it is more common to consider the formal Kac-Moody algebra
 \begin{eqnarray*}
 \g_{fp}:=\n^-\oplus\h\oplus\n_f \quad \mb{ with }\quad \n_f:=\prod_{\al\in \Delta^+}\g_\al \,.
\end{eqnarray*} 
For every object $V$ of ${\mathcal O}_{int}^p$ the action of the Kac-Moody algebra $\g$ on $V$ extends to an action of the formal Kac-Moody algebra $\g_{fp}$ on $V$. Now the {\it formal Kac-Moody group} $G_{fp}$ is the submonoid of $Nat({\mathcal O}_{int}^p)$ generated by the elements given in (1) and (2):

(1) For every $h\in H$, $s\in\F^\times$ there exists $t_h(s)\in Nat({\mathcal O}_{int}^p)$, such that for every object $V$ of ${\mathcal O}_{int}^p$
it holds
\begin{eqnarray*}
  t_h(s)v_\la = s^{\la(h)}v_\la \quad,\quad v_\la\in V_\la\;\,,\;\,\la\in P(V).
\end{eqnarray*}

(2) For every $x\in \n_f \cup\bigcup_{\al\in\Delta^-_{re}}\g_\al$ there exists $exp(x)\in Nat({\mathcal O}_{int}^p)$, such that for every object $V$ of ${\mathcal O}_{int}^p$ it holds
\begin{eqnarray*}
  exp(x)v = \sum_{k\in\Nn} \frac{x^k v}{k!}\quad,\quad v\in V.
\end{eqnarray*}
The formal Kac-Moody group $G_{fp}$ now contains the prounipotent group $U_f=exp(\n_f)$.

The group $G_{fp}$ identifies with the Kac-Moody group of Kapitel 5 in \cite{Sl} for a simply connected minimal free realization. It identifies with the Kac-Moody group of Section 6 in \cite{Ku}. In \cite{Sl}, \cite{Ku} the field is assumed to be $\C$, but it works in the same way over a field $\F$ of characteristic 0. The results for $G_{fn}$ stated above can be proved similarly as the corresponding results for $G_{fp}$.\vspace*{1ex}

{\bf The coordinate ring of the face monoid $\GD$ and its restriction to the (minimal) Kac-Moody group $G$:} The Kac Moody group $G$ acts by its definition on every $\g$-module $V$ contained in ${\mathcal O}^p_{int}$ and on its restricted dual $V^{(*)}$. We now change this action of $G$  on $V^{(*)}$ to an action of $G^{op}$ on $V^{(*)}$ by concatenation the corresponding representation of $G$ with the inverse map. For a $\g$-module $V$ contained in ${\mathcal O}^p_{int}$, for $v\in V$ and $\phi\in V^{(*)}$ the map
\begin{eqnarray*}
\phi(\ins v):\, G  & \to     & \;\F         \\
               g  & \mapsto & \phi(gv)
\end{eqnarray*}
is called the matrix coefficient of $\phi$ and $v$ on $G$.

The coordinate ring $\FK{G}$ of the Kac-Moody group $G$ is the set of matrix coefficients obtained in this way by the modules contained in ${\mathcal O}^p_{int}$ and their restricted duals. It is a $\F$-algebra on which $G^{op}\times G$ acts. 

As explained above, $G = G_I\rtimes T_{rest}$. Here, the derived group $G_I$ of $G$ identifies with the Kac-Moody group of \cite{KP1} and $T_{rest}$ is a certain subtorus of $T$. The multiplication map $m: G_I\times T_{rest}\to G$ induces an isomorphism $m^*:\FK{G}\to\FK{G_I}\otimes \FK{T_{rest}}$. The restriction $\FK{G_I}:=\FK{G}\res{G_I}$ coincides with the algebra of strongly regular functions of \cite{KP2}. The restriction $\FK{T_{rest}}:=\FK{G}\res{T_{rest}}$ is the classical coordinate ring of the torus $T_{rest}$, which identifies with the group algebra $\FK{P_{rest}}$ of the lattice $P_{rest}$. 

V. Kac and D. Peterson showed in Theorem 1 of \cite{KP2} a Peter-Weyl theorem for $\FK{G_I}$ . They showed in Corollary 2.2 and Remark 2.1 of \cite{KP2} a Borel-Weil theorem for $\F[G_I]^U$. It is easy to transcribe these theorems to theorems for $\FK{G}$ and $\F[G]^U$. (Use in addition: $P_I$ identifies with the weight lattice of \cite{KP2}. It is $P=P_I\oplus P_{rest}$ and $P^+=P_I^+\oplus P_{rest}$. The modules $L(N)$, $N\in P_{rest}$, are one-dimensional. In particular, $L(\La)\otimes L(N)=L(\La+N)$ for all $\La\in P^+$ and $N\in P_{rest}$.) Also the proofs of these theorems could be easily adapted to obtain the following theorems.

The Peter-Weyl theorem: A $G^{op}\times G$-equivariant bijective linear map
\begin{eqnarray*}
    PW:\;\bigoplus_{\La\in P^+}L(\La)^{(*)}\otimes L(\La) \to \FK{G}   
\end{eqnarray*}
is given by $PW(\phi\otimes v):=\phi(\ins v)$ for all $\phi\in L(\La)^{(*)}$, $v\in L(\La)$, $\La\in P^+$. 

The Borel-Weil theorem: Denote by $\F[G]^U$ the functions of $\FK{G}$ which take constant values on the elements of every coset $gU$, $g\in G$. Choose a highest weight vector $v_\La\in L(\La)_\La$ for every $\La\in P^+$. A $G^{op}$-equivariant linear bijective map
\begin{eqnarray*}
  BW:\; \bigoplus_{\La\in P^+}L(\La)^{(*)}\to \F[G]^U
\end{eqnarray*}
is given by $BW(\phi):=\phi(\ins v_\La)$ for all $\phi\in L(\La)^{(*)}$, $\La\in P^+$. Pulling back the algebra structure of $\F[G]^U$ by $BW$ gives the Cartan algebra structure on $\bigoplus_{\La\in P^+}L(\La)^{(*)}$, which will be explained in subsection \ref{subsectionCA}.

Note also that the coordinate ring $\FK{G}$, and therefore also $\F[G]^U$, is an integral domain. This follows with Theorem 2.1 of \cite{KP2}, or could be shown directly.

Similarly as above we may define matrix coefficients for the face monoid $\GD$, and the algebra of matrix coefficients $\FK{\GD}$, on which $\widehat{G}^{\,op}\times\GD$ acts. By definition, $\FK{G}=\FK{\GD}\res{G}$. It has been shown in \cite{M1} that the Kac-Moody group $G$ is the Zariski open dense unit group of the face monoid $\GD$. In particular, the coordinate ring $\FK{G}$ is isomorphic to the coordinate ring $\FK{\GD}$ by the restriction map.

The Peter-Weyl theorem and Borel-Weil theorem of above give corresponding theorems for $\FK{\GD}$ and $\F[\GD]^U$, the $G^{op}\times G$ and $G^{op}$ actions now extended to  $\widehat{G}^{\,op}\times \GD$ and $\widehat{G}^{\,op}$ actions, but we do not need it in the article. 

Note that this is a reinterpretation of the coordinate rings used in \cite{KP2}, \cite{Kas}, and in subsequent articles of other authors. These coordinate rings should be considered as the coordinate rings of face monoids $\GD_I$, $\GD$. Classically, i.e., if the generalized Cartan matrix has only components of finite type, it is $\GD=\GD_I=G$. This explains why the face monoid $\GD$ may turn up in algebraic geometric constructions, where classically only the group $G$ is involved.\vspace*{1ex}

{\bf Some remarks to the use of opposite monoids:}
Let $M$ be a monoid. The opposite monoid ($M^{op}$, $\cdot_{op}$) is defined as follows: It is $M^{op}:=M$ as a set, and
\begin{eqnarray*}
m_1\cdot_{op}m_2:=m_2 m_1 \quad \mb{ where }\quad m_1,\,m_2\in M.
\end{eqnarray*}
Now let $V$ be a linear space. If $M$ acts  by $\centerdot$ linearly from the right  on $V$, then this corresponds by
\begin{eqnarray*}
   m\centerdot' v := v\centerdot m\qquad\mb{where }\quad  m\in M,\,v\in V,
\end{eqnarray*}
to a linear (left) action $\centerdot' $ of the opposite monoid $M^{op}$.  It corresponds by
\begin{eqnarray*}
   \pi(m)v:=v\centerdot m\qquad\mb{where }\quad m\in M,\,v\in V,
\end{eqnarray*}
to a monoid anti-homomorphism $\pi :M \to End(V)$, resp. to a representation (monoid homomorphism) $\pi :M^{op}\to End(V)$. 

In the article we need to work with such maps $\pi$. We speak about actions of $M^{op}$ on $V$ / $M^{op}$ acts on $V$ / representations of $M^{op}$ on $V$. In proofs we work with the monoid anti-homomorphisms $\pi :M \to End(V)$, which keeps our notation simple.

%
%
%
\section{Projective spectra of graded algebras}
In this section we briefly work out the theory of spectra of homogeneous prime ideals of graded algebras as far as we will need it in the article. We call such spectra projective spectra for short. Proofs which are similar to the nongraded or the $\Nn$-graded case are omitted. These cases can be found in Chapter II, Section 1 - 3 of \cite{Ha}, and in detail in Chapter 2 of \cite{Li}.

To explain the difference to the projective spectra considered in algebraic geometry let $A=\bigoplus_{n\in\Nn}A_n$ be an $\Nn$-graded $\F$-algebra such that $\F \,1\subseteq A_0$. For simplicity assume equality. Denote by $\PA$ all $\Nn$-homogeneous prime ideals of $A$ and by $ m:=\bigoplus_{n\in\N}A_n$ the irrelevant ideal. Then
\begin{eqnarray*}
   \PA= \PA\setminus\{ m \} \,\dot{\cup}\, \{ m \} .
\end{eqnarray*}
Suppose that $\PA\setminus\{  m \}\neq\emptyset $. Then $\PA\setminus\{  m \}$ can be made into a projective scheme as for example explained in \cite{Ha}, Chapter II, Section 2. However, it is only an open dense part of $\PA$. Its boundary is given by the closed point $m$. 

The example investigated in the following sections \ref{psCa} and \ref{afmb} is the projective spectrum $\PCA$ of the $P^+$-graded Cartan algebra $CA$ associated to a symmetrizable Kac-Moody group $G$. The projective spectrum $\PCA$ is a locally ringed space. There is an open dense part of $\PCA$ which is a scheme, but now also its boundary is interesting. Actually, $\PCA$ is stratified by schemes. The algebraic geometric action of the face monoid $\GD$ on $\PCA$, which we will obtain naturally, does not respect the strata but the closures of the strata. The existence of such an action is not at all obvious if we would only consider the union of these strata.

In the introduction and in sections \ref{Preli}, \ref{psCa}, \ref{afmb}, and \ref{quest} the field $\F$ is of characteristic 0. However in this section the field $\F$ may be of arbitrary characteristic.
\subsection{Projective spectra and their $\F$-valued points as topological spaces}
In this article a $\mathcal M$-graded $\F$-algebra $A$, or graded algebra $A$ for short, consists of the following:
\begin{itemize}
\item[(a)] A commutative monoid $({\mathcal M},\,+)$ with the cancellation property, i.e., 
\begin{eqnarray*}
  \La_1+N=\La_2+N \quad\Rightarrow \quad\La_1=\La_2
\end{eqnarray*}
 for all $\La_1$, $\La_2$, $N\in {\mathcal M}$.
\item[(b)] As a consequence of (a), the commutative monoid $\mathcal M$ embeds into a commutative group ${\mathcal M}-{\mathcal M}$, minimal over $\mathcal M$. We require ${\mathcal M}-{\mathcal M}$ to be torsion-free.
\item[(c)] A commutative, associative, $\F$-linear algebra $A$ with unit 1. Furthermore, a $\F$-linear decomposition  
\begin{eqnarray*}
   A=\bigoplus_{\La\in {\mathcal M}} A_\La \quad\mb{ such that }\quad \F\,1\subseteq A_0 \quad \mb{ and }\quad A_\La  \prd A_N\subseteq A_{\La+N} \mb{ for all }\La,\,N\in{\mathcal M}.
\end{eqnarray*}
\end{itemize}

For the $\mathcal M$-graded algebras considered in sections \ref{psCa} and \ref{afmb} the commutative group ${\mathcal M}-{\mathcal M}$ is isomorphic to some $\Z^m$, $m\in\Nn$.\\

Let $A$ be a $\mathcal M$-graded algebra. The projective spectrum of $A$ (as a topological space) consists of the set
\begin{eqnarray*}
      \PA :=\Mkl{{\mathcal M}\mb{-homogeneous prime ideals of }A}
\end{eqnarray*}
equipped with the Zariski topology, whose closed sets are obtained by
\begin{eqnarray*}
    {\mathcal V}(I):=\Mklz{Q\in\PA}{Q\supseteq I},\quad I \mb{ a $\mathcal M$-homogeneous ideal of }A.\\
\end{eqnarray*}

Fix a point $Q\in\PA$. Since $Q$ is a $\mathcal M$-homogeneous prime ideal it has the decomposition
\begin{eqnarray*}
    Q=\bigoplus_{\La\in {\mathcal M}} Q_\La\quad  \mb{ with }\quad Q_\La:=Q\cap A_\La ,
\end{eqnarray*} 
and $A\setminus Q$ is a multiplicatively closed subset of $A$. Also the set $\bigcup_{\La\in {\mathcal M}} A_\La \setminus Q_\La$ of homogeneous elements of $A\setminus Q$ is a multiplicatively closed subset of $A$. Now
\begin{eqnarray*}
    A_{(Q)}:=\Mklz{\frac{a}{b}}{ a \in A_\La,\, b\in A_\La\setminus Q_\La,\, \La\in {\mathcal M} \mb{ such that } A_\La\setminus Q_\La\neq \emptyset  } \subseteq \left(\bigcup_{\La\in {\mathcal M}} A_\La\setminus Q_\La\right)^{-1} A,
\end{eqnarray*}
which is the zero-homogeneous part of the graded algebra $\left(\bigcup_{\La\in {\mathcal M}} A_\La\setminus Q_\La\right)^{-1} A$, is a commutative, associative, unital, local, $\F$-linear algebra with maximal ideal
\begin{eqnarray*}
    m_{(Q)}:=\Mklz{\frac{a}{b}}{a\in Q_\La,\, b\in A_\La\setminus Q_\La,\, \La\in {\mathcal M} \mb{ such that } A_\La\setminus Q_\La\neq \emptyset }.
\end{eqnarray*}
We define the residue field of $Q$ to be 
\begin{eqnarray*}
\F_{(Q)}:=A_{(Q)}/m_{(Q)}.
\end{eqnarray*}
We identify the field $\F$ with the corresponding subfield of $\F_{(Q)}$.\\

We call $Q\in \mb{Proj}\,(A)$ an $\F$-valued point of the projective spectrum $\PA$ if $\F_{(Q)}=\F$. We set
\begin{eqnarray*}
\PAF:=\Mklz{Q\in \mb{Proj\,}(A)}{\F_{(Q)}=\F}.
\end{eqnarray*}
The following characterization will be useful to determine $\F$-valued points.
\begin{thm}\label{dim} Let $Q\in \mb{Proj\,}(A)$. The following are equivalent:
\begin{itemize}
\item[(i)] $Q\in \mb{Proj}\,(A)(\F)$.
\item[(ii)] $\mb{dim}\,( A_\La/Q_\La ) \leq 1$ for all $\La\in {\mathcal M}$.
\end{itemize}
\end{thm}
\Proof For every $\La\in {\mathcal M}$ we choose a linear complement $R_\La$ of $Q_\La$ in $A_\La$, i.e., $A_\La= Q_\La\oplus R_\La$. We first show
\begin{eqnarray}\label{resf}
 \F_{(Q)}=\Mklz{\frac{r}{s}+ m_{(Q)}} {r\in R_\La,\,s\in R_\La\setminus\{0\},\, \La\in {\mathcal M} \mb{ such that }R_\La\neq\{0\}}.
\end{eqnarray}
Obviously, the inclusion "$\supseteq$" holds. Now let $\frac{a}{b}+ m_{(Q)}\in \F_{(Q)}$ where $a \in A_\La$, $b\in A_\La\setminus Q_\La$ and $\La\in {\mathcal M}$ such that $A_\La\setminus Q_\La\neq \emptyset$. Decompose
\begin{eqnarray*}
  a=q_a+r   & \quad \mb{ with } \quad &   q_a\in Q_\La,\; r\in R_\La,\\
  b=q_b+s   & \quad \mb{ with } \quad &   q_b\in Q_\La,\; s\in R_\La\setminus\{0\}.
\end{eqnarray*}  
We find
\begin{eqnarray*}
  \frac{a}{b}-\frac{r}{s}=\frac{a s - r  b}{b s}=\frac{q_a  s - r  q_b}{b s} \in m_{(Q)},
\end{eqnarray*}
because of $q_a s$, $r q_b\in Q_\La  \prd R_\La\subseteq Q_{2\La}$ and  $b s\in (A_\La\setminus Q_\La)  \prd (A_\La\setminus Q_\La)\subseteq A_{2\La}\setminus Q_{2\La}$.

Now suppose that (ii) holds. Let $\La\in {\mathcal M}$ such that $R_\La\neq\{0\}$. If $r\in R_\La$, $s\in R_\La\setminus\{0\}$ then $r=\mu s$ for some $\mu\in \F$. Hence $\frac{r}{s}=\mu \frac{1}{1}$. Now (i) follows from (\ref{resf}).

Suppose that (i) holds. Let $\La\in {\mathcal M}$ such that $R_\La\neq\{0\}$. Choose $s\in R_\La\setminus\{0\}$. By (\ref{resf}) we find that for every $r\in R_\La$ there exists some $\mu\in\F$ such that
\begin{eqnarray*}
 \frac{r}{s}-\mu \frac{1}{1}=\frac{r-\mu s}{s}\in m_{(Q)} \,.
\end{eqnarray*}
Therefore, there exist $a\in Q_N$, $b\in A_N\setminus Q_N$, $N\in {\mathcal M}$ such that
\begin{eqnarray*}
  \frac{r-\mu s}{s}=\frac{a}{b} \,.
\end{eqnarray*}
It follows that there exists $c\in A_M\setminus Q_M$, $M\in {\mathcal M}$ such that
\begin{eqnarray}\label{resfgl}
   (r-\mu s)   b   c = a   s  c \, .
\end{eqnarray}
Since $a\in Q$ the right hand side of (\ref{resfgl}) is contained in $Q$. Assume that $r-\mu s\neq 0$. Since $r-\mu s\in R_\La\setminus\{0\}\subseteq A\setminus Q$ and $b$, $c \in A\setminus Q$, the left hand side of (\ref{resfgl}) is contained in $A\setminus Q$, which is not possible. It remains $r=\mu s$.
Therefore, (ii) holds.
\qed

For a subset $S\subseteq \PA$ we define
\begin{eqnarray*}
  S(\F):=S\cap \PAF.
\end{eqnarray*}
We equip $\PAF$ with its relative topology as a subset of $\PA$. We denote the relative closure of $Z\subseteq \PAF$ by $\overline{Z}^{\,pts}$. Trivially, $\overline{Z}^{\,pts}=\overline{Z}(\F)$.\\

Now we adapt some results from the theory of group rings, which for example can be found in Chapter 2, \S 6 of \cite{La}, to monoid rings. The next corollary follows from a remark in the middle of page 90 of \cite{La}, and from Theorems 6.29 and 6.31 of \cite{La}.
\begin{cor}\label{mr} Let $\mathcal M$ be a commutative monoid. The following are equivalent:
\begin{itemize}
\item[(i)]  The cancellation property holds in $\mathcal M$, and ${\mathcal M}-{\mathcal M}$ is torsion-free.
\item[(ii)] There exists a total order $\leq$ on $\mathcal M$ such that 
\begin{eqnarray*}
  \La_1 < \La_2 \quad\Rightarrow \quad\La_1+N < \La_2+N.
\end{eqnarray*}
for all $\La_1$, $\La_2$, $N\in {\mathcal M}$.
\item[(iii)] For every integral domain $R$ the monoid ring $R\,[\mathcal M]$ is an integral domain.
\item[(iv)] There exists an integral domain $R$ such that the monoid ring $R\,[\mathcal M]$ is an integral domain.
\end{itemize}
\end{cor}
\Proof  The proof of the direction ''(ii) $\Rightarrow$ (iii)'' is completely similar to the proof of Theorem 6.29 of \cite{La}. The direction ''(iii) $\Rightarrow$ (iv)'' is trivial.

 To ''(iv) $\Rightarrow$ (i)'': Denote by $e_\La\in R\,[\mathcal M]$ the image of $\La\in{\mathcal M}$ under the canonical injection $({\mathcal M},+)\to (R\,[{\mathcal M}],\prd)$. We first show the cancellation property. Let $\La_1$, $\La_2$, $N\in {\mathcal M}$ such that $\La_1+N=\La_2+N$. We get
\begin{eqnarray*} 
     ( e_{\La_1} - e_{\La_2} ) e_N =   e_{\La_1}e_N - e_{\La_2}  e_N =   e_{\La_1+N} - e_{\La_2+N} = 0
\end{eqnarray*}
in  $R\,[\mathcal M]$. Since $R\,[\mathcal M]$ has no zero divisors and $e_N\neq 0$ we find $e_{\La_1} = e_{\La_2}$. Therefore, $\La_1=\La_2$.

To show that  ${\mathcal M}-{\mathcal M}$ is torsion-free it is sufficient to show $n\La_1\neq n\La_2$ for all $\La_1,\,\La_2\in {\mathcal M}$, $\La_1\neq\La_2$, and all $n\in \N$. To do this we vary a remark in the middle of page 90 of \cite{La}.

Let $\La_1,\,\La_2\in {\mathcal M}$, $\La_1\neq\La_2$. Suppose there exists an element $n\in\N$ such that $n\La_1 = n\La_2$. We choose the  smallest element $n\in\N$ with this property. We get
\begin{eqnarray*}
     ( e_{\La_1} - e_{\La_2} )  \,\sum_{k=0}^{n-1} \,  e_{(n-1-k)\La_1+k\La_2} \hspace{50ex}\\
            =   ( e_{\La_1} - e_{\La_2} )   \,\sum_{k=0}^{n-1} \, (e_{\La_1})^{n-1-k}   (e_{\La_2})^k   =  (e_{\La_1}) ^n - (e_{\La_2} )^n =  e_{n\La_1} - e_{n\La_2} = 0
\end{eqnarray*}
in  $R\,[\mathcal M]$. Now let $k,\,k'\in\{0,\,1,\,\ldots,\,n-1\}$, $k\leq k'$, such that
\begin{eqnarray*}
   (n-1-k)\La_1 +k\La_2 =  (n-1-k')\La_1 +k'\La_2\,.
\end{eqnarray*}
By the cancellation property we get $(k'-k)\La_1=(k'-k)\La_2$. Furthermore, $0\leq k'-k< n$. The minimality of $n$ implies $k'=k$. This shows that the sum of the elements $e_{(n-1-k)\La_1+k\La_2}$, $k=0,\,1,\,\ldots,\,n-1$, is nonzero. Since $R\,[\mathcal M]$ has no zero divisors we find $e_{\La_1} = e_{\La_2}$, which contradicts $\La_1\neq \La_2$.

 To  ''(i) $\Rightarrow$ (ii)'': ${\mathcal M}-{\mathcal M}$ is a torsion-free commutative group. By Theorem 6.31 of \cite{La} there exists a total order $\leq$ on  ${\mathcal M}-{\mathcal M}$ compatible with the addition. Restricting this order to $\mathcal M$ shows (ii).
\qed

Let  $A$ be a $\mathcal M $-graded algebra. If $J$ is an ideal of $A$ then
\begin{eqnarray*}
  J^h:=\bigoplus_{\La\in{\mathcal M}}J\cap A_\La
\end{eqnarray*}
is an $\mathcal M $-homogeneous ideal of $A$ such that $J^h\subseteq J$. We call $J^h$ the $\mathcal M$-homogeneous ideal associated to $J$. The following theorem generalizes  Lemma 3.35 (a) of \cite{Li}.
\begin{thm}\label{homid} Let  $A$ be a $\mathcal M $-graded algebra. If $Q\in \mbox{Spec\,}(A)$ then  $Q^h\in \mbox{Proj\,}(A)$.
\end{thm}
\Proof The theorem is trivial for ${\mathcal M}=\{0\}$. Let  ${\mathcal M}\neq\{0\}$. We choose a total order $\leq$ on $\mathcal M$ with the property described in Corollary \ref{mr} (ii). For $x\in A$ with homogeneous components $x_\La\in A_\La$, $\La\in{\mathcal M}$, we set
\begin{eqnarray*}
  cs (x):=\{\La\in {\mathcal M}\mid x_\La\not\in Q\}\,.
\end{eqnarray*}

It is $Q^h\subseteq Q\neq A$. Let $x,\,x'\in A\setminus Q^h$. Then $cs(x)$ and $cs(x')$ are nonempty finite subsets of $\mathcal M$. Let $M$ be the biggest element of $cs(x)$ and $M'$ the biggest element of $cs(x')$.
If $\La\in cs(x)$, $\La'\in cs(x')$ such that $\La < M$ or $\La' < M'$ then $\La+\La' < M+M'$. Therefore, the only possibility to write $M+M'$ as a sum 
\begin{eqnarray*}
    M+M'=\La+\La'  \quad \mb{ with }\quad \La\in cs(x) \;\mb{ and }\;\La'\in cs(x')
\end{eqnarray*} 
is  $M+M'= M+M'$. This is used for the computation of the $M+ M'$-homogeneous component of $x x'$ modulo the linear space $Q\cap A_{M+M'}$:
\begin{eqnarray*}
    (x x')_{M+M'}  &=&  \sum_ {\La,\,\La' \in {\mathcal M}  ,\, \La+ \La' =M+ M' } x_\La x'_{\La'}    \\
                       &\equiv &  \sum_ {\La \in cs(x),\,\La'\in cs(x'),\, \La+  \La' = M+ M' } x_\La x'_{\La'} 
                           \:= \:  x_M x'_{M'}  \quad\mb{ mod }\, Q\cap A_{M+ M'} \,.
\end{eqnarray*}
Because $Q$ is prime and $x_M, \,x'_{M'}\notin Q$ it is  $x_M x'_{M'}\notin Q\cap A_{M+ M'}$. Therefore, $x x'\in A\setminus Q^h$.
\qed
 
The previous theorem can be used to generalize certain results for spectra to projective spectra.
\begin{thm}\label{ISP} Let  $A$ be a $\mathcal M $-graded algebra. Let $I$ be a $\mathcal M$-homo\-ge\-neous ideal of $A$ and $S$ a multiplicatively closed subset of $A$ such that $I\cap S=\emptyset$. Then there exists a $\mathcal M$-homo\-ge\-neous prime ideal $P\in \mbox{Proj\,}(A)$ such that $I\subseteq P$ and $P\cap S=\emptyset$.
\end{thm}
\Proof By Theorem 3.44  of \cite{Sh} there exists a prime ideal $Q\in \mb{Spec\,}(A)$ such that $I\subseteq Q$ and $Q\cap S=\emptyset$. By Theorem \ref{homid} we get a $\mathcal M$-homo\-ge\-neous prime ideal $Q^h\in \mbox{Proj\,}(A)$ such that $Q^h\subseteq Q$. Trivially,  $Q^h\cap S=\emptyset$. Since $I$ is homogeneous, $I\subseteq Q^h $.
\qed
\begin{cor} Let  $A$ be a $\mathcal M $-graded algebra.
\begin{itemize} 
\item[(a)] Then ${\mathcal V}(I)\neq \emptyset$ for every proper  $\mathcal M$-homogeneous ideal $I$ of $A$. In particular, $\PA\neq\emptyset$.
\item[(b)] If $I$ is a $\mathcal M$-homogeneous ideal of $A$ then also the radical $\sqrt{I}$ is a $\mathcal M$-homogeneous ideal of $A$ and
\begin{eqnarray}\label{intv}
  \bigcap_{P\in {\mathcal V}(I) } P = \sqrt{I}  \,,
\end{eqnarray}
an intersection over the empty set defined to be $A$.
\end{itemize}
\end{cor}
\begin{rem} A similar statement as in (a) does not hold for the $\F$-valued points, as it does not hold in the nongraded case. Take for example the $\R$-algebra $\C$, then $\mb{Spec\,}(\C)(\R)=\emptyset$.
\end{rem}
\Proof (a) follows from Theorem \ref{ISP} for $S=\{1\}$. Formula  (\ref{intv}) of (b) follows also from Theorem \ref{ISP}. It can be proved in the same way as Lemma 3.48 of \cite{Sh}, which gives the formula in the nongraded case. The ideal $\sqrt{I}$ is $\mathcal M$-homogeneous as an intersection of $\mathcal M$-homogeneous ideals.
\qed

We obtain from part (b) of the previous corollary:
\begin{cor}\label{VIVJ} Let  $I$, $J$ be $\mathcal M$-homo\-ge\-neous  ideals of the $\mathcal M $-graded algebra $A$. Then ${\mathcal V}(I)\subseteq {\mathcal V}(J)$ if and only if $\sqrt{I}\supseteq \sqrt{J}$.
\end{cor}

Let $A$ be a ${\mathcal M}_A$-graded algebra and $B$ be a ${\mathcal M}_B$-graded algebra. A morphism $f^*:\,A\to B$ of graded algebras is a unital morphism of algebras $f^*:\,A\to B$ for which there exists a map $F^*:{\mathcal M}_A\to {\mathcal M}_B$ such that 
\begin{eqnarray*}
    f^*(A_\La) \subseteq B_{F^*(\La)} \quad \mb{ for all } \quad\La\in {\mathcal M}_A .
\end{eqnarray*}

Note that the map $F^*$ is uniquely determined only on $\Mklz{\La\in {\mathcal M}_A} {f^*(A_\La)\neq\{0\}}$. Outside this set the map $F^*$ is not relevant and can be chosen arbitrarily. The definition has been formulated in this way to avoid working with partial maps.

Let $f^*:\,A\to B$ be a morphism of graded algebras and $J$ a homogeneous ideal of $B$. If $F^*$ is injective on $\Mklz{\La\in {\mathcal M}_A} {f^*(A_\La)\neq\{0\}}$ then the inverse image $(f^*)^{-1}(J)$ is a homogeneous ideal of $A$, i.e.,
\begin{eqnarray*}
(f^*)^{-1}(J)= (f^*)^{-1}(J)^h.
\end{eqnarray*}
If $F^*$ is not injective on $\Mklz{\La\in {\mathcal M}_A} {f^*(A_\La)\neq\{0\}}$ then the inverse image $(f^*)^{-1}(J)$ does not need to be a homogeneous ideal of $A$, i.e.,
\begin{eqnarray*}
(f^*)^{-1}(J)\neq (f^*)^{-1}(J)^h
\end{eqnarray*}
is possible. For this reason we work with
\begin{eqnarray*}
   (f^*)^{-1}(J)^h = \bigoplus_{\La\in  {\mathcal M}_A }\left(f^*\mid_{A_\La}^{ B_{F^*(\La)} }\right)^{-1}(J_{F^*(\La)}) =
   \bigoplus_{\La\in  {\mathcal M}_A \atop  f^*(A_\La)=\{0\}}  A_\La  \;\; \oplus \bigoplus_{\La\in  {\mathcal M}_A \atop f^*(A_\La)\neq\{0\}}\left(f^*\mid_{A_\La}^{ B_{F^*(\La)} }\right)^{-1}(J_{F^*(\La)})
\end{eqnarray*}
instead of the inverse image $(f^*)^{-1}(J)$. It is easy to prove the following:
\begin{thm}\label{morph} Let $f^*:\,A\to B$ be a morphism of graded algebras. 
\begin{itemize}
\item[(a)]We get a continuous map by
\begin{eqnarray*}
  \begin{array}{cccc}
    f: &\PB & \to & \PA\\
         & Q  &\to & (f^*)^{-1}(Q)^h 
  \end{array}.
\end{eqnarray*}
It maps $\PBF$ into $\PAF$. 
\item[(b)] For every $Q\in \PB$ we get a morphism of algebras
\begin{eqnarray*}
  \begin{array}{cccc}
    f^*_Q: & A_{(f(Q))} & \to & B_{(Q)}\\
         & \frac{a}{b}  &\to & \frac{f^*(a)}{f^*(b)} 
  \end{array},
\end{eqnarray*}
which is local, i.e., $f^*_Q(m_{(f(Q))})\subseteq m_{(Q)}$. It induces an injective morphism $f^*_Q:\,\F_{(f(Q))}\to\F_{(Q)}$ of fields over $\F$. 
\end{itemize}
\end{thm}
\begin{rem} Let  $A$ be a $\mathcal M $-graded algebra. Let $f^*:\,A\to A$ be the identity map, the domain equipped with the $\mathcal M$-gradation, the target equipped with the trivial gradation. This map is a morphism of graded algebras, the associated map $f:\mb{Spec\,}(A)\to \mb{Proj\,}(A)$ given by $f(Q)=Q^h$, $Q\in \mb{Spec\,}(A)$.
\end{rem}

The graded algebras and its morphisms form a category, in which the concatenation of morphisms is the concatenation of maps. If we assign to a graded algebra its projective spectrum and to a morphism of graded algebras the map of Theorem \ref{morph} (a) we get a contravariant functor into the category of topological spaces. 
\begin{rem} If $f^*:A\to B$ is an isomorphism of graded algebras then a map $F^*$ as in the definition maps the set $\Mklz{\La\in {\mathcal M}_A} {A_\La\neq\{0\}}$ bijectively to $\Mklz{N\in {\mathcal M}_B} {B_N\neq\{0\}}$. But the monoids ${\mathcal M}_A$ and ${\mathcal M}_B$ may be non-isomorphic, even non-bijective. 

 A morphism $f^*:\, A\to B$ of graded algebras is an isomorphism of graded algebras if and only if the map $f^*:A\to B$ is bijective and the map  $F^*:\Mklz{\La\in {\mathcal M}_A} {f^*(A_\La)\neq\{0\}}\to {\mathcal M}_B$ is injective.
\end{rem}

Let $A$ be a $\mathcal M$-graded algebra. Let $S\subseteq A\setminus\{0\}$ be a multiplicatively closed subset of homogeneous elements. For $\La\in {\mathcal M}$ denote by $S_\La:=S\cap A_\La$ the $\La$-homogeneous elements of $S$. Then ${\mathcal M}_S:=\Mklz{\La\in {\mathcal M}}{S_\La\neq\emptyset}$ is a submonoid of $\mathcal M$ and ${\mathcal M}-{\mathcal M}_S:=\Mklz{\La-N\in{\mathcal M}-{\mathcal M}}{\La\in{\mathcal M},\,N\in {\mathcal M}_S}$ is a submonoid of the group ${\mathcal M}-{\mathcal M}$. Set
\begin{eqnarray*}
        D(S)    :=    \Mklz{Q\in \PA}{S\subseteq A\setminus Q}.
\end{eqnarray*}
The localization $S^{-1}A$ is a commutative associative unital $({\mathcal M}-{\mathcal M}_S)$-graded $\F$-algebra. The canonical morphism $i^*:\,A\to S^{-1}A$ is a morphism of graded algebras. As map $I^*: {\mathcal M}\to {\mathcal M}-{\mathcal M}_S$ we can take the canonical embedding.
It is not difficult to prove the following:
\begin{thm}\label{loc} 
\begin{itemize}
\item[(a)] The continuous map 
\begin{eqnarray*}
\begin{array}{ccc}
  i:\mb{Proj\,}(S^{-1}A) & \to & \PA\\
    Q &\mapsto & (i^*)^{-1}(Q)
\end{array}
\end{eqnarray*}
maps $\mb{Proj\,}(S^{-1}A)$ homeomorphically to $D(S)$. It also maps $\mb{Proj\,}(S^{-1}A)(\F)$ homeomorphically to $D(S)(\F)$.
\item[(b)] For every $Q\in \mb{Proj\,}(S^{-1}A)$ the local morphism $i^*_Q:\,A_{(i(Q))}\to (S^{-1}A)_{(Q)}$ is an isomorphism. 
\end{itemize}
\end{thm}

In particular, the theorem applies to the principal open sets
\begin{eqnarray*}
   D(a):=\Mklz{Q\in\PA}{a\notin Q} 
\end{eqnarray*}
with $a\in A$ homogeneous and not nilpotent, because it holds $D(a)=D(\Mklz{a^n}{n\in\Nn})$. These sets form a base of the open sets of the Zariski topology of $\PA$.\\

Let $A$ be a $\mathcal M$-graded algebra. Let $S\subseteq A\setminus\{0\}$ be a multiplicatively closed subset of homogeneous elements of $A$. We call $s_0\in S$ a {\it principal element} of the monoid $S$ if the following property holds: For every $s\in S$ there exist $\ti{s}\in S$ and $n\in\Nn$ such that $s\ti{s}=s_0^n$.

A subset $I\subseteq S$ is called a semigroup ideal of the monoid $S$ if $S\prd I\subseteq I$. Note that we also allow $I$ to be the empty set. A submonoid $F\subseteq S$, for which $S\setminus F$ is a semigroup ideal of $S$, is called a face of $S$. The relative interior $ri(S)$ of $S$ is the semigroup ideal of $S$ defined by
\begin{eqnarray*}
  ri(S):= S\setminus\bigcup_{F  \,\textrm{\footnotesize a face of } S,\; F\neq S}F .
\end{eqnarray*}

It is easy to prove the following:
\begin{thm}\label{poset} Suppose there exists a principal element $s_0\in S$.
\begin{itemize}
\item[(a)] Then $D(s_0)=D(S)$.
\item[(b)] An isomorphism of graded algebras is given by
\begin{eqnarray*}
\begin{array}{ccc}
   \Mklz{s_0^n}{n\in\Nn}^{-1}A &\to& S^{-1}A\\
         \frac{a}{s_0^n} &\mapsto &\frac{a}{s_0^n}
\end{array}.
\end{eqnarray*} 
\end{itemize}
Furthermore, the set of principal elements of $S$ is a semigroup ideal of $S$ contained in $ri(S)$.
\end{thm}

If there exists a principal element of $S$ we also call $D(S)$ a {\it principal open set}.
\begin{rem} If $S$ is generated by finitely many elements $s_1$, $s_2$,\ldots, $s_k$ then $s_0:=s_1s_2\cdots s_k$ is a principal element of $S$. It is easy to find examples where $S$ is not finitely generated and principal elements exist.
\end{rem}

Let $A$ be a $\mathcal M$-graded algebra. Let $I$ be a homogeneous ideal of $A$, $I\neq A$. The quotient algebra $A/I$ is $\mathcal M$-graded. The canonical projection $\pi^*:\,A\to A/I$ is a morphism of graded algebras.   As map $\Pi^*:{\mathcal M}\to {\mathcal M}$ we can take the identity map.
It is not difficult to prove the following:
\begin{thm}\label{csub} 
\begin{itemize}
\item[(a)] The continuous map 
\begin{eqnarray*}
\begin{array}{ccc}
  \pi:\mb{Proj\,}(A/I) & \to & \PA\\
    Q &\mapsto & (\pi^*)^{-1}(Q)
\end{array}
\end{eqnarray*}
maps $\mb{Proj\,}(A/I)$ homeomorphically to ${\mathcal V}(I)$. It also maps $\mb{Proj\,}(A/I)(\F)$ homeomorphically to ${\mathcal V}(I)(\F)$.
\item[(b)] For every $Q\in \mb{Proj\,}(A/I)$ the local morphism $\pi^*_Q:\,A_{(\pi(Q))}\to (A/I)_{(Q)}$ is surjective.\\
\end{itemize}
\end{thm}

If $A$ is a $\mathcal M$-graded algebra and $B$ is a $\mathcal N$-graded algebra then the tensor product $A\otimes B$ of the algebras $A$ and $B$ over $\F$ is a ${\mathcal M}\oplus {\mathcal N}$-graded algebra. The canonical injections 
\begin{eqnarray*}
\begin{array}{rccc}
  i_A^*: & A & \to & A\otimes B \\
       & a & \mapsto &a\otimes 1
\end{array}
\quad\mb{ and }\quad
\begin{array}{rccc}
  i_B^*: & B & \to & A\otimes B \\
       & b & \mapsto & 1\otimes b
\end{array}
\end{eqnarray*}
are morphisms of graded algebras. As maps $I_A^*:{\mathcal M}\to {\mathcal M}\oplus {\mathcal N}$ and $I_B^*:{\mathcal N}\to {\mathcal M}\oplus {\mathcal N}$ we can take the canonical embeddings. 
It is not difficult to prove the following:
\begin{thm}\label{prod}
\begin{itemize}
\item[(a)] The map 
\begin{eqnarray*}
\begin{array}{ccc}
    i_A:\mb{Proj\,}(A\otimes B) & \to & \PA\\
    Q &\mapsto & (i_A^*)^{-1}(Q) 
\end{array}
\end{eqnarray*}
is a continuous map, mapping $\mb{Proj\,}(A\otimes B)(\F)$ into $\PAF$.
\item[(b)] For every $Q\in \mb{Proj\,}(A\otimes B)$ the local morphism $(i_A^*)_Q:\,A_{( i_A(Q))}\to (A\otimes B)_{(Q)}$ is injective.\\
\end{itemize}
\end{thm}
Let $\mathcal N$ be a commutative monoid with the cancellation property, such that ${\mathcal N}-{\mathcal N}$ is torsion-free. Its monoid algebra 
\begin{eqnarray*}
      \FK{\mathcal N}=\bigoplus_{N\in {\mathcal N}} \F e_N,\quad   e_N\prd e_{N'}=e_{N+N'} \;\mb{ for all }\; N,\,N'\in {\mathcal N},
\end{eqnarray*} 
over the field $\F$ is a $\mathcal N$-graded algebra. By Corollary \ref{mr} it has no zero divisors.
\begin{thm}\label{prodlattice} Let $A$ be a $\mathcal M$-graded algebra. Let $\mathcal N$ be a  torsion-free commutative group. The map 
\begin{eqnarray*}
\begin{array}{cccc}
  i_A : & \mb{Proj\,}(A\otimes \FK{\mathcal N}) & \to     & \PA\\
        &    Q                                  & \mapsto & (i_A^*)^{-1}(Q)
\end{array}
\end{eqnarray*}
is a homeomorphism, mapping $\mb{Proj\,}(A\otimes \FK{\mathcal N})(\F)$ bijectively to $\PAF$.
Its inverse is
\begin{eqnarray*}
\begin{array}{cccc}
  (i_A)^{-1}: &\PA  &\to     & \mb{Proj\,}(A\otimes \FK{\mathcal N})\\
              & Q   &\mapsto & Q\otimes \FK{\mathcal N}
\end{array}.
\end{eqnarray*}
For every $Q\in \PA$ the local morphism $(i_A^*)_{ Q\otimes\,\F\,[{\mathcal N}] } : A_{ (Q) }\to (A\otimes \FK{\mathcal N} )_{  (Q\otimes\,\F\,[{\mathcal N}])  } $ is an isomorphism. 
\end{thm}
\begin{rem}  For $A=\F$ the theorem specializes to: $\mb{Proj\,}(\FK{\mathcal N})  =  \mb{Proj\,}(\FK{\mathcal N})(\F)  =  \left\{\{0\}\right\}$ for the $\mathcal N$-graded group algebra $\FK{\mathcal N}$ of a  torsion-free commutative group $\mathcal N$.
\end{rem}
\Proof The map $i_A$ is a well defined continuous map by Theorem \ref{prod} (a). We next show that the map $(i_A)^{-1}$ given in the theorem is well defined and $i_A\circ (i_A)^{-1}= id$. If $Q\in \PA$ then $Q\otimes \FK{\mathcal N}$ is a homogeneous ideal of $A\otimes \FK{\mathcal N}$ and 
\begin{eqnarray*}
   (i_A^*)^{-1}(Q\otimes \FK{\mathcal N})=\Mklz{a\in A}{a\otimes 1\in Q\otimes \FK{\mathcal N}} = Q.
\end{eqnarray*}
Because $A/Q$ has no zero divisors also
\begin{eqnarray*}
     (A\otimes \FK{\mathcal N})/(Q\otimes\FK{\mathcal N})\cong (A/Q)\otimes\FK{\mathcal N}\cong ((A/Q)\otimes\F)\,[\mathcal N]\cong (A/Q)\,[\mathcal N]
\end{eqnarray*} 
has no zero divisors by Corollary \ref{mr}. Therefore, the homogeneous ideal $Q\otimes \FK{\mathcal N}$ is prime.

Now let $J$ be a homogeneous ideal of $A\otimes \FK{\mathcal N}$. Then $(i_A^*)^{-1}(J)$ is a homogeneous ideal of $A$ because $i_A^*$ is a morphism of graded algebras. We show 
\begin{eqnarray}\label{jiastar}
  J = (i_A^*)^{-1}(J)\otimes\FK{\mathcal N}.
\end{eqnarray}
Since $J$ is an ideal, the inclusion "$\supseteq$" holds. Now let $x\in J$. Write $x$ in the form
\begin{eqnarray*}
    x = \sum_{N\in {\mathcal N}} a_N\otimes e_N \quad\mb{ with } \quad a_N\in A,\; a_N\neq 0 \mb{ for only finitely many }N.
\end{eqnarray*}
Because $J$ is homogeneous we get $a_N\otimes e_N\in J$ for all $N\in {\mathcal N}$. Since $J$ is an ideal we find
\begin{eqnarray*}
    a_N\otimes 1 =(a_N\otimes e_N)\prd (1\otimes e_{-N})\in J \quad\mb{ for all }\quad N\in {\mathcal N},
\end{eqnarray*}
which is equivalent to $a_N\in (i_A^*)^{-1}(J)$ for all $N\in{\mathcal N}$, which in turn is equivalent to $x\in (i_A^*)^{-1}(J)\otimes\FK{\mathcal N}$.

Equation (\ref{jiastar}) shows $(i_A)^{-1}\circ i_A=id$. Furthermore, it implies $i_A({\mathcal V}(J))={\mathcal V}((i_A^*)^{-1}(J))$ for every homogeneous ideal $J$ of $A\otimes \FK{\mathcal N}$, which shows that $(i_A)^{-1}$ is continuous.

Now let $Q\in\PA$. By Theorem \ref{prod} (b) the local morphism $(i_A^*)_{ Q\otimes\,\F\,[{\mathcal N}] }$ is injective. Let $M\in {\mathcal M}$ and $N\in {\mathcal N}$. For $a\otimes e_N\in A_M\otimes e_N$ and $b\otimes e_N\in A_M\otimes e_N\setminus Q_M\otimes e_N$ it holds
\begin{eqnarray*}
    (i_A^*)_{ Q\otimes\,\F\,[{\mathcal N}] }(\frac{a}{b})  =  \frac{a\otimes 1}{b\otimes 1}  =  \frac{a\otimes e_N}{b\otimes e_N}.
\end{eqnarray*}
Therefore, $(i_A^*)_{ Q\otimes\,\F\,[{\mathcal N}] }$ is also surjective.
It follows that $i_A$ maps $\mb{Proj\,}(A\otimes \FK{\mathcal N})(\F)$ bijectively to $\PAF$.
\qed

The $\F$-valued points of the projective spectrum of the tensor product of graded algebras can be described by the following theorem.
\begin{thm}\label{prodF} Let $A$ be a $\mathcal M$-graded algebra. Let $B$ a $\mathcal N$-graded algebra such that ${\mathcal N}-{\mathcal N}$ is tor\-sion-free. Then:
\begin{itemize}
\item[(a)] $\mb{Proj\,}(A\otimes B)(\F)\neq \emptyset \;\;\iff \; \; \PAF\neq\emptyset \;\mb{ and }\;\PBF\neq\emptyset$.
\item[(b)] The map 
\begin{eqnarray*}
\begin{array}{rccc}
  (i_A,\,i_B):&\mb{Proj\,}(A\otimes B)(\F) & \to & \PAF\times \PBF\\
    &Q &\mapsto & ((i_A^*)^{-1}(Q),\,(i_B^*)^{-1}(Q))
\end{array}
\end{eqnarray*}
is a bijective continuous map. Its inverse is given by
\begin{eqnarray*}
\begin{array}{rccc}
  (i_A,\,i_B)^{-1}:&\PAF\times\PBF &\to & \mb{Proj\,}(A\otimes B)(\F)\\
   &(R,\,S) &\mapsto & R\otimes B + A\otimes S  
\end{array}.
\end{eqnarray*}
\end{itemize}
\end{thm}
\Proof 
To (a): If $Q\in \mb{Proj\,}(A\otimes B)(\F)$ then by Theorem \ref{prod} (a) we obtain $(i_A^*)^{-1}(Q)\in \PAF$ and $(i_B^*)^{-1}(Q)\in \PBF$.

Now let $R\in \PAF$ and $S\in \PBF$. Obviously, $R\otimes B + A\otimes S$ is a homogeneous ideal of $A\otimes B$ with homogeneous parts
\begin{eqnarray*}
  (R\otimes B + A\otimes S)_{(M,\,N)}=R_M\otimes B_N +A_M\otimes S_N \quad\mb{ for }\quad  M\in {\mathcal M}, \;N\in {\mathcal N}.
\end{eqnarray*}
To show that the ideal  $R\otimes B + A\otimes S$ is prime it is sufficient to show that its quotient algebra
\begin{eqnarray*}
  A\otimes B/(R\otimes B + A\otimes S)&\cong & A/R\otimes B/S
\end{eqnarray*}
has no zero divisors. We do this by a variant of the proof of Theorem 6.29 in \cite{La}. We choose a total order $\leq$ on $\mathcal N$ with the properties described in Corollary \ref{mr} (ii). By Theorem \ref{dim} we have
\begin{eqnarray*}
     \mb{dim}((B/S)_N) = \mb{dim}(B_N/S_N)\leq 1 \quad\mb{ for all }\quad N\in\mathcal N\,. 
\end{eqnarray*}
Therefore, $u,\,v\in (A/R\otimes B/S)\setminus\{0\}$ can be written in the form
\begin{eqnarray*}
  u= c_1\otimes d_{N_1} + r_u \mb{ with } c_1\in( A/R )\setminus\{0\} ,\, d_{N_1}\in (B/S)_{N_1}\setminus\{0\}\mb{ and }r_u\in \bigoplus_{N\in {\mathcal N},\,N>N_1  } A/R\otimes (B/S)_N\,,\\
  v= c_2\otimes d_{N_2} + r_v \mb{ with } c_2\in (A/R)\setminus\{0\},\, d_{N_2}\in (B/S)_{N_2}\setminus\{0\}\mb{ and }r_v\in \bigoplus_{N\in {\mathcal N},\,N>N_2  } A/R\otimes (B/S)_N\,.
\end{eqnarray*}
It follows that
\begin{eqnarray*}
   uv= c_1c_2\otimes \underbrace{ d_{N_1}d_{N_2} }_{ \in (B/S)_{N_1+N_2} }  +\, r \quad \mb{ with } \quad r\in \bigoplus_{N\in {\mathcal N},\,N>N_1+N_2 } A/R\otimes (B/S)_N\,.
\end{eqnarray*}
The algebras $A/R$ and $B/S$ have no zero divisors because $R$ and $S$ are prime ideals. Therefore, $c_1c_2\neq 0$ and $ d_{N_1}d_{N_2}\neq 0$, which shows $uv\neq 0$.

Next we show that $R\otimes B + A\otimes S$ is an $\F$-valued point of $\mb{Proj\,}(A\otimes B)$. By Theorem \ref{dim} it is sufficient to show that for every $ M\in {\mathcal M}$ and $N\in {\mathcal N}$ the linear space
\begin{eqnarray*}
   (A\otimes B)_{(M,\,N)}/(R\otimes B + A\otimes S)_{(M,\,N)} 
      &\cong & A_M/R_M\otimes B_N/S_N   
\end{eqnarray*}
is at most one-dimensional. 
This follows from Theorem \ref{dim} applied to the $\F$-valued point $R$ of $\PA$ and the $\F$-valued point $S$ of $\PB$: 
\begin{eqnarray*}
    \mb{dim}(A_M/R_M\otimes B_N/S_N)=\mb{dim}(A_M/R_M)\mb{dim}(B_N/S_N)\leq 1\,.
\end{eqnarray*}

To (b): By Theorem \ref{prod} (a) the map $(i_A,\,i_B)$ is continuous. It remains to show that the maps of the theorem are inverse maps.
For $R\in\PAF$ and $S\in \PBF$ we get
\begin{eqnarray*}
  (i_A^*)^{-1}(R\otimes B + A\otimes S)= \Mklz{a\in A}{a\otimes 1\in R\otimes B + A\otimes S}=R,
\end{eqnarray*}
because of $1\notin S$. Similarly, we have $(i_B^*)^{-1}(R\otimes B + A\otimes S)=S$. Therefore, $(i_A,\,i_B)\circ(i_A,\,i_B)^{-1}=id$. 

Now we show $(i_A,\,i_B)^{-1}\circ(i_A,\,i_B)=id$, which is equivalent to 
\begin{eqnarray*}
    i_A(Q)\otimes B + A\otimes i_B(Q) = Q   \quad\mb{ for all }\quad Q\in \mb{Proj\,}(A\otimes B)(\F).
\end{eqnarray*}
For all $M\in {\mathcal M}$ and $N\in {\mathcal N}$ we have the inclusions
\begin{eqnarray}\label{pspab}
  i_A(Q)_M\otimes B_N + A_M\otimes i_B(Q)_N \subseteq Q_{M+N}\subseteq A_M\otimes B_N
\end{eqnarray}
where $\mb{dim}\left(A_M/i_A(Q)_M\right)\leq 1$ and $\mb{dim}\left(B_N/i_B(Q)_N \right)\leq 1$  by Theorem \ref{dim} .
If $i_A(Q)_M=A_M$ or $i_B(Q)_N =B_N$ then in ($\ref{pspab}$) both inclusions are equalities.
Now suppose that $i_A(Q)_M\neq A_M$, $i_B(Q)_N\neq B_N$, and $i_A(Q)_M\otimes B_N + A_M\otimes i_B(Q)_N \subsetneqq Q_{M+N}$. Then $Q_{M+N}=A_M\otimes B_N$. Furthermore, there exist elements $a\in A_M\setminus i_A(Q)_M$, $b\in B_N\setminus i_B(Q)_N$. We obtain
\begin{eqnarray*}
  \underbrace{(a\otimes 1)}_{\notin Q}  \prd \underbrace{(1\otimes b)}_{\notin Q} =a\otimes b\in A_M \otimes B_N = Q_{M+N},
\end{eqnarray*}
which contradicts that $Q$ is prime.
\qed
\subsection{The structure sheaves of projective spectra} 
For every $\mathcal M$-graded algebra $A$ we now introduce a locally ringed space ($\PA$, ${\mathcal O}_{\textrm{\footnotesize Proj\,}(A)}$). 
For a nongraded algebra $A$ it coincides with the one introduced in the second section of Chapter II of \cite{Ha}. For an $\Nn$-graded algebra $A$ its restriction to the open set $\mb{Proj\,}(A)\setminus{\mathcal V}(m)$, $m$ the irrelevant ideal of $A$, coincides with the  locally ringed space introduced in the second section of Chapter II of \cite{Ha}. (Please note: The set $\mb{Proj\,}(A)\setminus{\mathcal V}(m)$ is denoted by $\mb{Proj\,}(A)$ in \cite{Ha}.) 

The following proposition is not difficult to prove.
\begin{prop}\label{presheaf1} Let $A$ be a $\mathcal M$-graded algebra.
\begin{itemize}
\item[(a)] A presheaf $\widetilde{\mathcal O}_{\textrm{\footnotesize Proj\,}(A)}$ of $\F$-algebras on $\PA$ is obtained as follows: For $\emptyset\neq U\subseteq \PA$, $U$ open, let $\widetilde{\mathcal O}_{\textrm{\footnotesize Proj\,}(A)}(U)$ be the commutative, associative, unital, $\F$-linear algebra
\begin{eqnarray*}
   \widetilde{\mathcal O}_{\textrm{\footnotesize Proj\,}(A)}(U)
   := \Mklz{\frac{a}{b}\,}{ {a \in A_\La,\, b\in \bigcap_{Q\in U} A_\La\setminus Q_\La,\,\atop
      \La\in {\mathcal M} \mb{ such that } \bigcap_{Q\in U} A_\La\setminus Q_\La\neq \emptyset } } 
    \subseteq     \left(\bigcup_{\La\in{\mathcal M}}\bigcap_{Q\in U} A_\La\setminus Q_\La  \right)^{-1}A.
\end{eqnarray*}
For $\emptyset\neq U_1\subseteq U_2\subseteq \PA$, $U_1$ and $U_2$ open, let 
\begin{eqnarray*} 
\begin{array}{rccc}
 \widetilde{res}^{U_2}_{U_1}:&\widetilde{\mathcal O}_{\textrm{\footnotesize Proj\,}(A)}(U_2) &\to & \widetilde{\mathcal O}_{\textrm{\footnotesize Proj\,}(A)}(U_1) \\
    &\frac{a}{b} &\mapsto & \frac{a}{b}
\end{array}
\end{eqnarray*}
be the restriction morphism.
\item[(b)] The stalk in $Q\in\PA$ and its restriction morphisms are given by the local $\F$-algebra $A_{(Q)}$ and the morphisms
\begin{eqnarray*} 
\begin{array}{rccc}
 \widetilde{res}^U_Q:& \widetilde{\mathcal O}_{\textrm{\footnotesize Proj\,}(A)}(U) &\to & A_{(Q)} \\
    &\frac{a}{b} &\mapsto & \frac{a}{b}
\end{array},
\end{eqnarray*}
where $U$ is open in $\PA$ with $Q\in U$.
\end{itemize}
\end{prop}

The localisation $ \left(\bigcup_{\La\in{\mathcal M}}\bigcap_{Q\in U} A_\La\setminus Q_\La  \right)^{-1}A$ in Proposition \ref{presheaf1} (a) is a graded algebra. The algebra $\widetilde{\mathcal O}_{\textrm{\footnotesize Proj\,}(A)}(U)$ coincides with its zero-homogeneous part.
\begin{prop}\label{presheaf2}  Let $A$ be a $\mathcal M$-graded algebra.  Let $a\in A$ be homogeneous, not nilpotent. We get an isomorphism of graded algebras by
\begin{eqnarray*} 
\begin{array}{ccc}
    \{a^n|n\in\Nn\}^{-1}A  &\to &  \left(\bigcup_{\La\in{\mathcal M}}\bigcap_{Q\in D(a)} A_\La\setminus Q_\La  \right)^{-1}A \\
    \frac{b}{a^n} &\mapsto & \frac{b}{a^n}
\end{array}.
\end{eqnarray*}
It restricts to an isomorphism of algebras from $\left(\{a^n|n\in\Nn\}^{-1}A \right)_0$, the zero-homogeneous part of $\{a^n|n\in\Nn\}^{-1}A$, to $\widetilde{\mathcal O}_{\textrm{\footnotesize Proj\,}(A)}(D(a))$.
\end{prop}
\Proof Set $S:= \bigcup_{\La\in{\mathcal M}}\bigcap_{Q\in D(a)} A_\La\setminus Q_\La $. Trivially, $a\in S$. Because of Theorem \ref{poset} it is sufficient to show that $a$ is a principal element of $S$. 
The idea for the following proof has been isolated from a consideration in the proof of Proposition 2.2 (b) in Chapter II of \cite{Ha}, the last paragraph on page 71. For $x\in A$ we denote by $(x)$ the ideal generated by $x$.

Let $s\in S$. By the definition of $S$ we have $s\notin Q$ for all $Q\in D(a)$, which is equivalent to $D(a)\subseteq D(s)$. This in turn is equivalent to $\mathcal{V}((a))\supseteq \mathcal{V}((s))$, from which we find $a\in\sqrt{(a)}\subseteq\sqrt{(s)}$ by Corollary \ref{VIVJ}. Therefore, there exist
$n\in \N$, $\ti{s}\in A$ such that $a^n = \ti{s} s$. Since $a$ and $s$ are homogeneous we can find a homogeneous $\ti{s}\in A$ satisfying this equation. For every $Q\in D(a)$ we have $s\ti{s}=a^n\in A\setminus Q$. Since $Q$ is an ideal we obtain $\ti{s}\in A\setminus Q$. This shows $\ti{s}\in S$. \qed

The structure sheaf ${\mathcal O}_{\textrm{\footnotesize Proj\,}(A)}$ on $\PA$ is the sheafification of the presheaf $\widetilde{\mathcal O}_{\textrm{\footnotesize Proj\,}(A)}$ on $\PA$. The stalk in $Q\in\PA$ identifies with the local $\F$-algebra $A_{(Q)}$.

The following theorem generalizes Proposition 2.2 (b) as well as a part of Proposition 2.5 (b) in Chapter II of \cite{Ha}. 
\begin{thm}  Let $A$ be a $\mathcal M$-graded algebra.  Let $a\in A$ be homogeneous, not nilpotent. Then
\begin{eqnarray*}
   {\mathcal O}_{\textrm{\footnotesize Proj\,}(A)}(D(a)) \cong  \widetilde{\mathcal O}_{\textrm{\footnotesize Proj\,}(A)}(D(a)) \cong  \left(\{a^n|n\in\Nn\}^{-1}A \right)_0 ,
\end{eqnarray*}
where the first isomorphism is given by the canonical morphism $ \widetilde{\mathcal O}_{\textrm{\footnotesize Proj\,}(A)}(D(a)) \to {\mathcal O}_{\textrm{\footnotesize Proj\,}(A)}(D(a)) $ of the sheafification, and the second isomorphism $\left(\{a^n|n\in\Nn\}^{-1}A \right)_0\to  \widetilde{\mathcal O}_{\textrm{\footnotesize Proj\,}(A)}(D(a))$ is described in Proposition \ref{presheaf2}.

In particular, ${\mathcal O}_{\textrm{\footnotesize Proj\,}(A)}(\mb{Proj\,}(A)) \cong  \widetilde{\mathcal O}_{\textrm{\footnotesize Proj\,}(A)}(\mb{Proj\,}(A)) \cong A_0$.
\end{thm}
\Proof The  proof of Proposition 2.2 (b) in Chapter II of \cite{Ha}, the last paragraph on page 71 omitted, can be adapted to the graded case to prove the theorem. (To do this note: The structure sheaf ${\mathcal O}_{\textrm{\footnotesize Spec\,}(A)}$ on $\mb{Spec\,}(A)$ as well as the morphisms  $\{a^n|n\in\Nn\}^{-1}A \to{\mathcal O}_{\textrm{\footnotesize Spec\,}(A)}(D(a))$, $a\in A$ not nilpotent, are defined directly in \cite{Ha}. The presheaf $\widetilde{\mathcal O}_{\textrm{\footnotesize Spec\,}(A)}$ has not been introduced.)
\qed

For every morphism $f^*:\,A\to B$ of graded algebras we now construct a morphism $(f,\,f^*)$ of the locally ringed spaces 
($\PB$, ${\mathcal O}_{\textrm{\footnotesize Proj\,}(B)}$) and ($\PA$, ${\mathcal O}_{\textrm{\footnotesize Proj\,}(A)}$). For nongraded algebras it coincides with the one introduced in Proposition 2.3 (b) in Chapter II of \cite{Ha}.

Let $f^*:\,A\to B$ be a morphism of graded algebras. Let $f$ be the continuous map
\begin{eqnarray*}
  \begin{array}{cccc}
    f: &\PB & \to & \PA\\
         & Q  &\mapsto & (f^*)^{-1}(Q)^h
  \end{array}
\end{eqnarray*}
of Theorem \ref{morph}. It is easy to prove the following:
\begin{prop} We obtain a morphism of presheaves $\widetilde{f}^*: \,\widetilde{\mathcal O}_{\textrm{\footnotesize Proj\,}(A)}\to f_*(\widetilde{\mathcal O}_{\textrm{\footnotesize Proj\,}(B)})$ by 
\begin{eqnarray*}
   \begin{array}{cccc}
    \widetilde{f}^*(U): & \widetilde{{\mathcal O}}_{\textrm{\footnotesize Proj\,}(A)}(U) & \to & \widetilde{{\mathcal O}}_{\textrm{\footnotesize Proj\,}(B)}(f^{-1}(U))\\
         & \frac{a}{b}  &\mapsto & \frac{f^*(a)}{f^*(b)} 
  \end{array},
\end{eqnarray*}
where $U\subseteq \PA$ is open and nonempty. The morphisms induced on the stalks are given by the local morphisms
\begin{eqnarray*}
  \begin{array}{cccc}
    f^*_Q: & A_{(f(Q))} & \to & B_{(Q)}\\
         & \frac{a}{b}  &\to & \frac{f^*(a)}{f^*(b)} 
  \end{array},\qquad Q\in\PB,
\end{eqnarray*}
of Theorem \ref{morph}.
\end{prop}

Hence there exists uniquely  a morphism of sheaves $f^*:\, {\mathcal O}_{\textrm{\footnotesize Proj\,}(A)}\to f_*({\mathcal O}_{\textrm{\footnotesize Proj\,}(B)})$ such that the diagram
\begin{eqnarray*}
\xymatrix{
    \widetilde{\mathcal O}_{\textrm{\footnotesize Proj\,}(A)} \ar[r]^{\mbox{\footnotesize  $\widetilde{f}^*$\  \ }} \ar[d]_{\mbox{\footnotesize  $\tau_A$}} &  f_*(\widetilde{\mathcal O}_{\textrm{\footnotesize Proj\,}(B)})  \ar[d]^{\mbox{\footnotesize  $f_*(\tau_B)$}}  \\
      {\mathcal O}_{\textrm{\footnotesize Proj\,}(A)} \ar[r]^{\mbox{\footnotesize  $f^*$\  \ }}  &  f_*({\mathcal O}_{\textrm{\footnotesize Proj\,}(B)}) 
}
\end{eqnarray*} 
commutes, $\tau_A: \widetilde{\mathcal O}_{\textrm{\footnotesize Proj\,}(A)}\to {\mathcal O}_{\textrm{\footnotesize Proj\,}(A)} $ and $\tau_B:\widetilde{\mathcal O}_{\textrm{\footnotesize Proj\,}(B)}\to {\mathcal O}_{\textrm{\footnotesize Proj\,}(B)} $ the canonical morphisms of the sheafifications. 
The local morphisms  $f^*_Q:\,  A_{(f(Q))} \to  B_{(Q)}$, $Q\in \PB$, are the induced morphisms on the stalks of the structure sheaves.\vspace*{1ex}

If we assign to a graded algebra its associated locally ringed space and to a morphism of graded algebras its associated morphism of locally ringed spaces we get a contravariant functor from the category of graded algebras into the category of locally ringed spaces spaces. \vspace*{1ex}

The next two results are immediate corollaries of Theorem \ref{loc} and Theorem \ref{prodlattice}.
\begin{cor}\label{principalscheme} Let $A$ be a  $\mathcal M$-graded algebra. Let $a\in A$ be homogeneous, not nilpotent. The canonical morphism $i^*:A\to \{a^n|n\in\Nn\}^{-1}A$ induces an isomorphism of the restriction ($D(a)$, ${\mathcal O}_{\textrm{\footnotesize Proj\,}(A)}\res{D(a)}$) and ($\mb{Proj\,}(\{a^n|n\in\Nn\}^{-1}A)$, ${\mathcal O}_{\textrm{\footnotesize Proj\,}(\{a^n|n\in\Nn\}^{-1}A)}$). 
\end{cor}
\begin{cor}\label{prodlatticescheme} Let $A$ be a $\mathcal M$-graded algebra. Let $\mathcal N$ be a torsion-free commutative group, and let $\FK{\mathcal N}$ be its $\mathcal N$-graded group algebra. The canonical morphism $i_A^*:A\to A\otimes\FK{\mathcal N}$ induces an isomorphism of ($\PA$, ${\mathcal  O}_{\textrm{\footnotesize Proj\,}(A)}$) and ($\mb{Proj\,}(A\otimes \FK{\mathcal N})$, ${\mathcal O}_{\textrm{\footnotesize Proj\,}(A\otimes \F [{\mathcal N}])}$).
\end{cor}
%
%

%
%
%
%
%
\section{The projective spectrum of the Cartan algebra and its $\F$-valued points\label{psCa}}
In this section we describe the projective spectrum $\PCA$ of the $P^+$-graded Cartan algebra $CA$ associated to a symmetrizable Kac-Moody group $G$. We show that its $\F$-valued points $\PCAF$ identify with a completion $\Omega_{fn}$ of the building $\Omega$ associated to $G$.
\subsection{The Cartan algebra\label{subsectionCA}} 
The Cartan algebra $CA$ associated to the Kac-Moody group $G$ is a commutative associative unital $P^+$-graded algebra over the field $\F$ without zero-divisors. It is obtained by the following construction: 
Take as $P^+$-graded $\F$-linear space
\begin{eqnarray*}
     CA:=\bigoplus_{\La\in P^+}L(\La)^{(*)},
\end{eqnarray*}
where $ L(\La)^{(*)}:=\bigoplus_{\la\in P(\La)} L(\La)_\la^*\subseteq L(\La)^*$ is the restricted $\F$-linear dual of $L(\La)$.
To define the product of $CA$ fix for every $\La\in P^+$ a highest weight vector $v_\La\in L(\La)_\La\setminus\{0\}$. 
For $\La,\,N\in P^+$ the $G$-module $L(\La)\otimes L(N)$ decomposes in the form
\begin{eqnarray*}
   L(\La)\otimes L(N)= \underbrace{L_{ \textrm{\footnotesize high} }}_{\cong L(\La + N)}  \oplus \underbrace{L_{  \textrm{\footnotesize low} }}_{ \textrm{\scriptsize only isotypical components} \atop  \textrm{\scriptsize of type }\,L(M),\,M<\La+N}.
\end{eqnarray*}
Since $v_\La\otimes v_N$ is a highest weight vector of $L_{ \textrm{\footnotesize high} }$, it is possible to define a $G$-equivariant linear map
\begin{eqnarray*}
  \Phi:\,L(\La+N)\to L(\La)\otimes L(N) \quad\mb{ by }\quad \Phi(v_{\La+N}):=v_\La\otimes v_N.
\end{eqnarray*} 
The product $\prdca$ of the Cartan algebra $CA$ between the parts $L(\La)^{(*)}$, $L(N)^{(*)}$ is now obtained dually to $\Phi$:
\begin{eqnarray*}
    \prdca : L(\La)^{(*)}\times L(N)^{(*)}\to L(\La)^{(*)}\otimes L(N)^{(*)}\to\left(L(\La)\otimes L(N)\right)^{(*)}\stackrel{\Phi^{(*)}}{\to} L(\La+N)^{(*)}.
\end{eqnarray*}
The unit of the Cartan algebra $CA$ is given by the element $\delta_0\in L(0)_0^*\subseteq CA$ which satisfies $\delta_0(v_0)=1$.

Algebras obtained in this way by different choices of highest weight vectors in $L(\La)_\La\setminus\{0\}$, $\La\in P^+$, are isomorphic.

The Borel-Weil map $BW: CA\to\F[G]^U$ of Section \ref{Preli} is a $G$-equivariant isomorphism of algebras.\vspace*{1ex}

\subsection{Some actions on the projective spectrum of the Cartan algebra} 
Fix $\La\in P^+$. The face monoid $\GD$ acts by its definition on the highest weight module $L(\La)$. This action induces an action of the opposite monoid $\widehat{G}^{\,op}$ on $L(\La)^{(*)}$. 

The $\g$-module $L(\La)^{(*)}$ is a lowest weight module of lowest weight $-\La$. The formal Kac-Moody group $G_{fn}$ acts by its definition on $L(\La)^{(*)}$. Dually, we get an action of $G_{fn}$ on 
\begin{eqnarray*}
   L(\La)_f := \prod_{\la\in P(\La)} L(\La)_\la = \prod_{\la\in P(\La)} (L(\La)_\la^*)^* = (L(\La)^{(*)})^*  \,,
\end{eqnarray*} 
where $L(\La)_\la$ has been identified canonically with $(L(\La)_\la^*)^*$, $\la\in P(\La)$. 

The Kac-Moody group $G$ sits inside $\GD$ and $G_{fn}$. The actions introduced before restrict to the same action of $G$ on $L(\La)$. They restrict to  actions of $G^{op}$ and $G$ on $L(\La)^{(*)}$, which differ by the inverse map. This is cumbersome. Therefore, in this article we work with the action $(G_{fn})^{op}$ on $L(\La)^{(*)}$ obtained from the action of $G_{fn}$ on $L(\La)^{(*)}$ by concatenation the corresponding representation with the inverse map. Similarly, we work with the action of $(\g_{fn})^{op}\supseteq \g^{op}$ on $L(\La)^{(*)}$ obtained from the action of $\g_{fn}\supseteq \g$ on $L(\La)^{(*)}$ by concatenation the corresponding representation with the multiplication by $-1$. 

The actions of $\widehat{G}^{\,op}$ and $(G_{fn})^{op}$ on all $L(\La)^{(*)}$, $\La\in P^+$, combine to actions of $\widehat{G}^{\,op}$ and $(G_{fn})^{op}$ on the Cartan algebra $CA$ by morphisms of graded algebras. 

To keep our notation clear and simple we denote the representations of the opposite monoid $\widehat{G}^{\,op}$, the opposite group  $(G_{fn})^{op}$ , the opposite Lie algebra $(\g_{fn})^{op}$ by the same symbol $\pi$.

Now, for $x\in \GD$ resp. $x \in G_{fn}$ we get a continuous map on the spectrum $\PCA$ of all $P^+$-homo\-ge\-neous prime ideals of $CA$ by 
\begin{eqnarray*}
   x Q:= \pi(x)^{-1}(Q) =\bigoplus_{\La\in P^+}\pi(x)^{-1}(Q_\La)\quad\mb{ where }\quad Q=\bigoplus_{\La\in P^+}Q_\La\in\PCA.
\end{eqnarray*}
It leaves the spectrum of $\F$-valued points $\PCAF$ invariant. Therefore, we get an action of $\GD$ as well as an action of $G_{fn}$ on the spectra $\PCA$ and $\PCAF$.
(We could have combined the actions of $\widehat{G}$ and $G_{fn}$ to the action of a bigger monoid. For the aims of this article it is not advantageous.)
\subsection{The completed building} 

The Kac-Moody group $G$ has the opposite BN-pairs $(B^\pm,\,N)$ with corresponding standard parabolic subgroups $P_J^\pm$, $J\subseteq I$.
 
The formal Kac-Moody group $G_{fn}$ has the BN-pair $(B_f^-,\,N)$ with corresponding standard parabolic subgroups $(P^-_{fn})_J$, $J\subseteq I$. 
The group  $(P^-_{fn})_J$ is an extension of $P_J^-$, i.e.,  $(P^-_{fn})_J\cap G= P_J^-$. Furthermore, for every $\La\in P^+\cap F_J$ it holds
\begin{eqnarray*}
  N_{G_{fn}}(L(\La)_\La^*) := \Mklz{g\in G_{fn}}{\pi(g) L(\La)_\La^* = L(\La)_\La^*}=(P^-_{fn})_J.
\end{eqnarray*}
Now we similarly extend the parabolic subgroups $P_J$, $J\subseteq I$, of $G$ to subgroups of $G_{fn}$. For $J\subseteq I$ set
\begin{eqnarray*}
  (P_{fn})_J:=  (U_f^-)_J P_J = (U_f^-)_J (\We_J T) U.
\end{eqnarray*}
In particular, $(P_{fn})_\emptyset=B$ and $(P_{fn})_I=G_{fn}$.

\begin{thm}\label{ngfnllala} Let $J\subseteq I$.
\begin{itemize}
\item[(a)] $(P_{fn})_J$ is a subgroup of $G_{fn}$ such that $(P_{fn})_J\cap G=P_J$.
\item[(b)] For every $\La\in P^+\cap F_J$ we have
\begin{eqnarray*}
  N_{G_{fn}}(L(\La)_\La):=\Mklz{g\in G_{fn}}{ g L(\La)_\La=  L(\La)_\La}=(P_{fn})_J .
\end{eqnarray*}
\end{itemize}
\end{thm}
\Proof We first show (b). We have
\begin{eqnarray*}
 && N_{\bf g}(L(\La)_\La):=\Mklz{x\in \g}{ x L(\La)_\La\subseteq  L(\La)_\La} = \n_J^-\oplus \h \oplus \n\,,\\
 && N_{G}(L(\La)_\La):=\Mklz{g\in G}{g L(\La)_\La=L(\La)_\La}= P_J  \, . 
\end{eqnarray*}
In particular, $\n_J^-L(\La)_\La\subseteq L(\La)_\La$, which implies $(\n_f^-)_J L(\La)_\La\subseteq L(\La)_\La$, which in turn implies $(U_f^-)_J L(\La)_\La = L(\La)_\La$. Therefore, we find $(P_{fn})_J=(U_f^-)_J P_J\subseteq N_{G_{fn}}(L(\La)_\La)$.

Now let $g\in N_{G_{fn}}(L(\La)_\La)$. Decompose $g$ in the form $g=u^J u_J n_\sigma v$ with $u^J\in (U_f^-)^J$, $u_J\in
(U_f^-)_J$, $n_\sigma\in N$ projecting to $\sigma\in\We$, and $v\in U^+$. From
\begin{eqnarray*}
   L(\La)_\La = g L(\La)_\La = u^J u_J n_\sigma L(\La)_\La = u^J u_J  L(\La)_{\sigma\La}
\end{eqnarray*}
we obtain $\sigma\La=\La$, which is equivalent to $\sigma\in\We_J$, and $u^J L(\La)_\La = L(\La)_\La$. Now $u^J$ acts as
\begin{eqnarray*} 
 u^J=exp(x) & \mb{ for some }& x=\sum_{\al\in(\Delta^-)^J}x_\al\in\prod_{\al\in(\Delta^-)^J}\g_\al.
\end{eqnarray*}
Assume that $u^J\neq 1$. Choose a nonzero homogeneous component $x_\beta$ with $\beta$ of maximal height. Then $x_\beta v_\La$ is the $(\La+\beta)$-homogeneous part of $u^J v_\La$. It is nonzero because of $x_\beta\notin {\bf p}_J = N_{\g}(L(\La)_\La)$. This contradicts $u^J L(\La)_\La=L(\La)_\La$. Therefore, $u^J=1$ and $g=u_J n_\sigma v\in (U_f^-)_J  (\We_J T) U  = (P_{fn})_J$.

Now (a) follows from (b) and $N_{G_{fn}}(L(\La)_\La)\cap G=N_G(L(\La)_\La)=P_J$.
\qed

The following easy Lemma will be used several times.
\begin{lem}\label{prodint} Let $C$ be a subgroup of $U_f^-$, let $D$ be a subgroup of $N$, and let $E$ be a subgroup of $U$. Then
\begin{eqnarray*}
   CDE\cap  U_f^- = C\quad\mb{ and }\quad CDE\cap N = D \quad\mb{ and }\quad CDE\cap U = E.
\end{eqnarray*}
\end{lem}
\Proof The inclusions "$\supseteq$" hold, because the unit 1 is contained in $C$, $D$, and $E$. The inclusion "$\subseteq $" of the second equality follows directly from the Birkhoff decomposition of $G_{fn}$. The inclusion "$\subseteq $" of the first and third equality  are shown similarly; let us illustrate this for the first equality. Let $c\in C$, $d\in D$, $e\in E$, and $u\in U_f^-$ such that $cde=u$. Then  
$(u^{-1}c)de=1$, from which we obtain $u^{-1}c=d=e=1$. In particular, $u=c\in C$.
\qed

\begin{prop} Let $J\subseteq I$. Then $(P_{fn})_J=P_J$ if and only if $J=J^0$.
\end{prop}
\Proof It is $(U_f^-)_J=U^-_J$ if and only if $J=J^0$. Now the proposition follows from the definition of $(P_{fn})_J$ and Lemma \ref{prodint}.\qed

Let $J\subseteq I$. The parabolic subgroup $(P_{fn}^-)_J$ has the Levi decomposition
\begin{eqnarray*}
     (P_{fn}^-)_J =  (U^-_f)^J \rtimes   (L_{fn})_J
  \end{eqnarray*}
with $(L_{fn})_J :=  (U^-_f)_J (\We_J T) (U^-_f)_J = (U^-_f)_J (\We_J T) U_J^- = (U^-_f)_J (\We_J T) U_J $. Similarly, $(P_{fn})_J$ has a Levi decomposition:
\begin{thm}\label{levipfn} For $J\subseteq I$ it is $ (P_{fn})_J=(L_{fn})_J \ltimes U^J$.
\end{thm}
\Proof This holds trivially for $J=I$. Let $J\neq I$. $(L_{fn})_J$ and $U^J$ are subgroups of $G_{fn}$ such that 
\begin{eqnarray*}
  (P_{fn})_J = (U_f^-)_J (\We_J T) U = (U_f^-)_J  (\We_J T) U_J U^J =(L_{fn})_J U^J. 
\end{eqnarray*}
By Lemma \ref{prodint}  we obtain $(L_{fn})_J\cap U^J = (L_{fn})_J\cap U \cap U^J = U_J\cap U^J=\{1\}$.

It remains to show that $(L_{fn})_J  = (U^-_f)_J (\We_J T) (U^-_f)_J$ normalizes $U^J$. The group $\We_J T$ normalizes $U^J=\bigcap_{\sigma\in {\cal W}_J} \sigma U\sigma^{-1}$. The proof that $(U^-_f)_J$ normalizes $U^J$ is divided into several steps.

(a) We first show 
\begin{eqnarray}\label{Uhgw}
   U = \Mklz{ g\in G_{fn} }{ gv_\La=v_\La \mb{ for all }\La\in P^+\cap C }.
\end{eqnarray}
Trivially, the inclusion "$\subseteq $" holds. To show "$\supseteq $" fix $\La\in P^+\cap C$ and let $g\in G_{fn}$ such that $g v_\La = v_\La$. Then $g L(\La)_\La = L(\La)_\La$, from which we obtain $g\in B = U T$  by Theorem \ref{ngfnllala} (b). Write $g$ in the form
\begin{eqnarray*}
  g = u\prod_{i=1}^{2n-l}t_i(s_i)  \quad\mb{ with }\quad u\in U  \mb{ and }  s_1,\,s_2,\,\ldots,\,s_{2n-l}\in\F^\times \,.
\end{eqnarray*} 
From $g v_\La= v_\La$ we get $\prod_{i=1}^{2n-l}(s_i)^{\La(h_i)} =1$. For $k\in \{1,\,2,\,\ldots,\,2n-l\}$ it is $\La+\La_k \in P^+\cap C$. Suppose that in addition $g v_{\La+\La_k}=v_{\La+\La_k}$. Then also $\prod_{i=1}^{2n-l}(s_i)^{(\La+\La_k)(h_i)}=s_k\prod_{i=1}^{2n-l}(s_i)^{\La(h_i)}=1$. Hence $s_k=1$.

For every $\sigma\in \We_J$ choose $n_\sigma\in N$ projecting to $\sigma$. By (\ref{Uhgw}) we find
\begin{eqnarray}\label{UJhgw}
   U^J  =\bigcap_{\sigma\in {\cal W}_J} \sigma U\sigma^{-1} = \Mklz{ g\in G_{fn}\, }{\, g\, n_\sigma v_\La = n_\sigma v_\La \mb{ for all }\sigma\in\We_J \mb{ and all }\La\in P^+\cap C }.
\end{eqnarray}

(b) It is $(U^-_f)_J =exp((\n_f^-)_J)$ with $(\n_f^-)_J=\bigoplus_{\al\in \Delta_J^-}\g_\al$. For $v_\la\in L(\La)_\la$, $\la\in P$, we get
\begin{eqnarray*}
  (U^-_f)_J \,v_\la \subseteq \prod_{\mu\,\in\, \la - (Q_0^+)_J} L(\La)_\mu\,.
\end{eqnarray*}

The group $U^J$ is also a subgroup of the group $G_{fp}$, which acts on $L(\La)$. In particular, we have $U^J\subseteq (U^+_f)^J =exp((\n_f)^J)$ with $(\n_f)^J=\bigoplus_{\al\in\Delta^+\setminus \Delta_J^+}\g_\al$. For $v_\la\in L(\La)_\la$, $\la\in P$, we find
\begin{eqnarray*}
  U^J v_\la \subseteq v_\la+ \bigoplus_{\mu\,\in\, \la + Q_0^+\setminus (Q_0^+)_J} L(\La)_\mu\,.
\end{eqnarray*}

(c) Now let $u\in U^J$ and  $u^-\in (U_f^-)_J$. We use (\ref{UJhgw}) to show $u^- u (u^-)^{-1}\in U^J$. Let $\sigma\in\We_J$ and $\La\in P^+\cap C $. By (b) we obtain
\begin{eqnarray*}
  u (u^-)^{-1}n_\sigma v_\La \;\in\; (u^-)^{-1}n_\sigma v_\La + \prod_{\mu\,\in \,\sigma\La- (Q^+_0)_J + Q^+_0\setminus(Q^+_0)_J} L(\La)_\mu \,.
\end{eqnarray*}
By Theorem 3.12 d) of \cite{K} applied to the fundamental chamber in $\h_\R^*$ it is $\La-\sigma \La\in Q_0^+$.  Because of $\sigma\in \We_J$ we find $\La-\sigma \La\in (Q_0^+)_J$. Therefore, we get
\begin{eqnarray*}
    \sigma\La- (Q^+_0)_J + Q^+_0\setminus(Q^+_0)_J \,\subseteq\, \La- (Q^+_0)_J + Q^+_0\setminus(Q^+_0)_J.
\end{eqnarray*}
An element $\mu$ of the set on the right is of the form
\begin{eqnarray*}
  \mu = \La - \sum_{i\in J} n_i\al_i + \sum_{i\in I} m_i\al_i  = \La- \sum_{i\in J} (n_i-m_i)\al_i  +\sum_{i\in I\setminus J} m_i\al_i 
\end{eqnarray*}
with $n_i\in \Nn$, $i\in J$, with $m_i\in\Nn$, $i\in I$, and $m_i\neq 0$ for at least one $i\in I\setminus J$. Therefore, $\mu\not\in \La- Q^+_0\supseteq P(\La)$, from which we conclude that $L(\La)_\mu=\{0\}$. This shows $u^- u (u^-)^{-1}n_\sigma v_\La = u^- ( (u^-)^{-1}n_\sigma v_\La + 0) = n_\sigma v_\La $.
\qed

There are the following modified generalized Birkhoff decompositions:
\begin{thm} Let $J\subseteq I$. Then
\begin{eqnarray*}
  G_{fn} = \dot{\bigcup_{\sigma\in {\cal W}^J}}  U_f^- \,\sigma \,(P_{fn})_J = \dot{\bigcup_{\sigma\in {\cal W}^J}} \big( U_f^-\cap \sigma (U_f^-)^J\sigma^{-1} \big) \,\sigma \, (P_{fn})_J \,.
\end{eqnarray*}
Furthermore, for $\sigma\in \We^J$ and $u_1,\,u_2\in U_f^-\cap \sigma (U_f^-)^J\sigma^{-1}$ we have
\begin{eqnarray*}
   u_1 \sigma (P_{fn})_J = u_2 \sigma (P_{fn})_J \quad\iff\quad u_1=u_2\,.
\end{eqnarray*}
\end{thm}
\Proof (a) It is $G_{fn}= G_{fn} (P_{fn})_J  = U_f^- G (P_{fn})_J $. By the generalized Birkhoff decomposition of $G$ we get
\begin{eqnarray*}
  G_{fn} = U_f^- \Big( \dot{\bigcup_{\sigma\in {\cal W}^J}} U^- \,\sigma \,P_J \Big) (P_{fn})_J 
       = \bigcup_{\sigma\in {\cal W}^J} U_f^- \,\sigma \,(P_{fn})_J \,.
\end{eqnarray*}
This union is also disjoint. Suppose that $U_f^- \,\sigma \,(P_{fn})_J \,\cap\, U_f^- \,\sigma' \,(P_{fn})_J\neq\emptyset $. Then there exist $u,\,u'\in U_f^-$, $p,\,p'\in (P_{fn})_J$, and $n_\sigma,\,n'_{\sigma'}\in N$ projecting to $\sigma,\,\sigma'$ such that $u n_\sigma p  = u'n'_{\sigma'} p'$. Apply this element to $L(\La)_\La$, where $\La\in P^+\cap F_J$ is chosen arbitrarily. By Theorem \ref{ngfnllala} (b) we obtain $u L(\La)_{\sigma\La} = u' L(\La)_{\sigma'\La}$, from which in turn we find $\sigma\La=\sigma'\La$, which is equivalent to $\sigma\We_J=\sigma'\We_J$. Because $\sigma$, $\sigma'$ are minimal coset representatives, we get $\sigma=\sigma'$.

(b) Let $\sigma\in\We^J$. It holds $U_f^- = (U_f^-\cap\sigma (U_f^-)\sigma^{-1})(U_f^- \cap\sigma U\sigma^{-1})$, which is a property of refined Tits systems. Therefore, we find 
\begin{eqnarray}
  && U_f^- \,\sigma \,(P_{fn})_J = (U_f^-\cap\sigma (U_f^-)\sigma^{-1})(U_f^- \cap\sigma U\sigma^{-1}) \sigma (P_{fn})_J \nonumber\\
  &&   = (U_f^-\cap\sigma (U_f^-)\sigma^{-1})\sigma (\sigma^{-1}U_f^- \sigma\cap U)  (P_{fn})_J 
     =  (U_f^-\cap\sigma (U_f^-)\sigma^{-1})\sigma  (P_{fn})_J \,. \label{birktr0}
\end{eqnarray}

Next we show
\begin{eqnarray}\label{birktr1}
  U_f^- \cap \sigma (U^-_f)\sigma^{-1} = \big(U_f^-\cap\sigma (U_f^-)^J \sigma^{-1}\big)\sigma (U^-_f)_J\sigma^{-1} \,.
\end{eqnarray}
Since $\sigma$ has minimal length in $\sigma\We_J$ we have $l(\sigma\sigma_j)>l(\sigma)$ for all $j\in J$,  which implies $\sigma\Delta_J^-\subseteq \Delta^-$ by Lemma 3.11 (a) in \cite{K}. Choose $n_\sigma\in N$ projecting to $\sigma$. Then 
\begin{eqnarray*}
  Ad(n_\sigma)(\n_f^-)_J = Ad(n_\sigma)\prod_{\al\in \Delta_J^-}\g_\al = \prod_{\al\in \Delta_J^-}\g_{\sigma\al}\subseteq \n_f^-\,,
\end{eqnarray*}
from which we find
\begin{eqnarray}\label{birktr2}
  \sigma (U^-_f)_J\sigma^{-1} = n_\sigma\,exp((\n_f^-)_J)\,n_\sigma^{-1} = exp(Ad(n_\sigma) (\n_f^-)_J)\subseteq U_f^-\,.
\end{eqnarray}
This shows the inclusion "$\supseteq $" of (\ref{birktr1}). To show the reverse inclusion let $a=bc$ with $a\in U_f^-$ and $b\in \sigma (U^-_f)^J\sigma^{-1}$, $c\in \sigma (U^-_f)_J\sigma^{-1}$. By  (\ref{birktr2}) we get $b=a c^{-1}\in U_f^-\cap\sigma (U^-_f)^J\sigma^{-1} $.

Inserting  (\ref{birktr1}) in (\ref{birktr0}) we find
\begin{eqnarray*}
    U_f^- \,\sigma \,(P_{fn})_J   
         = \big(U_f^-\cap\sigma (U_f^-)^J \sigma^{-1}\big)\sigma (U^-_f)_J (P_{fn})_J 
         =  \big(U_f^-\cap\sigma (U_f^-)^J \sigma^{-1}\big)\sigma (P_{fn})_J .
\end{eqnarray*}

(c) Let $\sigma\in \We^J$ and $u_1,\,u_2\in U_f^-\cap \sigma (U_f^-)^J\sigma^{-1}$ that $ u_1 \sigma (P_{fn})_J = u_2 \sigma (P_{fn})_J $. By Lemma \ref{prodint} we obtain
\begin{eqnarray*}
  u_1^{-1}u_2 \,\in\, \sigma  (P_{fn})_J \sigma^{-1}\cap (U_f^-\cap \sigma (U_f^-)^J\sigma^{-1})\subseteq \sigma \big( (P_{fn})_J \cap (U_f^-)^J \big)\sigma^{-1}= \sigma \big( (U_f^-)_J\cap (U_f^-)^J \big)\sigma^{-1} =\{1\}.
\end{eqnarray*}
\qed

The next two theorems show: The system of nonparabolic subgroups $(P_J)_{fn}$, $J\subseteq I$, of $G_{fn}$ has many properties in common with systems of parabolic subgroups. 

\begin{thm} \label{spfpp} Let $J$, $K\subseteq I$. The following hold:
\begin{itemize}
\item[(a)] $(P_{fn})_J\subseteq (P_{fn})_K$ if and only if $J\subseteq K$.
\item[(b)] $(P_{fn})_J\cap (P_{fn})_K = (P_{fn})_{J\cap K}$. 
\item[(c)] In general, $(P_{fn})_{J\cup K}$ is not generated by $(P_{fn})_J$ and $(P_{fn})_K$ as a group. 
\end{itemize}
\end{thm}
\Proof (a) For $J\subseteq K$ we have $(U_f^-)_J\subseteq (U_f^-)_K$, from which we get $(P_{fn})_J\subseteq (P_{fn})_K$. If $(P_{fn})_J\subseteq (P_{fn})_K$ we find $P_J\subseteq P_K$ by Theorem \ref{ngfnllala} (a). Hence $J\subseteq K$.

The inclusion "$\supseteq$" of part (b) follows from part (a).  Now let $g\in (P_{fn})_J\cap (P_{fn})_K $. By the Birkhoff decomposition of $G_{fn}$, and by the definitions of $(P_{fn})_J$, $(P_{fn})_K$, the element $g$ can be written in the forms
\begin{eqnarray*} 
   g = u_J n_\sigma u = \ti{u}_K n_\sigma \ti{u}
\end{eqnarray*}
with $u,\,\ti{u}\in U$, $n_\sigma\in N$ projecting to $\sigma \in \We_J\cap\We_ K = \We_{J\cap K}$ and $u_J\in (U_f^-)_J$, $\ti{u}_K\in (U_f^-)_K$. We find
\begin{eqnarray*}
  (u_J)^{-1} \ti{u}_K  \in U_f^-\cap \sigma U\sigma^{-1} = U^-\cap \sigma U\sigma^{-1}=\prod_{\al\in \Delta^-_{re}\cap \sigma\Delta^+_{re}} U_\al.
\end{eqnarray*}
Because of $\sigma\in\We_{J\cap K}$ we have $\Delta^-_{re}\cap \sigma\Delta^+_{re}\subseteq (\Delta_{J\cap K})^-_{re}$. Therefore, we get
\begin{eqnarray*}
   \ti{u}_K\in (U_f^-)_K\cap (U_f^-)_J U_{J\cap K}^-  = (U_f^-)_K\cap (U_f^-)_J = exp((\n_f^-)_K)\cap exp((\n_f^-)_J). 
\end{eqnarray*}
Since the exponential function $exp:\,\n^-_f\to U_f^-$ is bijective we obtain
\begin{eqnarray*}
   exp((\n_f^-)_K)\cap exp((\n_f^-)_J) = exp((\n_f^-)_K\cap (\n_f^-)_J) = exp((\n_f^-)_{J\cap K}) = (U_f^-)_{J\cap K}.
\end{eqnarray*} 
We conclude that $g=\ti{u}_K n_\sigma \ti{u}\in (U_f^-)_{J\cap K} (\We_{J\cap K}T) U=(P_{fn})_{J\cap K}$.

To show (c) suppose that for all $J,\, K\subseteq I$ the group $(P_{fn})_{J\cup K}$ is generated by the groups $(P_{fn})_J$ and $(P_{fn})_K$. Then $G_{fn}=(P_{fn})_I$ is generated by the groups $(P_{fn})_{\{i\}}=P_{\{i\}}$, $i\in I$. Hence $G_{fn}=G$, which is only possible if all components of the generalized Cartan matrix $A$ are of finite type.
\qed

\begin{thm} Let $J,\,K\subseteq I$ and $g\in G_{fn}$. The following are equivalent:
\begin{itemize}
\item[(i)] $g(P_{fn})_J g^{-1}\subseteq (P_{fn})_K$.
\item[(ii)] $J\subseteq K$ and $g\in (P_{fn})_K$.
\end{itemize}
In particular: The subgroups $(P_{fn})_J$, $(P_{fn})_K$ of $G_{fn}$ are conjugated if and only if $J=K$. The subgroup $(P_{fn})_J$ of $G_{fn}$ coincides with its normalizer. 
\end{thm}
\Proof It only remains to show that $(i)$ implies $(ii)$. Here a well-known argument to compute  representation theoretically the normalizers of parabolics also applies:

(a) Fix $\La\in P^+$. Suppose that $V$ is a one-dimensional subspace of $L(\La)_f$ such that $B\,V\subseteq V$. We show $V=L(\La)_\La$:
 
Let $v\in V\setminus\{0\}$. Decompose
\begin{eqnarray*}
  v=\sum_{\la\in  \textrm{\footnotesize  supp} (v)} v_\la \quad \mb{ with } \quad \mb{supp}(v):=\Mklz{\la\in P(\La)}{v_\la\neq 0}\neq\emptyset
\end{eqnarray*}
according to the decomposition $L(\La)_f=\prod_{\la\in P(\La)}L(\La)_\la$.

We first show $|\mb{supp}(v)|=1$. Suppose there exist $\la_1,\,\la_2\in\mb{supp}(v)$, $\la_1\neq\la_2$. Choose $h\in H$ such that $\la_1(h)\neq \la_2(h)$. Choose $s\in\N\subseteq \F^\times$ such that $s^{\la_1(h)}\neq s^{\la_2(h)}$. We get
\begin{eqnarray*}
  t_h(s)v = \sum_{\la\in \textrm{\footnotesize supp}(v)}s^{\la(h)}v_\la \,\notin\, \F v\,,
\end{eqnarray*}
which contradicts  $T\,V\subseteq V$.

Now let $v = v_\la$. Suppose that $\n\, v_\la \neq \{ 0\}$. Since $\n$ is generated by $e_i$, $i\in I$, there exists an index $i$ such that $e_i v_\la \neq 0$. We find
\begin{eqnarray*}
  exp(e_i) v_\la = \underbrace{v_\la}_{\in L(\La)_\la \setminus \{0\}}  
                       + \underbrace{e_i v_\la}_ {\in L(\La)_{\la +\al_i}\setminus \{0\}}  
                        +\sum_{k\in \Nn,\,k\geq 2}\underbrace{\frac{1}{k!}(e_i)^k v_\la}_{\in L(\La)_{\la+ k\al_i} }  
         \;\notin \; \F v_\la,
\end{eqnarray*}
which contradicts $U\,V\subseteq V$.

We have shown $v\in L(\La)_\la$ for some $\la\in P(\La)$ and $\n\, v=\{0\}$. By Proposition 9.3 b) in \cite{K} we get $v \in L(\La)_\La$.

(b) Now suppose that (i) holds. Choose $\La\in P^+\cap F_K$. By (i) and Theorem \ref{ngfnllala} (b) we find
\begin{eqnarray*}
   B g^{-1} L(\La)_\La \subseteq (P_{fn})_J g^{-1} L(\La)_\La  \subseteq g^{-1}(P_{fn})_K L(\La)_\La = g^{-1}L(\La)_\La .
\end{eqnarray*}
Now $g^{-1}L(\La)_\La$ is a one-dimensional subspace of $L(\La)_f$. By (a) we obtain $g^{-1}L(\La)_\La = L(\La)_\La$. Again by Theorem 
\ref{ngfnllala} (b) we find $g\in (P_{fn})_K$. With (i) we get $(P_{fn})_J \subseteq (P_{fn})_K$, which is equivalent to $J\subseteq K$  by Theorem (\ref{spfpp}) (a).
\qed

We take as {\it completed building} $\Omega_{fn}$ associated to the formal Kac-Moody group $G_{fn}$ the set
\begin{eqnarray*}
    \Omega_{fn} := \dot{\bigcup_{J\subseteq I}}G_{fn}/(P_{fn})_J =\Mklz{g (P_{fn})_J}{g\in G_{fn},\;J\subseteq I}
\end{eqnarray*}
partially ordered by the reverse inclusion, i.e., for $g,\,g'\in G_{fn}$ and $J,\,J'\subseteq I$, 
\begin{eqnarray*}
    g (P_{fn})_J\leq g' (P_{fn})_{J'} \;:\iff\; g (P_{fn})_J\supseteq g' (P_{fn})_{J'}. 
\end{eqnarray*}
The group $G_{fn}$ acts order preservingly on $\Omega_{fn}$ by multiplication from the left. 
We denote by
\begin{eqnarray*}
   {\mathcal A}_{fn}:=\Mklz{n (P_{fn})_J}{n\in N,\,J\subseteq I}
\end{eqnarray*} 
the {\it standard apartment} of $\Omega_{fn}$. The completed building $\Omega_{fn}$ is covered by the apartments $g{\mathcal A}_{fn}$, $g\in G_{fn}$.
\begin{rem} (a) At first, the construction of $\Omega_{fn}$ by the nonparabolic subgroups $(P_{fn})_J$, $J\subseteq I$, of $G_{fn}$ may look strange. We give a hand-waving motivation, which indicates that it is natural: 

Let $J\subseteq I$. Consider first the classical case, i.e., all components of the generalized Cartan matrix $A$ are of finite type, and take $\F=\C$. Here $G/P_J$ is a well-known manifold, a partial flag manifold. It is covered by big cells, which are charts of the manifold. The standard big cell $BC(J)$ of $G/P_J$ is obtained as the image of $U^-$ in $G/P_J$. As a set it is given by
\begin{eqnarray*}
 BC(J)\cong U^-/(U^-\cap P_J)=U^-/U_J^-\cong (U^-)^J.
\end{eqnarray*}

Now consider the case of an arbitrary generalized Cartan matrix. By analogy, the standard big cell of $G_{fn}/(P_{fn})_J$ 
is obtained as the image of $U^-_f$ in $G_{fn}/(P_{fn})_J$. As a set it is given by
\begin{eqnarray*}
 BC(J)\cong U_f^-/(U_f^-\cap (P_{fn})_J)=U^-_f/(U_f^-)_J\cong (U_f^-)^J.
\end{eqnarray*}

(b) The completed building $\Omega_{fn}$ should be considered as an equivariant completion of the building $\Omega$ associated to $G$ as described in Section \ref{afmb}. In the nonclassical case, i.e., if not all components of the generalized Cartan matrix $A$ are of finite type, it does not satisfy all properties of a building. This can be seen as follows: Let $g_1,\,g_2\in G_{fn}$. It is easy to check that the following are equivalent.
\begin{itemize}
\item[(i)] $g_1 B$ and $g_2 B$ are contained in a common apartment.
\item[(ii)] $(g_1)^{-1}g_2\in G$.
\end{itemize} 
In the nonclassical case there exists an element $g\in G_{fn}\setminus G$. Thus, $B$ and $gB$ are not contained in a common apartment.
\end{rem}
\subsection{The projective spectrum of the Cartan algebra and its $\mathbb{F}$-valued points}
We first show that the completed building $\Omega_{fn}$ embeds $G_f$-equivariantly into the $\F$-valued points $\PCAF$ of the spectrum $\PCA$.

For $\La\in P^+$ let $\delta_\La\in L(\La)^*_\La\subseteq CA$ be defined by $\delta_\La(v_\La):=1$. It is easy to check that $\delta_\La  \prdca \delta_N=\delta_{\La+N}$ for all $\La$, $N\in P^+$. Therefore,  
\begin{eqnarray*}
   \widetilde{P^+}:=\Mklz{\delta_\La}{\La\in P^+}
\end{eqnarray*} 
is a multiplicatively closed subset of the Cartan algebra $CA$, isomorphic to $(P^+,\,+)$. For $M\subseteq P^+$ we set
\begin{eqnarray*}
 \widetilde{M}:=\Mklz{\delta_\La}{\La\in M}\subseteq \widetilde{P^+}.
\end{eqnarray*}

For $\La\in P^+$ we set
\begin{eqnarray*}
     L(\La)^{(*)}_{\neq\La}:=\Mklz{\phi\in L(\La)^{(*)}}{\phi(v_\La)= 0}=\bigoplus_{\la\in P(\La)\setminus\{\La\}}L(\La)_\la^* \,.
\end{eqnarray*}
Here, as always, a sum over the empty set is defined to be $\{0\}$. 
\begin{thm}\label{CAB1} For $J\subseteq I$ define
\begin{eqnarray*}
   P(J):=\bigoplus_{\La\in P^+\cap \overline{F_J}} L(\La)^{(*)}_{\neq\La}\oplus \bigoplus_{\La\in P^+\setminus\overline{F_J}} L(\La)^{(*)}   \subseteq CA .
\end{eqnarray*}
A $G_{fn}$-equivariant embedding of the $G_{fn}$-set $\Omega_{fn}$ into the $G_{fn}$-set $\PCAF$ is given by
\begin{eqnarray*}
\begin{array}{rccc}
  \omega: & \Omega_{fn} &  \to      & \PCAF \\
          & g(P_{fn})_J   &  \mapsto  & g P(J)
\end{array}.
\end{eqnarray*}
Furthermore, for $g(P_{fn})_J,\,h(P_{fn})_K\in\Omega_{fn}$ we have
\begin{eqnarray*}
  g(P_{fn})_J\leq h(P_{fn})_K   \quad\iff\quad   \overline{\{g P(J)\}}^{\,pts}\subseteq \overline{\{h P(K)\}}^{\,pts}.
\end{eqnarray*}
\end{thm}
\Proof Obviously, $P(J)$ is a $P^+$-homogeneous linear space different from $CA$. We first show that $P(J)$ is an ideal. Let $N\in P^+$. If $\La\in P^+\setminus\overline{F_J}$ then also $N+\La \in P^+\setminus\overline{F_J}$, because $P^+\cap\overline{F_J}$ is a face of $P^+$. We get
\begin{eqnarray*}
    L(N)^{(*)}  \prdca P(J)_\La= L(N)^{(*)}  \prdca L(\La)^{(*)}\subseteq L(N+\La)^{(*)}=P(J)_{N+\La}.
\end{eqnarray*}
Now let $\La\in P^+\cap\overline{F_J}$. For every $\phi\in L(N)^{(*)}$ and every $\psi\in L(\La)^{(*)}_{\neq \La}$ we have
\begin{eqnarray*}
   (\phi  \prdca \psi)(v_{N+\La})=\phi(v_N)\underbrace{\psi(v_\La)}_{=0}=0.
\end{eqnarray*}
We find
\begin{eqnarray*}
     L(N)^{(*)}  \prdca P(J)_\La= L(N)^{(*)}  \prdca L(\La)^{(*)}_{\neq \La}\subseteq L(N+\La)^{(*)}_{\neq \La}\subseteq P(J)_{N+\La}.
\end{eqnarray*} 
To show that $P(J)$ is prime it is sufficient to show that the algebra $CA/P(J)$ has no zero divisors. Obviously,  $CA/P(J)$ is isomorphic to the monoid algebra $\FK{\widetilde{P^+\cap\overline{F_J}}}$, which identifies with the monoid algebra $\FK{P^+\cap\overline{F_J}}$. Now $P^+\cap\overline{F_J}$ is a submonoid of the lattice $P$. It follows from Corollary \ref{mr} that  $\FK{P^+\cap\overline{F_J}}$ has no zero divisors.
Because of $\mb{dim\,}(L(\La)^{(*)}/P(J)_\La)\leq 1$ for all $\La\in P^+$ we get $P(J)\in\PCAF$ by Theorem \ref{dim}.

$\omega$ is a well-defined $G_{fn}$-equivariant embedding if and only if for all $J\subseteq I$ it holds 
\begin{eqnarray*}
  \mb{Stab}_{G_{fn}}(P(J)) = (P_{fn})_J.
\end{eqnarray*} 
Let $g\in G_{fn}$. We have $P(J)=gP(J)=\pi(g)^{-1}(P(J))$ if and only if the following equations are satisfied:
\begin{eqnarray}
  \pi(g)^{-1}L(\La)^{(*)}=L(\La)^{(*)} &\mb{ for all }& \La\in P^+\setminus\overline{F_J},\label{eqstab1}\\
  \pi(g)^{-1}L(\La)^{(*)}_{\neq \La}=L(\La)^{(*)}_{\neq \La} &\mb{ for all }& \La\in P^+\cap\overline{F_J} \label{eqstab2}.   
\end{eqnarray}
Equations (\ref{eqstab1}) always hold. Written differently, equations (\ref{eqstab2}) are 
\begin{eqnarray*}
    \Mklz{\phi\in L(\La)^{(*)}}{ \phi(g v_\La)= 0} = \Mklz{\phi\in L(\La)^{(*)}}{ \phi(v_\La)= 0}   &\mb{ for all }& \La\in P^+\cap\overline{F_J} .  
\end{eqnarray*}
Here $\F g v_\La$ and $\F v_\La$ are finite dimensional subspaces of $L(\La)_f$, which is paired nondegenerately with $L(\La)^{(*)}$. Therefore, these equations are equivalent to $\F g v_\La=\F v_\La$ for all $\La\in P^+\cap\overline{F_J}$. (To check this directly use a  base adopted to the decomposition $L(\La)^{(*)}=\bigoplus_{\la\in P(\La)}L(\La)^*_\la$,  $\La\in P^+\cap\overline{F_J}$.)  By Theorem \ref{ngfnllala} (b) and Theorem \ref{spfpp} (a) this is in turn equivalent to $g\in \bigcap_{K\supseteq J}(P_{fn})_K=(P_{fn})_J$. 

Now, $g(P_{fn})_J\leq h(P_{fn})_K$ is by definition equivalent to $g(P_{fn})_J\supseteq h(P_{fn})_K$, which is equivalent $(P_{fn})_J\supseteq g^{-1}h(P_{fn})_K$, which in turn is equivalent to $J\supseteq K$ and $g^{-1}h\in (P_{fn})_J$. On the other hand, we have $\overline{\{gP(J)\}}^{\,pts}\subseteq \overline{\{hP(K)\}}^{\,pts}$ if and only if $gP(J)\in\overline{\{hP(K)\}}^{\,pts}=\overline{\{hP(K)\}}\cap \PCAF$, if and only if $g P(J)\supseteq h P(K)$, if and only if $ P(J)\supseteq g^{-1}h P(K)=\pi(g^{-1}h)^{-1}(P(K))$, which is equivalent to the inclusions
\begin{eqnarray}
\bigoplus_{\La\in P^+\setminus \overline{F_J}} L(\La)^{(*)}  & \supseteq & \bigoplus_{\La\in P^+\setminus \overline{F_K}} L(\La)^{(*)} ,\label{rel1}\\
 \bigoplus_{\La\in P^+\cap \overline{F_J}} L(\La)^{(*)}_{\neq \La} & \supseteq &  \bigoplus_{\La\in P^+\cap \overline{F_K}\cap\overline{F_J}} \pi(g^{-1}h)^{-1} L(\La)^{(*)}_{\neq \La} .\label{rel2}
\end{eqnarray}
Inclusion (\ref{rel1}) is equivalent to $P^+\setminus \overline{F_J}\supseteq P^+\setminus \overline{F_K}$, which is equivalent to $\overline{F_J}\subseteq \overline{F_K}$, which in turn is equivalent to $J\supseteq K$. Therefore,
inclusion (\ref{rel2}) is equivalent to the inclusions
\begin{eqnarray*}
  L(\La)^{(*)}_{\neq \La}\supseteq \pi(g^{-1}h)^{-1} L(\La)^{(*)}_{\neq \La} \quad\mb{ for all }\quad\La\in P^+\cap\overline{F_J}.
\end{eqnarray*}
Equality always holds, because $L(\La)^{(*)}_{\neq \La}$ and $ \pi(g^{-1}h)^{-1} L(\La)^{(*)}_{\neq \La} =\pi(h^{-1}g)L(\La)^{(*)}_{\neq \La}$ are 1-co\-di\-men\-sion\-al subspaces of $L(\La)^{(*)}$. As we have seen by the calculation of the stabilizers, these equations are equivalent to $g^{-1}h\in (P_{fn})_J$.\qed

It remains to show that the map $\omega$ of Theorem \ref{CAB1} is surjective, which is not straightforward. It is reached in Corollary \ref{surjofOmega} as a consequence of our description of the full projective spectrum $\PCA$, which we investigate next. For that we first introduce a $G_{fn}$-invariant stratification of $\PCA$, whose strata extend the orbits $G_{fn} P(J)$, $J\subseteq I$. 

The Cartan algebra $CA$ is a $P^+$-graded algebra without zero divisors. The faces of $P^+$ are $P^+\cap\overline{F_J}$, $J\subseteq I$. It follows that for every $J\subseteq I$ there is the semidirect decomposition 
\begin{eqnarray}\label{semi1}
  CA = \bigoplus_{\La\in P^+\cap \overline{F_J}}L(\La)^{(*)} \oplus \bigoplus_{\La\in P^+\setminus \overline{F_J}} L(\La)^{(*)}
\end{eqnarray}
with graded subalgebra $\bigoplus_{\La\in P^+\cap \overline{F_J}}L(\La)^{(*)}$ and graded prime ideal $\bigoplus_{\La\in P^+\setminus \overline{F_J}} L(\La)^{(*)}$, a sum over the empty set defined to be $\{0\}$.

For every $J\subseteq I$ set
\begin{eqnarray*}
    \overline{Or(J)}:= {\mathcal V}(\bigoplus_{\La\in P^+\setminus \overline{F_J}} L(\La)^{(*)}) = \overline{\{ \bigoplus_{\La\in P^+\setminus \overline{F_J}} L(\La)^{(*)}\}} \, ,
\end{eqnarray*}
which is a $G_{fn}$-invariant closed subset of $\PCA$. In particular, $\overline{Or(\emptyset)}=\PCA$ and $\overline{Or(I)}=\{\bigoplus_{\La\in P^+\setminus \overline{F_I}} L(\La)^{(*)}\}$. The following properties are trivial to check:

\begin{prop}\label{OrbitAbschluss1} Let $K,\, J\subseteq I$.
\begin{itemize}
\item[(a)] $\overline{Or(K)}\subseteq \overline{Or(J)}$ if and only if $K\supseteq J$. 
\item[(b)] $\overline{Or(J)}\cap\overline{Or(K)}=\overline{Or(J\cup K)}$.
\end{itemize}
\end{prop}

As a closed set $\overline{Or(J)}$ can be described as a projective spectrum:
\begin{prop}\label{OrbitAbschluss2} Let $J\subseteq I$. The map
\begin{eqnarray*}
      \mb{Proj\,}(\bigoplus_{\La\in P^+\cap \overline{F_J}}L(\La)^{(*)})    &\to &\quad\quad\overline{Or(J)}\\
        Q \quad\quad\quad &\mapsto &  Q\oplus \bigoplus_{\La\in P^+\setminus\overline{F_J}}L(\La)^{(*)}
\end{eqnarray*}
is a homeomorphism, mapping $\mb{Proj\,}(\bigoplus_{\La\in P^+\cap \overline{F_J}}L(\La)^{(*)})(\F)$ bijectively to $\overline{Or(J)}(\F)$.
\end{prop}
\Proof The restriction of the canonical map $CA\to CA/\bigoplus_{\La\in P^+\setminus \overline{F_J}} L(\La)^{(*)}$ to the graded subalgebra $\bigoplus_{\La\in P^+\cap \overline{F_J}}L(\La)^{(*)}$ is an isomorphism of graded algebras. The proposition now follows from Theorem \ref{csub}.
\qed

For every $J\subseteq I$ set
\begin{eqnarray*}
  Or(J):=\overline{Or(J)}\setminus\bigcup_{I\supseteq K\supsetneqq J}\overline{Or(K)}
      =\overline{Or(J)}\setminus\bigcup_{i\in I\setminus J}\overline{Or(J\cup\{i\})}.
\end{eqnarray*}

\begin{prop}\label{Orbit2} Let $J\subseteq I$.
\begin{itemize}
\item[(a)] $Or(J)$ is $G_{fn}$-invariant and $P(J)\in Or(J)$.
\item[(b)] $Or(J)$ is open and dense in $\overline{Or(J)}$. In particular, $Or(\emptyset)$ is open and dense in $\PCA$.
\item[(c)] It holds
\begin{eqnarray*}
     \overline{Or(J)} = \dot{\bigcup_{I\supseteq K\supseteq J}} Or(K).
\end{eqnarray*}
In particular, $\PCA = \dot{\bigcup}_{K\subseteq I} Or(K)$.
\end{itemize}
\end{prop}
\Proof To (a) and (b): By its definition, $Or(J)$ is open in $\overline{Or(J)}$.  It is a $G_{fn}$-invariant set because the sets $\overline{Or(K)}$, $K\supseteq J$, are $G_{fn}$-invariant. It is trivial to check that $P(J)$ as well as $\bigoplus_{\La\in P^+\setminus \overline{F_J}} L(\La)^{(*)}$ are contained in $Or(J)$. It follows that $Or(J)$ is dense in $\overline{Or(J)}$ .

To (c): By Proposition \ref{OrbitAbschluss1} we have
\begin{eqnarray*}
    \overline{Or(J)}=\bigcup_{I\supseteq K\supseteq J}\overline{Or(K)}\supseteq\bigcup_{I\supseteq K\supseteq J} Or(K).
\end{eqnarray*} 
Let $Q\in \overline{Or(J)}$. Choose a set $I\supseteq K_{max}\supseteq J$, maximal with respect to the inclusion, such that $Q\in \overline{Or(K_{max})}$. Then $Q\in Or(K_{max})$. 

Let $K_1,\, K_2\subseteq I$ and $K_1 \neq K_2$. By Proposition \ref{OrbitAbschluss1} we get
\begin{eqnarray}\label{ordisj1}
  Or(K_1)\cap Or(K_2) \subseteq  \overline{Or(K_1)}\cap \overline{Or(K_2)} \subseteq  \overline{Or(K_1\cup K_2)}.
\end{eqnarray}
At least one of the sets $K_1$, $K_2$ is nonempty. Assume $K_2$ is nonempty. Then $K_1\cup K_2\supsetneqq K_1$. By the definition of $Or(K_1)$ we get
\begin{eqnarray}\label{ordisj2}
   Or(K_1)  \cap \overline{Or(K_1\cup K_2)}  =  \emptyset.
\end{eqnarray}
From (\ref{ordisj1}) and (\ref{ordisj2}) we find
\begin{eqnarray*}
  Or(K_1)\cap Or(K_2)  =  Or(K_1)\cap Or(K_2) \cap \overline{Or(K_1\cup K_2)}  =  \emptyset. 
\end{eqnarray*}
\qed

Let $J\subseteq I$. The Or-stratification at $P(J)$ will be investigated by restricting to a principal open subset at $P(J)$. Define
\begin{eqnarray*}
   D(J):=\Mklz{Q\in \PCA}{\delta_N\notin Q \mb{ for all } N\in P^+\cap \overline{F_J}}.
\end{eqnarray*}
In particular, $D(I)= \PCA$. This follows because $\delta_N$ is a unit of $CA$ for every $N\in P^+\cap\overline{F_I}=\sum_{i=n+1}^{2n-l}\mathbb{Z}\La_i$. 

To describe $D(J)$ as a projective spectrum we identify $CA$ with its image in $(\widetilde{P^+\cap\overline{F_J}})^{-1}CA$, which is possible because $CA$ has no zero divisors. By Theorem \ref{loc} the following holds:
\begin{prop}\label{Sp3} Let $J\subseteq I$. The map 
\begin{eqnarray} \label{locsp}
\begin{array}{ccc}
       \mb{Proj\,}((\widetilde{P^+\cap\overline{F_J}})^{-1}CA )  &\to & D(J)\\
           Q &\mapsto & Q\cap CA
\end{array}
\end{eqnarray}
is a homeomorphism, mapping $\mb{Proj\,}( (\widetilde{P^+\cap\overline{F_J}})^{-1}CA ) (\F)$ bijectively to $D(J)(\F)$.
\end{prop}

\begin{rem}\label{PrincipalOpen1} It is easy to check that $ri(P^+\cap\overline{F_J})=P^+\cap F_J$ is the set of principal elements of the monoid $P^+\cap\overline{F_J}$. By Theorem \ref{poset} the set $D(J)$ is principal open, i.e., for every $\La\in P^+\cap F_J$ we have
\begin{eqnarray*}
\begin{array}{ccc}
     D(J)&=& D(\delta_\La),\\
       \widetilde{(P^+\cap \overline{F_J})}^{-1}CA &\cong & \Mklz{\delta_\La^n}{n\in\Nn}^{-1}CA.
\end{array}
\end{eqnarray*}
We will use the description which is  independent of the choice of $\La\in P^+\cap F_J$ in the formulation of theorems and propositions. In proofs we will make use of both descriptions.
\end{rem}

$D(J)$ decomposes by the Or-stratification as follows:
\begin{prop}\label{PrincipalOpenStrata} Let $J\subseteq I$. Then
\begin{eqnarray*}
  D(J) \;=\; \dot{\bigcup_{K\subseteq J}} \Big( D(J)\cap Or(K) \Big) \qquad \mb{ and }\qquad P(K)\in D(J)\cap Or(K)\; \mb{ for all }\; K\subseteq J\,.
\end{eqnarray*}
\end{prop}
\Proof It is $P(K)\in Or(K)$ for all $K\subseteq I$ by Proposition \ref{Orbit2} (a). It is trivial to check that also $P(K)\in D(J)$ for all $K\subseteq J$. 

Let $K\subseteq I$. If $\emptyset \neq D(J)\cap Or(K) \subseteq D(J) \cap \overline{Or(K)}$ there exists $Q\in\PCA$ such that
\begin{eqnarray*}
  Q\supseteq \bigoplus_{\La\in P^+\setminus \overline{F_K}}L(\La)^{(*)} \quad \mb{ and }\quad \delta_\La\not\in Q \;\mb{ for all }\; \La\in P^+\cap\overline{F_J}\,.
\end{eqnarray*}
It follows that $\overline{F_J}\subseteq \overline{F_K}$, which is equivalent to $J\supseteq K$. Now use Proposition \ref{Orbit2} (c).
\qed

Next we describe the open set obtained as the $G_{fn}$-orbit of $D(J)$. For the theorem we use the definition ${\mathcal V}(S):=\{Q\in\PCA | Q\supseteq S\}$ for subsets $S\subseteq CA$. Note that ${\mathcal V}(S)={\mathcal V}((S_{hom}))$, where $(S_{hom})$ is the ideal generated by the set $S_{hom}$ of homogeneous components of the elements of $S$.
\begin{thm}\label{PrincipalOpen2} Let $J\subseteq I$.
\begin{itemize} 
\item[(a)] Then
\begin{eqnarray*}
    \bigcup_{g\in G_{fn}} g D(J)\; = \; \PCA\setminus{\mathcal V}(\bigoplus_{\La\in P^+\cap F_J}L(\La)^{(*)})  .
\end{eqnarray*}
\item[(b)] We have $(P_{fn}^-)_J\, D(J)=D(J)$. In particular, the union in (a) can be taken over sets of coset representatives
of $G_{fn}/(P_{fn}^-)_J$.
\end{itemize}
\end{thm}
\Proof We first show (a). Let $Q\in \PCA$. Then 
\begin{eqnarray*}
    Q\notin \bigcup_{g\in G_{fn}} g D(J)
\end{eqnarray*} 
is equivalent to $\pi(g)Q\notin D(J)$ for all $g\in G_{fn}$. By Remark \ref{PrincipalOpen1} this is equivalent to
\begin{eqnarray}\label{uebhovoll}
    \delta_\La\in \pi(g)Q \quad\mb{ for all }\quad  g\in G_{fn} \;\mb{ and } \;\La\in P^+\cap F_J.
\end{eqnarray}
The  irreducible $(G_{fn})^{op}$-module $L(\La)^{(*)}$ is spanned by $\pi(G_{fn})\delta_\La$. Therefore, since $Q$ is homogeneous, (\ref{uebhovoll}) is equivalent to 
\begin{eqnarray*}
   Q\supseteq \bigoplus_{\La\in P^+\cap F_J}L(\La)^{(*)}.
\end{eqnarray*}

Let $\La\in P^+\cap F_J$. By Remark \ref{PrincipalOpen1} we have $D(J)=\Mklz{Q\in \PCA}{\delta_\La\notin Q}$. Now (b) follows from $\pi((P_{fn}^-)_J)\delta_\La\in\F\delta_\La\setminus\{0\}$.
\qed

Let $J\subseteq I$. From the description of $\overline{Or(J)}$ in Proposition \ref{Orbit2} (c) and from Proposition \ref{PrincipalOpenStrata} we get
\begin{eqnarray*}
 BC(J) := \overline{Or(J)}\cap D(J) =  Or(J)\cap D(J)\subseteq Or(J).
\end{eqnarray*}
Furthermore, $P(J)\in BC(J)$. We call $BC(J)$ the {\it standard big cell} of $Or(J)$. We call every $gBC(J)$, $g\in G_{fn}$, a {\it big cell} of $Or(J)$.

\begin{thm} \label{BigCell1}For $J\subseteq I$ the following hold:
\begin{itemize}
\item[(a)] Every big cell $gBC(J)$, $g\in G_{fn}$, is principal open in $\overline{Or(J)}$ and dense in $\overline{Or(J)}$.
\item[(b)] $Or(J)$ can be covered by big cells:
\begin{eqnarray*}
   Or(J)\:=\bigcup_{g\in G_{fn}} g BC(J)\,.
\end{eqnarray*}
We have $(P_{fn})_J^-\, BC(J)=BC(J)$. In particular, the union can be taken over sets of coset representatives
of $G_{fn}/(P_{fn})_J^-$.
\end{itemize}
\end{thm}
\Proof It is sufficient to show part (a) for $BC(J)$ because $\overline{Or(J)}$ is $G_{fn}$-invariant. $BC(J)$ is principal open in $\overline{Or(J)}$ because $D(J)$ is principal open in $\PCA$ by Remark \ref{PrincipalOpen1}. By the definition of $\overline{Or(J)}$ and $BC(J)$ we get
\begin{eqnarray*}
   \bigoplus_{\La\in P^+\setminus \overline{F_J}}L(\La)^{(*)} \in BC(J)=\overline{Or(J)}\cap D(J)\subseteq \overline{Or(J)}=
\overline{\{\bigoplus_{\La\in P^+\setminus \overline{F_J}}L(\La)^{(*)}\}}.
\end{eqnarray*}
Therefore, $BC(J)$ is dense in $\overline{Or(J)}$.

For part (b) it is sufficient to show
\begin{eqnarray}\label{Orbit3}
     Or(J)=\overline{Or(J)}\cap \PCA\setminus{\mathcal V}(\bigoplus_{\La\in P^+\cap F_J}L(\La)^{(*)}) .
\end{eqnarray}
Then (b) follows from Theorem \ref{PrincipalOpen2} (a), (b) and the $G_{fn}$-invariance of $\overline{Or(J)}$.
Taking the complement in $\overline{Or(J)}$, the equation (\ref{Orbit3}) is equivalent to
\begin{eqnarray*}
   \bigcup_{K\supsetneqq J}\overline{Or(K)} =  \overline{Or(J)} \cap {\mathcal V}(\bigoplus_{\La\in P^+\cap F_J}L(\La)^{(*)}) .
\end{eqnarray*}
Inserting the definition of $\overline{Or(K)}$, $K\supseteq J$, this is equivalent to
\begin{eqnarray}\label{Voffuebfgl}
\bigcup_{K\supsetneqq J} {\mathcal V}(\bigoplus_{\La\in P^+\setminus \overline{F_K}}L(\La)^{(*)}) = {\mathcal V}(\bigoplus_{\La\in (P^+\cap F_J)\cup P^+\setminus \overline{F_J}}L(\La)^{(*)}) .
\end{eqnarray}
The inclusion "$\subseteq$" is valid, because for all $K\supsetneqq J$ we have
\begin{eqnarray*}
 P^+\setminus\overline{F_K}= P^+\setminus\overline{F_J}\cup (P^+\cap (\overline{F_J}\setminus\overline{F_K})) \supseteq P^+\setminus\overline{F_J} \cup (P^+\cap F_J).
\end{eqnarray*}
Now we show "$\supseteq$". Let $Q$ be in the right hand side of (\ref{Voffuebfgl}). Suppose that $Q$ is not contained in the left hand side of (\ref{Voffuebfgl}). Then for every $K\supsetneqq J$ there exists a weight $\La_K\in P^+\setminus \overline{F_K}$ and an element $\phi_K\in L(\La_K)^{(*)}$ such that $\phi_K\notin Q$. Since $Q$ is prime we obtain
\begin{eqnarray}\label{orclcon}
 Q\not\ni \prod_{K\supsetneqq J} \phi_K\in L(\sum_{K\supsetneqq J}\La_K)^{(*)}.
\end{eqnarray}
Since $P^+\setminus\overline{F_K}$ is a semigroup ideal of $P^+$ we find
\begin{eqnarray*}
   \sum_{K\supsetneqq J}\La_K &\in & \bigcap_{K\supsetneqq J} P^+\setminus\overline{F_K}
  \, = \, \bigcap_{K\supsetneqq J}\left( P^+\setminus\overline{F_J}\cup (P^+\cap (\overline{F_J}\setminus\overline{F_K})) \right)\\
   &&=\, P^+\setminus\overline{F_J} \cup \bigcap_{K\supsetneqq J}\left(  P^+\cap (\overline{F_J}\setminus\overline{F_K}) \right)
    \,= \,P^+\setminus\overline{F_J} \cup (P^+\cap F_J).
\end{eqnarray*}
Therefore, (\ref{orclcon}) contradicts that $Q$ is contained in the right hand side of (\ref{Voffuebfgl}).
\qed

As a principal open set of a closed set the standard big cell $BC(J)$ can be described as a projective spectrum:

\begin{prop}\label{BigCell2a} Let $J\subseteq I$. Then
\begin{eqnarray*}
      \mb{Proj\,}((\widetilde{P^+\cap\overline{F_J}})^{-1}\bigoplus_{\La\in P^+\cap \overline{F_J}}L(\La)^{(*)})    &\to &\quad\quad\quad\quad\quad\quad BC(J)\\
        Q \quad\quad\quad \quad\quad\quad &\mapsto &  (Q\cap \bigoplus_{\La\in P^+\cap \overline{F_J}}L(\La)^{(*)})\oplus \bigoplus_{\La\in P^+\setminus\overline{F_J}}L(\La)^{(*)}
\end{eqnarray*}
is a homeomorphism, mapping $\mb{Proj\,}((\widetilde{P^+\cap\overline{F_J}})^{-1}\bigoplus_{\La\in P^+\cap \overline{F_J}}L(\La)^{(*)})(\F)$ bijectively to $BC(J)(\F)$.
\end{prop}
\Proof The projective spectrum $\mb{Proj\,}(\bigoplus_{\La\in P^+\cap \overline{F_J}}L(\La)^{(*)}) )$ maps homeomorphically onto $\overline{Or(J)}$ by the map given in Proposition \ref{OrbitAbschluss2}. The preimage of $BC(J)=\overline{Or(J)}\cap D(J)$ is the set
\begin{eqnarray*}
    \Big\{\, Q\in \mb{ Proj\,} (\bigoplus_{\La\in P^+\cap \overline{F_J}}L(\La)^{(*)} ) \, \Big| \,\delta_\La \notin Q  \mb{ for all } \La\in P^+\cap\overline{F_J} \,\Big\}\,.
\end{eqnarray*}
Now the proposition follows from Theorem \ref{loc}.
\qed

Next we introduce for $J\subseteq I$ a closed subset of $D(J)$, which is a transversal slice to $Or(J)$ at the point $P(J)$ as we will see later. We need some easy preparations, before we can give its definition.  

For $\emptyset\neq J\subseteq I$ the Lie subalgebra $\g_J$ of $\g$ identifies with the derived Kac-Moody algebra $\g(A_J)'$ to the generalized Cartan submatrix $A_J$ of $A$. The sublattice $Q_J$ of $Q$ identifies with the root lattice, the sublattice $P_J$ of $P$ identifies with the weight lattice of $\g(A_J)'$. The following theorem is well-known. Since the author has not found a reference, the nontrivial parts of the proof are added.
\begin{thm}\label{sm} Let $\emptyset\neq J\subseteq I$ and $\La\in P^+$. Decompose $\La=\La_J+\La^J$ with $\La_J\in P_J^+$ and $\La^J\in P^+\cap\overline{F_J}$. Regard $L(\La)$ as a $\g_J$-module. 
\begin{itemize}
\item[(a)] Define an equivalence relation on $P(\La)$ by 
\begin{eqnarray*}
\la\sim\mu \quad :\iff \quad \la-\mu\in Q_J.
\end{eqnarray*}
If $C$ is an equivalence class, then $L(\La)_C:=\bigoplus_{\la\in C}L(\La)_\la$ is a $\g_J$-module of $L(\La)$. It decomposes into a direct sum of irreducible highest weight modules of $\g_J$ with highest weights in $P_J^+$.
\item[(b)] The equivalence class $[\La]$ of $\La$ is given by $[\La]= P(\La)\cap \left(\La-(Q_J)_0^+\right)$. The $\g_J$-module $L(\La)_{[\La]}$ is an irreducible highest weight module of $\g_J$ of highest weight $\La_J$ and highest weight space $L(\La)_\La$.
\end{itemize}
\end{thm}
\Proof The decomposition of $L(\La)_C$  into a direct sum of irreducible highest weight modules of $\g_J$ with highest weights in $P_J^+$ stated in (a) follows from Theorem 10.7 a) of \cite{K}. For (b) we only show $U(\n_J^-)L(\La)_\La=\bigoplus_{\la\in [\La]} L(\La)_\la$, where  $[\La]= P(\La)\cap \left(\La-(Q_J)_0^+\right)$.  We have
\begin{eqnarray*}
    L(\La)=U(\n^-)L(\La)_\La= U((\n^-)^J)U(\n^-_J)L(\La)_\La
    = U(\n^-_J)L(\La)_\La + \sum_{\beta\in Q^+\setminus Q_J^+} U((\n^-)^J)_{-\beta}\,U(\n^-_J)L(\La)_\La.
\end{eqnarray*} 
Here
\begin{eqnarray} 
  U(\n^-_J)L(\La)_\La&\subseteq& \bigoplus_{\la\in P(\La)\cap (\La-(Q_J)^+_0)}L(\La)_\la , \label{LJmod1}\\
  \sum_{\beta\in Q^+\setminus Q_J^+} U((\n^-)^J)_{-\beta}\, U(\n^-_J)L(\La)_\La &\subseteq& \bigoplus_{\la\in P(\La)\cap (\La-(Q_J)^+_0-Q^+\setminus Q_J^+)}L(\La)_\la .\label{LJmod2}
\end{eqnarray}
Now $L(\La)$ is the direct sum of the sums in (\ref{LJmod1}) and (\ref{LJmod2}) on the right. It follows that in (\ref{LJmod1}) and (\ref{LJmod2}) equality holds.
\qed

For $J\subseteq I$ and $\La\in P^+$ set 
\begin{eqnarray*}
   \quad  P_J(\La):=P(\La)\cap \left(\La-(Q_J)_0^+\right)\,.
\end{eqnarray*}
\begin{rem}\label{PJp} In particular, we have:
\begin{itemize}
\item[(a)] $P_I(\La)=P(\La)$ for all $\La\in P^+$.
\item[(b)] $P_J(\La)=\{\La\}$ for all $\La\in P^+\cap\overline{F_J}$ and $J\subseteq I$.
\end{itemize}
Here (a) and the case $J=\emptyset$ in (b) follow from the definition. (b) also holds for $J\neq\emptyset$, because the $\g_J$-module $L(\La)_{[\La]}$ is trivial.
\end{rem}

For $J\subseteq I$ and $\La\in P^+$ set
\begin{eqnarray*}
  && L_J(\La):=\bigoplus_{\la\in P_J(\La)} L(\La)_\la  \;\subseteq \;L(\La) \quad\mb{ and }\quad  L_J(\La)_f:=\prod_{\la\in P_J(\La)} L(\La)_\la  \;\subseteq \;L(\La)_f\,, \\
  && L_J(\La)^{(*)} := \bigoplus_{\la\in P_J(\La)}L(\La)_\la^*\;\subseteq\; L(\La)^{(*)}\,.    
\end{eqnarray*}
Here $L_J(\La)^{(*)}$ is the restricted dual of the $\g_J$-module $L_J(\La)$. Note that we have
\begin{eqnarray}\label{rwts}
   P(\La)\setminus P_J(\La) \,=\, P(\La)\cap \Big(\La- \big(Q_0^+\setminus (Q_J)_0^+\big)\Big)\,,
\end{eqnarray}
where we adopted the convention $B\pm\emptyset := \pm\emptyset + B := \emptyset$ for $B\subseteq P$. Now set
\begin{eqnarray*}
   && R_J(\La):= \bigoplus_{\la\in P(\La)\setminus P_J(\La)} L(\La)_\la  \;\subseteq \;L(\La) \quad  \mb{ and }  \quad R_J(\La)_f:= \prod_{\la\in P(\La)\setminus P_J(\La)} L(\La)_\la  \;\subseteq \;L(\La)_f\,,\\
  && R_J(\La)^{(*)} := \bigoplus_{\la\in P(\La)\setminus P_J(\La)}L(\La)_\la^*\;\subseteq\; L(\La)^{(*)}\,,
\end{eqnarray*} 
a sum or product over the empty set defined to be $\{0\}$.

$L_J(\La)^{(*)}$, $R_J(\La)^{(*)}$ are $(G_{fn})_J^{op}$-invariant subspaces of $L(\La)^{(*)}$, and $L_J(\La)_f$, $R_J(\La)_f$ are $(G_{fn})_J$-invariant subspaces of $L(\La)_f$. Furthermore, $L_J(\La)$, $R_J(\La)$ are $G_J$-invariant subspaces of $L(\La)$. 
\begin{prop}\label{smdu} Let $J\subseteq I$. The following inclusions hold in $CA$ for all $\La,\,N\in P^+$: 
\begin{itemize}
\item[(a)] $ L_J(\La)^{(*)}\prdca L_J(N)^{(*)} \,\subseteq\, L_J(\La+N)^{(*)}$.
\item[(b)] $R_J(\La)^{(*)}\prdca L(N)^{(*)}\,\subseteq \,R_J(\La+N)^{(*)}$ and $ L(\La)^{(*)}\prdca R_J(N)^{(*)} \,\subseteq \,R_J(\La+N)^{(*)}$. 
\end{itemize}
\end{prop}
\Proof We first show 
\begin{eqnarray}\label{wCP}
  L(\La)_\la^*\prdca L(N)_\mu^*\subseteq L(\La+N)_{\la+\mu}^* \quad \mb{ for all }\quad \la,\,\mu\in P.
\end{eqnarray}
Let $\phi_\la\in L(\La)^*_\la$, $\psi_\mu\in L(N)_\mu^*$, $\la,\,\mu\in P$. With the notation of Subsection \ref{subsectionCA} the Cartan product $\phi_\la\prdca \psi_\mu\in L(\La+N)^{(*)}$ is given by
\begin{eqnarray*}
    (\phi_\la \prdca \psi_\mu)(w) :=(\phi_\la\otimes\psi_\mu)(\Phi(w)) \quad\mb{ where }\quad w\in L(\La+N).
\end{eqnarray*}
The $G$-equivariant linear map $\Phi$ preserves the weights. We get immediately $\phi_\la \prdca \psi_\mu\in L(\La+N)_{\la+\mu}^*$. 

As an example we show $R_J(\La)^{(*)}\prdca L(N)^{(*)}\subseteq R_J(\La+N)^{(*)} $. The other inclusions can be treated similarly. 
This inclusion holds trivially for $J=I$, because $R_I(\La)^{(*)}=\{0\}$ and $R_I(\La +N)^{(*)}=\{0\}$.
For $J\neq I$ use (\ref{rwts}) to rewrite the definitions of $R_J(\La)^{(*)}$ and $R_J(\La+N)^{(*)}$. Then the inclusion follows from
\begin{eqnarray*}
   \left(\La-\big(Q^+_0\setminus (Q_J)^+_0\big)\right) + \left( N -Q ^+_0\right) \;\subseteq \;(\La+ N)-\big(Q^+_0\setminus (Q_J)^+_0\big)\,.
\end{eqnarray*} 
with (\ref{wCP}). 
\qed

Let $J\subseteq I$. The Cartan algebra $CA$ has no zero divisors. With the previous proposition we get the semidirect decomposition
\begin{eqnarray}\label{semi2}
  CA= \bigoplus_{\La\in P^+} L_J(\La)^{(*)} \oplus \bigoplus_{\La\in P^+} R_J(\La)^{(*)} 
\end{eqnarray}
with graded subalgebra $\bigoplus_{\La\in P^+} L_J(\La)^{(*)}$ and graded prime ideal $\bigoplus_{\La\in P^+} R_J(\La)^{(*)}$. We call
\begin{eqnarray*}
     S(J) \;:=\; {\mathcal V}( \bigoplus_{\La\in P^+} R_J(\La)^{(*)} ) \cap D(J)
\end{eqnarray*}
the standard transversal slice to $Or(J)$ at the point $P(J)\in Or(J)$. This name will be justified later. We now only show:

\begin{thm}\label{Slice1} Let $J\subseteq I$. Then
\begin{eqnarray*}
     Or(J) \cap S(J) \; = \; BC(J) \cap S(J) \; = \; \{ P(J) \}.
\end{eqnarray*}
\end{thm}
\Proof Because of $S(J)\subseteq D(J)$ we find
\begin{eqnarray*}
  BC(J)\cap S(J) \subseteq Or(J)\cap S(J)\subseteq \overline{Or(J)}\cap S(J) = \overline{Or(J)}\cap \cap D(J)\cap S(J) = BC(J) \cap S(J).
\end{eqnarray*}
By the definition of $BC(J)$, $S(J)$, and by Remark \ref{PJp} (b) we obtain
\begin{eqnarray*}
  BC(J) \cap S(J) &=& {\mathcal V}(\bigoplus_{\La\in P^+}R_J(\La)^{(*)}) \cap {\mathcal V}(\bigoplus_{\La\in P^+\setminus\overline{F_J}}L(\La)^{(*)}) \cap D(J) \\
  &=& {\mathcal V}(\bigoplus_{\La\in P^+}R_J(\La)^{(*)}\oplus \bigoplus_{\La\in P^+\setminus\overline{F_J}}L(\La)^{(*)}) \cap D(J)\\
  &=& {\mathcal V}(\bigoplus_{\La\in P^+\cap\overline{F_J}}L(\La)_{\neq\La}^{(*)}\oplus \bigoplus_{\La\in P^+\setminus\overline{F_J}}L(\La)^{(*)}) \cap D(J) \;=\; {\mathcal V}(P(J))\cap D(J).
\end{eqnarray*}
We get $P(J)\in {\mathcal V}(P(J))\cap D(J)$ from Proposition \ref{PrincipalOpenStrata}. Now let $Q\in {\mathcal V}(P(J))\cap D(J)$. By definition we have $Q\supseteq P(J)$ and $\delta_\La\notin Q$ for all $\La\in P^+\cap\overline{F_J}$. We show $Q\subseteq P(J)$. Every $\phi\in Q$ is of the form
\begin{eqnarray*}
\phi=\sum_{\La\in P^+\cap\overline{F_J}} c_\La \delta_\La + \psi
\end{eqnarray*}
with $\psi\in P(J)$ and $c_\La\in\F$, $c_\La\neq 0$ for at most finitely many $\La\in P^+\cap\overline{F_J}$. It follows that $\sum_{\La\in P^+\cap\overline{F_J}} c_\La \delta_\La \in Q$. Since $Q$ is graded we get $c_\La\delta_\La\in Q$, which implies $c_\La=0$, $\La\in P^+\cap\overline{F_J}$. We found $\phi=\psi\in P(J)$.
\qed

$S(J)$ decomposes with respect to the Or-stratification as follows:
\begin{prop}\label{stratsl} Let $J\subseteq I$. Then
\begin{eqnarray*}
  S(J) \;=\; \dot{\bigcup_{K\subseteq J}} \Big( S(J)\cap Or(K) \Big) \qquad \mb{ and }\qquad P(K)\in S(J)\cap Or(K)\; \mb{ for all }\; K\subseteq J\,.
\end{eqnarray*}
\end{prop}

\Proof It is trivial to check that $P(K)\supseteq \bigoplus_{\La\in P^+}R_J(\La)^{(*)}$ for all $K\subseteq I$. Now the proposition follows from the definition of $S(J)$ and Proposition \ref{PrincipalOpenStrata}.
\qed

As a principal open set of a closed set the standard transversal slice $S(J)$ can be described as a projective spectrum. The proof is completely similar as for $BC(J)$, the semidirect decomposition (\ref{semi1}) replaced by the semidirect decomposition (\ref{semi2}):
\begin{prop}\label{Slice2} Let $J\subseteq I$. The map
\begin{eqnarray*}
      \mb{Proj\,}((\widetilde{P^+\cap\overline{F_J}})^{-1}\bigoplus_{\La\in P^+}L_J(\La)^{(*)})    &\to &\quad\quad\quad\quad\quad\quad S(J)\\
        Q \quad\quad\quad \quad\quad\quad &\mapsto &  (Q\cap \bigoplus_{\La\in P^+ }L_J(\La)^{(*)})\oplus \bigoplus_{\La\in P^+}R_J(\La)^{(*)}
\end{eqnarray*}
is a homeomorphism, mapping $\mb{Proj\,}((\widetilde{P^+\cap\overline{F_J}})^{-1}\bigoplus_{\La\in P^+}L_J(\La)^{(*)})(\F)$ bijectively to $S(J)(\F)$.
\end{prop}

The projective spectrum of the graded algebra $(\widetilde{P^+\cap\overline{F_J}})^{-1}CA$ gives the principal open set $D(J)$ by Proposition \ref{Sp3}. Projective spectra of graded subalgebras of  $(\widetilde{P^+\cap\overline{F_J}})^{-1}CA$ give the standard big cell $BC(J)$ and the standard transversal slice $S(J)$ by Propositions  \ref{BigCell2a} and \ref{Slice2}.
The next aim is to describe the graded algebra $(\widetilde{P^+\cap\overline{F_J}})^{-1}CA$ explicitly as a tensor product of three graded algebras. 
Some preparations are required, before it can be reached in Theorem \ref{Sp4}.  \\

Let $J\subseteq I$. We equip $(\n^-_f)^J$ with the unital coordinate ring of functions $\FK{(\n^-_f)^J}$, which is generated by the linear functions contained in the subspace
\begin{eqnarray*}
   ((\n^-)^J)^{(*)}:=\bigoplus_{\al\in (\Delta^-)^J}\g_\al^*\;\subseteq \;((\n_f^-)^J)^*.
\end{eqnarray*}
It is easy to check that this coordinate ring is the symmetric algebra
\begin{eqnarray*}
        \FK{ (\n^-_f)^J}= \mb{Sym} \left(((\n^-)^J)^{(*)}\right)\,.
\end{eqnarray*}
We push this coordinate ring forward by the bijective exponential map $exp:\,(\n^-_f)^J\to (U^-_f)^J$ and obtain a coordinate ring 
$\FK{(U^-_f)^J}$ on $(U^-_f)^J$ such that the comorphism
\begin{eqnarray*}
   exp^*:\FK{ (U^-_f)^J}\to \FK{ (\n^-_f)^J}
\end{eqnarray*}
exists and is an isomorphism of algebras. We equip $\FK{(U_f^-)^J}$ with the trivial gradation $\{0\}$.

Trivially,  $(\n_f^-)^I=\{0\}$, $\FK{(\n_f^-) ^I}\cong \F$ and $(U_f^-)^I=\{1\}$, $\FK{(U_f^-)^I}\cong\F$.
\begin{rem} $\FK{(\n^-_f)^J}$ is the coordinate ring of the pronilpotent Lie algebra $(\n^-_f)^J$ as follows:  For $k\in\Nn$ it is
\begin{eqnarray*}
  (\n^-_f)^J(k) :=\prod_{\al\in(\Delta^-)^J,\,-ht(\al)\geq k+1}\g_\al
\end{eqnarray*} 
an ideal of $(\n^-_f)^J$, such that the quotient $(\n^-_f)^J/(\n^-_f)^J(k)$ is a finite dimensional nilpotent Lie algebra. There is the system of projections
\begin{eqnarray*}
  p_{kl}:\, (\n^-_f)^J/(\n^-_f)^J(k) \to (\n^-_f)^J/(\n^-_f)^J(l) \quad \mb{ where }\quad k,\,l\in\Nn,\,k\geq l\,.   
\end{eqnarray*}
It is easy to check that a projective limit is given by the projections
\begin{eqnarray*}
   p_k:\, (\n^-_f)^J \to (\n^-_f)^J/(\n^-_f)^J(k) \quad\mb{ where }\quad k\in\Nn\,. 
\end{eqnarray*}

Equip $(\n^-_f)^J/(\n^-_f)^J(k)$ with its coordinate ring $\F\,[(\n^-_f)^J/(\n^-_f)^J(k)]$ as a finite dimensional $\F$-linear space. There is the system of comorphisms
\begin{eqnarray*}
    p_{kl}^*:\, \F\,[(\n^-_f)^J/(\n^-_f)^J(l)] \to \F\,[(\n^-_f)^J/(\n^-_f)^J(k)] \quad \mb{ where }\quad l,\,k\in\Nn,\,l\leq k\,.   
\end{eqnarray*}
It is easy to check that a inductive limit is given by the comorphisms
\begin{eqnarray*}
  p_k^*:\, \F\,[(\n^-_f)^J/(\n^-_f)^J(k)] \to \F\,[(\n^-_f)^J ] \quad\mb{ where }\quad k\in\Nn\,. 
\end{eqnarray*}

Similarly, $\FK{(U^-_f)^J}$ is the coordinate ring of the prounipotent group $(U^-_f)^J$.  
\end{rem}

\begin{prop}\label{UniSpecF} Let $J\subseteq I$. Then
\begin{eqnarray*}
    \mb{Spec\,} \left(\FK{(U_f^-)^J}\right)(\mathbb{F}) =\Mklz{I(u)}{u\in (U_f^-)^J },
\end{eqnarray*}
where $I(u)$ denotes the vanishing ideal in $u\in (U_f^-)^J$.
\end{prop}
\Proof By the linear isomorphisms
\begin{eqnarray*}
  (\,(\, (\n^-)^J\,)^{(*)}\,)^* = (\,\bigoplus_{\al\in (\Delta^-)^J}\g_\al^* \;)^* \cong \prod_{\al\in (\Delta^-)^J}(\g_\al^*)^* \cong \prod_{\al\in (\Delta^-)^J}\g_\al=(\n_f^-)^J
\end{eqnarray*}
the elements of the dual of $((\n^-)^J)^{(*)}$ are given by the evaluations in the elements of $(\n_f^-)^J$. By the universal property of the symmetric algebra  we get $\mb{Spec\,}( \FK{(\n_f^-)^J} )(\mathbb{F}) = \{I(x)| x\in (\n_f^-)^J \}$, where $I(x)$ is the vanishing ideal in $x\in (\n_f^-)^J$. Since the comorphism $ exp^*:\FK{ (U^-_f)^J}\to \FK{ (\n^-_f)^J}$ is an isomorphism of coordinate rings the proposition follows.
\qed

For a $\g$-module $V$ contained in the category ${\mathcal O}^p_{int}$, for $v\in V_f$ and $\phi\in V^{(*)}$ the map
\begin{eqnarray*}
\phi(\ins v):\, U_f^-  & \to     & \qquad\quad\F         \\
                  u \;\;\, & \mapsto & (\pi(u)\phi)(v)=\phi(uv)
\end{eqnarray*}
is called the matrix coefficient of $\phi$ and $v$ on $U_f^-$. To keep our notation easy we often denote its restriction to nonempty subsets of $U_f^-$ by the same symbol.
\begin{thm}\label{UniRing} Let $J\subseteq I$. For every $\g$-module $V$ from the category ${\mathcal O}_{int}^p$, for every $v\in V_f$ and $\phi\in V^{(*)}$ it is
\begin{eqnarray}\label{MkOint}
   \phi(\ins v) \in  \FK{(U_f^-)^J} .
\end{eqnarray}
Furthermore, 
\begin{eqnarray}
  && \FK{(U_f^-)^J} = \mb{span}\left\{\phi(\ins v_N)\,\left|\,\phi\in L(N)^{(*)},\;N\in P^+\cap\overline{F_J}\right.\right\}.
\end{eqnarray}
For every $\La\in P^+$ with $\La(h_i)>0$ for all $i\in I\setminus J$ the linear map
\begin{eqnarray*}
  \n^J & \to     & \FK{(U_f^-)^J}\\
   x   & \mapsto & \left(\pi(x)\delta_\La\right)(\ins v_\La)
\end{eqnarray*}
is injective and 
\begin{eqnarray}
  && \FK{(U_f^-)^J} = \mb{Sym}( \;\Mklz{(\pi(x)\delta_\La)(\ins v_\La)\,}{\,x\in \n^J}\;) .
\end{eqnarray}
\end{thm}
\Proof The theorem holds trivially for $J=I$. Let $J\neq I$. For every $\al\in(\Delta^-)^J$ choose a base $y_{\al i}\in\g_\al$, $i=1,\,\ldots,\,m_\al$. Let $\ti{\psi}_{\al i}\in\g_\al^*$, $i=1,\,\ldots,\,m_\al$, be the dual base. Interpret these elements as elements of $((\n^-)^J)^{(*)}\subseteq ((\n_f^-)^J)^*$ and set
\begin{eqnarray*} 
\psi_{\al i}:=(\exp^*)^{-1}(\ti{\psi}_{\al i}) \in \FK{(U_f^-)^J}, \quad i=1,\,\ldots,\,m_\al .  
\end{eqnarray*}
Note that 
\begin{eqnarray*}
\ti{\psi}_{\al i}(y_{\beta j})=\delta_{\al\beta}\delta_{ij} \quad \mb{ for }\quad\al,\,\beta\in (\Delta^-)^J,\;  i=1,\,\ldots,\,m_\al,\;j=1,\,\ldots,\,m_\beta  .
\end{eqnarray*}
For every $\al\in (\Delta^+)^J$ the pairing $\iBl:\g_\al\times\g_{-\al}\to\F$ is nondegenerate. Let $z_{\al i}\in\g_\al$, $\al\in (\Delta^+)^J$, $i=1,\,\ldots,\,m_\al$, be the dual base of $y_{-\al i}\in\g_{-\al}$, $i=1,\,\ldots,\,m_{-\al}$.
Because of $\iB{\g_\al}{\g_{-\beta}}=\{0\}$ for $\al,\,\beta\in (\Delta^+)^J$, $\al\neq\beta$, we find
\begin{eqnarray*}
 \iB{ z_{\al i} }{ y_{-\beta j} }=\delta_{\al\beta}\delta_{ij} \quad \mb{ for }\quad\al,\,\beta\in(\Delta^+)^J,\;  i=1,\,\ldots,\,m_\al,\;j=1,\,\ldots,\,m_{-\beta}  .
\end{eqnarray*}

(a) Let $v=\sum_{\la\in P(V)}v_\la\in V_f$ with $v_\la\in V_\la$ for all $\la\in P(V)$. Let $\phi_\mu\in V_\mu^*$, $\mu\in P(V)$. We have
\begin{eqnarray*}
    \phi_\mu(exp(x)v) = \sum_{\la\in P(V)}\phi_\mu(exp(x)v_\la) = \sum_{\la\in P(V),\,\la\geq\mu}\phi_\mu(exp(x)v_\la) =  \Big(\sum_{\la\in P(V),\,\la\geq\mu}\phi_\mu(\ins v_\la)\Big)(exp(x))
\end{eqnarray*}
for all $x\in (\n_f^-)^J$, a sum over the empty set defined to be zero. The sum on the right is finite, because $V$ is contained in
${\mathcal O}^p_{int}$.

Therefore, it is sufficient to show (\ref{MkOint}) for $v_\la\in V_\la$, $\phi_\mu\in V_\mu^*$, $\la,\,\mu\in P(V)$, $\la\geq\mu$. For every $x=\sum_{\al\in(\Delta^-)^J}x_\al\in (\n_f^-)^J$ we have
\begin{eqnarray}
   \phi_\mu(exp(x)v_\la)=\sum_{k=0}^{ht(\la-\mu)}\frac{1}{k!}\phi_\mu(x^k v_\la) 
   =\phi_\mu(v_\la)+ \sum_{k=1}^{ht(\la-\mu)}\frac{1}{k!}\sum_{\beta_1,\,\ldots,\,\beta_k\in (\Delta^-)^J \atop\beta_1+\cdots+\beta_k=\mu-\la } \phi_\mu(x_{\beta_1}\cdots x_{\beta_k}v_\la)  . \label{mkcontained}
\end{eqnarray}
For every $\beta\in (\Delta^-)^J$ we develop $x_\beta$ in the basis $y_{\beta i}\in\g_{\beta}$, $i=1,\,\ldots,\,m_\beta$, as follows:
\begin{eqnarray*}
 x_\beta =  \sum_{i=1}^{m_\beta}\ti{\psi}_{\beta i}(x_\beta)y_{\beta i} 
         =  \sum_{i=1}^{m_\beta}\ti{\psi}_{\beta i}(x)y_{\beta i}
         =  \sum_{i=1}^{m_\beta}\psi_{\beta i}(exp(x))y_{\beta i} .
\end{eqnarray*}
Inserting in the inner sumand of (\ref{mkcontained}) we get:
\begin{eqnarray*}
   \phi_\mu(x_{\beta_1}\cdots x_{\beta_k}v_\la)
   = \sum_{i_1,\,\ldots ,\,i_k}\psi_{\beta_1 i_1}(exp(x))\cdots \psi_{\beta_k i_k}(exp(x)) \,\phi_\mu(y_{\beta_1 i_1}\cdots y_{\beta_k i_k}v_\la) .
\end{eqnarray*}
Inserting back in (\ref{mkcontained}) we find
\begin{eqnarray}\label{uniringmatent}
   \phi_\mu(\ins v_\la)
    =  \phi_\mu(v_\la)1 + \sum_{k=1}^{ht(\la-\mu)}\frac{1}{k!}
       \sum_{\beta_1,\,\ldots,\,\beta_k\in (\Delta^-)^J,\,i_1,\,\ldots ,\,i_k \atop\beta_1+\cdots+\beta_k=\mu-\la } 
       \phi_\mu(y_{\beta_1 i_1}\cdots y_{\beta_k i_k}v_\la) \,\psi_{\beta_1 i_1}\cdots \psi_{\beta_k i_k}  \in \FK{(U_f^-)^J} . 
\end{eqnarray}

(b) Note that
\begin{eqnarray*}
   \Mklz{N\in P^+\,}{\,N(h_i)>0 \mb{ for all } i\in I\setminus J} = P^+\cap\bigcup_{K\subseteq J} F_K .
\end{eqnarray*} 
Fix $K\subseteq J$ and $\La\in P^+\cap F_K $. We get
\begin{eqnarray*}
  \mb{span}\left\{\phi(\ins v_N)\left|\phi\in L(N)^{(*)},\;N\in P^+\cap\overline{F_K}\right.\right\} \subseteq  \FK{(U_f^-)^J} 
\end{eqnarray*}
from (\ref{MkOint}). This span is even a subalgebra of $\FK{(U_f^-)^J}$. It contains $\delta_0(\ins v_0)$, which is the unit of $\FK{(U_f^-)^J}$. Furthermore, if  $\phi\in L(N)^{(*)}$, $N\in P^+\cap\overline{F_K} $ and $\psi\in L(M)^{(*)}$, $M\in P^+\cap\overline{F_K}$, then
\begin{eqnarray*}
   \phi(\ins v_N) \psi(\ins v_M)=(\phi\otimes \psi) (\ins (v_N\otimes v_M))=(\phi \prdca \psi)(\ins v_{N+M})
\end{eqnarray*}
with $\phi \prdca \psi\in L(N+M)^{(*)}$ and $N+M\in P^+\cap\overline{F_K}$.

Let $\mb{Alg}(\Mklz{(\pi(x)\delta_\La)(\ins v_\La)}{x\in \n^J})$ be the subalgebra of $\FK{(U_f^-)^J}$ generated by $\Mklz{(\pi(x)\delta_\La)(\ins v_\La)}{x\in \n^J}$. We have shown the inclusions
\begin{eqnarray*}
        \mb{Alg}\left( \Mklz{(\pi(x)\delta_\La)(\ins v_\La)}{x\in \n^J} \right)
           \subseteq  \mb{span}\left\{\phi(\ins v_N)\left|\phi\in L(N)^{(*)},\;N\in P^+\cap\overline{F_K}\right.\right\}
           \subseteq \FK{(U_f^-)^J}.
\end{eqnarray*}

(c) Now we show 
\begin{eqnarray*}
  \psi_{-\al i}\in \mb{Alg}\left( \Mklz{(\pi(x)\delta_\La)(\ins v_\La)}{x\in \n^J} \right)\quad \mb{ for all }\quad\al\in (\Delta^+)^J, \quad i=1,\,\ldots,\,m_{-\al},
\end{eqnarray*} 
which implies $\FK{(U_f^-)^J}\subseteq\mb{Alg}\left( \Mklz{(\pi(x)\delta_\La)(\ins v_\La)}{x\in \n^J} \right)$.

We first specialize formula (\ref{uniringmatent}) to the matrix coefficient $(\pi(z_{\al i})\delta_{\La})(\ins v_\La)$: Because $z_{\al i}$ kills $v_\La$ we have
\begin{eqnarray*} 
       (\pi(z_{\al i})\delta_{\La})(v_\La) = \delta_\La(z_{\al i}v_\La)=0.
\end{eqnarray*}
The summand for $k=1$ of the sum in (\ref{uniringmatent}) transforms as follows:
\begin{eqnarray*}
   \sum_{j}\delta_\La(z_{\al i}y_{-\al j}v_\La)\psi_{-\al j} = \sum_{j}\delta_\La( [z_{\al i},\, y_{-\al j}]v_\La)\psi_{-\al j} 
= \sum_{j}\delta_\La( \iB{z_{\al i}}{y_{-\al j}} \nu^{-1}(\al)v_\La)\psi_{-\al j} = 
    \iB{\La}{\al}\psi_{-\al i}.
\end{eqnarray*}
We have found:
\begin{eqnarray}
  && (\pi(z_{\al i})\delta_{\La})(\ins v_\La)\nonumber\\
  &&  = \iB{\La}{\al}\psi_{-\al i}
     + \sum_{k=2}^{ht(\al)}\frac{1}{k!}
       \sum_{\beta_1,\,\ldots,\,\beta_k\in (\Delta^-)^J,\,i_1,\,\ldots ,\,i_k \atop\beta_1+\cdots\beta_k=-\al } 
       \delta_\La( z_{\al i}y_{\beta_1 i_1}\cdots y_{\beta_k i_k}v_\La) \psi_{\beta_1 i_1}\cdots \psi_{\beta_k i_k} . \label{uniringind1}
\end{eqnarray}
Here $\al$ is of the form $\al=\sum_{i\in I}n_i\al_i$ with $n_i\in\Nn$ and $n_{j_0}\neq 0$ for at least one $j_0\in J$. With $\nu^{-1}(\al_i)=\frac{1}{2}\iB{\al_i}{\al_i}h_i$ and $\iB{\al_i}{\al_i}>0 $, $i\in I$, we obtain
\begin{eqnarray}\label{uniringind2}
  \iB{\La}{\al} = \La(\nu^{-1}(\al))=\sum_{i\in I,\,i\neq j_0}\underbrace{n_i}_{\geq 0} \frac{\iB{\al_i}{\al_i}}{2}\underbrace{\La(h_i)}_{\geq 0} +
  \underbrace{n_{j_0}}_{>0} \frac{  \iB{\al_{j_0}}{\al_{j_0}}  }{2}  \underbrace{\La(h_{j_0})}_{>0} > 0.
\end{eqnarray}

Let $\al\in (\Delta^+)^J$ with $\mb{ht}(\al)=1$. Then $m_\al = m_{-\al} =1$ and from (\ref{uniringind1}) and (\ref{uniringind2}) we get
\begin{eqnarray*}
   \psi_{-\al 1}=\frac{1}{ \iB{\La}{\al}}(\pi(z_{\al 1})\delta_{\La})(\ins v_\La)\in \mb{Alg}\left( \Mklz{(\pi(x)\delta_\La)(\ins v_\La)}{x\in \n^J} \right).
\end{eqnarray*}
Suppose $\psi_{-\beta i}\in \mb{Alg}\left( \Mklz{(\pi(x)\delta_\La)(\ins v_\La)}{x\in \n^J} \right)$ for all $\beta\in (\Delta^+)^J$ with $\mb{ht}(\beta)<p\in\N $ and $i=1,\,\ldots,\,m_{-\beta}$. Let $\al\in (\Delta^+)^J$  with $\mb{ht}(\al)=p$. Then from (\ref{uniringind1}) and (\ref{uniringind2}) we get
\begin{eqnarray*}
\psi_{-\al i}\in \mb{Alg}\left( \Mklz{(\pi(x)\delta_\La)(\ins v_\La)}{x\in \n^J} \right)\quad \mb{ for all } i=1,\,\ldots,\,m_\al.
\end{eqnarray*}

(d) At last we show that the linear map $\n^J\to \FK{(U_f^-)^J}$ of the theorem is injective and that $\FK{(U_f^-)^J} = \mb{Alg}\left( \Mklz{(\pi(x)\delta_\La)(\ins v_\La)}{x\in \n^J} \right)$ is a symmetric algebra in $\Mklz{(\pi(x)\delta_\La)(\ins v_\La)}{x\in \n^J}$.

Choose a counting $(\beta_m,\,i_m)$, $m\in\N$, of the set
\begin{eqnarray*}
  \Mklz{ (\beta,i)\, }{ \,\beta\in (\Delta^+)^J,\;i=1,\,\ldots,\,m_\beta }.
\end{eqnarray*}
Since $z_{\beta_m i_m}$, $m\in\N$, is a base of $\n^J$ it is sufficient to show that the monomials in
\begin{eqnarray*}
    f_m:=\left(\pi(z_{\beta_m i_m})\delta_\La\right)(\ins v_\La),\quad m\in\N,
\end{eqnarray*}
are linearly independent.

Note that by (\ref{uniringind1}) and (\ref{uniringind2}) for all $m\in\N$ we have
\begin{eqnarray}\label{btuniring}
    f_m = \underbrace{ \iB{\La}{\beta_m} }_{\neq 0} \psi_{-\beta_m i_m} + r_m\,,
\end{eqnarray}
where $r_m$ is a sum of homogeneous components of degree $\geq 2$. 

Now let $r\in\N$ and $L\in\Nn$ and suppose that
\begin{eqnarray}\label{uniringlk}
  0 = \sum_{k_1,\,\ldots,\,k_r \in\Nn  \atop  k_1+\cdots + k_r\leq L} c_{k_1\cdots k_r} f_1^{k_1}\cdots f_r^{k_r} 
\end{eqnarray}
with coefficients $c_{k_1\cdots k_r}\in\F$. 

First we compare the 0-homogeneous components of both sides of equation (\ref{uniringlk}). Because of (\ref{btuniring}) we get 
$c_{0\cdots 0} = 0$. Let $l\in\N$, $l\leq L$, and suppose we have shown
\begin{eqnarray*}
     c_{k_1\cdots k_r} = 0 \quad\mb{ for all }\quad k_1,\,\ldots,\,k_r \in\Nn \;\mb{ with }\; k_1+\cdots + k_r\leq l-1.
\end{eqnarray*}
Now we compare the $l$-homogeneous components of both sides of equation (\ref{uniringlk}). Because of (\ref{btuniring}) we get 
\begin{eqnarray*}
    0 =  \sum_{k_1,\,\ldots,\,k_r \in\Nn  \atop  k_1+\cdots + k_r  = l } c_{k_1\cdots k_r} 
               \iB{\La}{\beta_1}^{k_1}\cdots \iB{\La}{\beta_r}^{k_r}  \psi_{-\beta_1 i_1}^{k_1}\cdots \psi_{-\beta_r i_r}^{k_r}  .
\end{eqnarray*}
Since the monomials in $\psi_{-\beta_m i_m}$, $m\in\N$, are linearly independent and $\iB{\La}{\beta_m}\neq 0$ for all $m\in\N$ we obtain
\begin{eqnarray*}
    c_{k_1\cdots k_r} = 0 \quad\mb{ for all }\quad k_1,\,\ldots,\,k_r \in\Nn \;\mb{ with }\; k_1+\cdots+ k_r=l.
\end{eqnarray*}
\qed

Next we make some observations for the groups $(U^-)^J$ and its coordinate rings  $\FK{(U_f^-)^J}$, $J\subseteq I$.
\begin{cor}\label{denseU} Let $J\subseteq I$. The subgroup $(U^-)^J$ of $(U^-_f)^J$ is dense. In particular, its coordinate ring 
\begin{eqnarray*}
  \FK{(U^-)^J}:= \FK{(U_f^-)^J}\res{(U^-)^J}
\end{eqnarray*}
is isomorphic to $\FK{(U_f^-)^J}$ by the restriction map.
\end{cor}
\Proof This holds trivially for $J=I$. Let $J\neq I$. Let $f\in \FK{(U_f^-)^J}$ such that $f\res{(U^-)^J}=0$. By the previous theorem we can write $f$ in the form
\begin{eqnarray*}
  f=\sum_{N\in P^+\cap \overline{F_J}}\phi_N(\ins v_N)
\end{eqnarray*}
with $\phi_N\in L(N)^{(*)}$, $\phi_N\neq 0$ for at most finitely many $N\in P^+\cap\overline{F_J}$. Because $U_J^-$ stabilizes $L(N)_N$ pointwise for all $N\in P^+\cap\overline{F_J}$ and because of $U^-=(U^-)^J U_J^-$ we find
\begin{eqnarray*}
  0 = \sum_{N\in P^+\cap \overline{F_J}}\phi_N(u v_N) \quad\mb{ for all }\quad u\in U^-\,.
\end{eqnarray*}
For $p\in\N$, $\beta_1,\,\ldots,\,\beta_p\in \Delta^-_{re}$, and $x_{\beta_1}\in\g_{\beta_1}$, \ldots, $x_{\beta_p}\in\g_{\beta_p}$ we have
\begin{eqnarray*}
  0 = \sum_{N\in P^+\cap \overline{F_J}}\phi_N(exp(t_1 x_{\beta_1})\cdots exp(t_1 x_{\beta_1}) v_N) 
    = \sum_{k_1,\,\ldots,\,k_p\in\Nn}\frac{t_1^{k_1}\cdots t_p^{k_p}}{k_1!\cdots k_p!}
         \sum_{N\in P^+\cap \overline{F_J}}\phi_N(x_{\beta_1}^{k_1}\cdots x_{\beta_p}^{k_p} v_N)
\end{eqnarray*}
for all $t_1,\,\ldots,\,t_p\in\F$. The last expression is well defined and polynomial in $t_1,\,\ldots,\,t_p$ because $x_{\beta_1}$, \ldots, $x_{\beta_p}$ act locally nilpotent. The coefficients of the polynomial vanish. In particular, we find
\begin{eqnarray*}
    0=\sum_{N\in P^+\cap \overline{F_J}}\phi_N(v_N) \quad\mb{ and }\quad  0=\sum_{N\in P^+\cap \overline{F_J}}\phi_N( x_{\beta_1}\cdots x_{\beta_p}v_N)\,. 
\end{eqnarray*}
By \S 1.3 of \cite{K} the Lie algebra $\n^-$ is generated by the elements $f_i\in\g_{-\al_i}$, $i\in I$. Therefore, also the universal enveloping algebra $U(\n^-)$ is generated by $f_i\in\g_{-\al_i}$, $i\in I$. Because of $-\al_i\in \Delta^-_{re}$, $i\in I$, we have found
\begin{eqnarray}\label{univ0}
  0 = \sum_{N\in P^+\cap \overline{F_J}}\phi_N(x v_N) \quad\mb{ for all }\quad x\in U(\n^-)\,.
\end{eqnarray}
Now let $y=\sum_{\beta\in(\Delta^-)^J}y_\beta \in (\n_f^-)^J$. With (\ref{univ0}) we get
\begin{eqnarray*}
  && f(exp(y)) = \sum_{N\in P^+\cap \overline{F_J}}\phi_N(exp(y)v_N) \\
  && = \sum_{N\in P^+\cap \overline{F_J}}\phi_N(v_N) 
      + \sum_{\gamma\in Q^- } \frac{1}{k!}\sum_{k=1}^{-ht(\gamma)}\sum_{\beta_1,\,\ldots,\,\beta_k\in (\Delta^-)^J\atop \beta_1+\cdots+\beta_k=\gamma} \sum_{N\in P^+\cap \overline{F_J}}\phi_N(y_{\beta_1}\cdots y_{\beta_k }v_N)\,=\,0\,. 
\end{eqnarray*}
We have shown $f=0$.
\qed

Recall from Section \ref{Preli} that the Kac-Moody group $G$ is equipped with a coordinate ring of matrix coefficients such that
a Peter-Weyl theorem holds: A $G^{op}\times G$-equivariant bijective linear map
\begin{eqnarray*}
    PW:\;\bigoplus_{\La\in P^+}L(\La)^{(*)}\otimes L(\La) \to \FK{G}   
\end{eqnarray*}
is given by $PW(\phi\otimes v):=\phi(\ins v)$ for all $\phi\in L(\La)^{(*)}$, $v\in L(\La)$, $\La\in P^+$. Here the matrix coefficient $\phi(\ins v)$ denotes the function which assigns to $g\in G$ the value $\phi(gv)\in\F$. To keep our notation easy we often denote its restriction to nonempty subsets of $G$ by the same symbol.

Now $(U^-)^J$ is also contained in $G$. It holds:
\begin{cor}\label{UniRes} Let $J\subseteq I$. The restriction of the coordinate ring $\FK{G}$ to $(U^-)^J$ coincides with the restriction of the coordinate ring $\FK{(U_f^-)^J}$ to $(U^-)^J$, which we denoted by $\FK{(U^-)^J}$.
\end{cor}
\Proof The Peter-Weyl theorem implies
\begin{eqnarray*}
  \FK{G}\res{U_J^-}= \mb{span}\left\{\phi(\ins v)\res{U_J^-}\left|\phi\in L(\La)^{(*)},\;v\in L(\La),\;\La\in P^+\right.\right\}.
\end{eqnarray*}
Now the corollary follows from Theorem \ref{UniRing}.
\qed

To avoid trivialities let $J\subsetneqq I$. The pairing $\iBl:\n^J\times (\n^-)^J\to\F$ is nondegenerate. Therefore, we get a linear bijective map
\begin{eqnarray*}
\ti{\psi}: \n^J &\to & ((\n^-)^J)^{(*)}\\
        x \;\, &\mapsto & \quad\ti{\psi}_x
\end{eqnarray*}
where $\ti{\psi}_x\in ((\n^-)^J)^{(*)}\subseteq ((\n_f^-)^J)^{*}$ is defined by
\begin{eqnarray*}
  \ti{\psi}_x(y) := \sum_{\al\in (\Delta^-)^J}\iB{x}{y_\al}\quad\mb{ for }\quad y=\sum_{\al\in (\Delta^-)^J}y_\al\in (\n_f^-)^J.
\end{eqnarray*}
Set $\psi_x := (exp^*)^{-1}(\ti{\psi}_x) \in \FK{ (U_f^-)^J }$ for every $x\in \n^J$. By definition it is
\begin{eqnarray*} 
      \FK{ (U_f^-)^J } = \mb{Sym}( \;\Mklz{\psi_x}{x\in \n^J}\;)\,.
\end{eqnarray*}
Fix $\La\in P^+$ with $\La(h_i)>0$ for all $i\in I\setminus J$. By Theorem \ref{UniRing} we have
\begin{eqnarray*}
      \FK{ (U_f^-)^J } = \mb{Sym}( \;\Mklz{(\pi(x)\delta_\La)(\ins v_\La)}{x\in \n^J}\;).
\end{eqnarray*}
The linear map which maps $\psi_x$ to $(\pi(x)\delta_\La)(\ins v_\La)$, $x\in\n^J$, induces an automorphism of the algebra $\FK{ (U_f^-)^J }$. The next theorem describes this automorphism as a comorphism.

\begin{thm}\label{UniIso} Let $J\subsetneqq I$ and $\La\in P^+$ with $\La(h_i)>0$ for all $i\in I\setminus J$. Set $h:=\nu^{-1}(\La)\in\h$. 
\begin{itemize}
\item[(a)] We get a bijective map $d\,\Xi_h: (\n_f^-)^J \to (\n_f^-)^J $ by
\begin{eqnarray*}
  (d\,\Xi_h)(x):= Ad(exp(-x))h-h = \sum_{k\in\N}\frac{1}{k!}(-ad(x))^k h\, ,\quad x\in (\n_f^-)^J .
\end{eqnarray*}
\item[(b)] The map
\begin{eqnarray*}
  \begin{array}{rccc}
      \Xi_h: & (U_f^-)^J &\to &(U_f^-)^J\\
              & u & \mapsto & exp(Ad(u^{-1})h-h)
  \end{array}
\end{eqnarray*}
is bijective. Its comorphism $\Xi_h^*:\FK{(U_f^-)^J} \to \FK{(U_f^-)^J}$ exists. It is an isomorphism of algebras for which it holds
\begin{eqnarray*}
    \Xi_h^*(\psi_x) = (\pi(x)\delta_\La)(\ins v_\La)\,,\quad x\in\n^J.
\end{eqnarray*}
\end{itemize}
\end{thm}

\Proof To (a): Let $x=\sum_{\beta\in(\Delta^+)^J}x_{-\beta}\in (\n_f^-)^J $. For $\gamma \in \Delta\cup\{0\}$ the $(-\gamma)$-homogeneous component of the defining sum of $(d\,\Xi_h)(x)$ is given by
\begin{eqnarray*}
    \sum_{k\in\N}\frac{1}{k!}\sum_{\beta_1,\,\beta_2,\,\ldots,\,\beta_k\in(\Delta^+)^J\atop \beta_1+\beta_2+\cdots+\beta_k =\gamma}
      (-ad\, x_{-\beta_1})(-ad\, x_{-\beta_2})\cdots(-ad\, x_{-\beta_k})h\,,
\end{eqnarray*}
a sum over the empty set defined to be zero. For $\gamma \notin (\Delta^+)^J$ this expression is zero. For $\gamma \in (\Delta^+)^J$ it is a well defined expression: The roots in $(\Delta^+)^J$ have at least height 1. Therefore, it is sufficient to sum from $k=1$ to $k= \mb{ht}(\gamma)$ in the first sum. 

Let $\gamma \in (\Delta^+)^J$ and write this component as follows:
\begin{eqnarray}\label{dbco1}
   (d\,\Xi_h)(x)_{-\gamma} =  -\gamma(h)x_{-\gamma} + \sum_{k\in\N,\,k\geq 2}\frac{(-1)^k}{k!}\sum_{\beta_1,\,\beta_2,\,\ldots,\,\beta_k\in(\Delta^+)^J\atop \beta_1+\beta_2+\cdots+\beta_k =\gamma}
       [x_{-\beta_1},[x_{-\beta_2},[\cdots [x_{-\beta_k},h]\cdots ]]]\,.
\end{eqnarray}
Here $\gamma$ is of the form $\gamma=\sum_{i\in I}n_i\al_i$ with $n_i\in\Nn$ and $n_{j_0}\neq 0$ for at least one $j_0\in J$. With $\nu^{-1}(\al_i)=\frac{1}{2}\iB{\al_i}{\al_i}h_i$ and $\iB{\al_i}{\al_i}>0 $, $i\in I$, we obtain
\begin{eqnarray}\label{dbco2}
  \gamma(h)=\gamma(\nu^{-1}(\La))=\iB{\gamma}{\La}=\La(\nu^{-1}(\gamma))=\sum_{i\in I,\,i\neq j_0}\underbrace{n_i}_{\geq 0} \frac{\iB{\al_i}{\al_i}}{2}\underbrace{\La(h_i)}_{\geq 0} +
  \underbrace{n_{j_0}}_{>0} \frac{  \iB{\al_{j_0}}{\al_{j_0}}  }{2}  \underbrace{\La(h_{j_0})}_{>0} > 0.
\end{eqnarray}

Let $y=\sum_{\beta\in(\Delta^+)^J}y_{-\beta}\in (\n_f^-)^J $. We show that there exists a unique $x=\sum_{\beta\in(\Delta^+)^J}x_{-\beta}\in (\n_f^-)^J $ such that $(d\,\Xi_h)(x)= y$ or equivalently
\begin{eqnarray*}
  (d\,\Xi_h)(x)_{-\gamma} = y_{-\gamma}\;\;\mb{ for all }\gamma\in (\Delta^+)^J.
\end{eqnarray*}
We do this by induction over $\mb{ht}(\gamma)$. If $\mb{ht}(\gamma)=1$, then from (\ref{dbco1}) and (\ref{dbco2}) we get $x_{-\gamma}=\frac{1}{-\gamma(h)}y_{-\gamma}$ uniquely.
Now let $\mb{ht}(\gamma) =m+1$, $m\in\N$,  and suppose that $x_{-\gamma'}$ has been determined uniquely for all $\gamma'\in (\Delta^+)^J$ with $\mb{ht}(\gamma')\leq m $. Then from (\ref{dbco1}) and (\ref{dbco2}) we get
 \begin{eqnarray*}
   x_{-\gamma} =  \frac{1}{-\gamma(h)}y_{-\gamma} - \frac{1}{-\gamma(h)}\sum_{k\in\N,\,k\geq 2}\frac{(-1)^k}{k!}\sum_{\beta_1,\,\beta_2,\,\ldots,\,\beta_k\in(\Delta^+)^J\atop \beta_1+\beta_2+\cdots+\beta_k =\gamma}
       [x_{-\beta_1},[x_{-\beta_2},[\cdots [x_{-\beta_k},h]\cdots ]]]
\end{eqnarray*}
uniquely.

To (b): The map $d\,\Xi_h: (\n_f^-)^J \to (\n_f^-)^J $ and the map $exp:(\n^-_f)^J\to (U^-_f)^J$ are bijective. Therefore, by its definition also the map $\Xi_h:  (U_f^-)^J \to (U_f^-)^J$ is bijective. Since $\psi_x$, $x\in\n^J$, as well as $(\pi(x)\delta_\La)(\ins v_\La)$, $x\in\n^J$, generate the algebra $\FK{U_f^-}$ it only remains to show  $\psi_x\circ \Xi_h = (\pi(x)\delta_\La)(\ins v_\La)$ for $x\in\n^J$.

For every $y\in (\n_f^-)^J $ we have
\begin{eqnarray}
&&  \psi_x(\Xi_h(exp(y))) = \Big(  (exp^*)^{-1}\ti{\psi}_x  \Big)  \Big(exp\big((d\,\Xi_h)(y)\big)\Big) = \ti{\psi}_x \Big((d\,\Xi_h)(y)\Big)=\ti{\psi}_x\Big(\sum_{k\in\N}\frac{1}{k!}(-ad(y))^k h\Big) \nonumber\\
 &&  = \sum_{\gamma\in(\Delta^+)^J} \sum_{k\in\N}\frac{(-1)^k}{k!}\sum_{\beta_1,\,\beta_2,\,\ldots,\,\beta_k\in(\Delta^+)^J\atop \beta_1+\beta_2+\cdots+\beta_k =\gamma}
  \ti{\psi}_x\left( [y_{-\beta_1},[y_{-\beta_2},[\cdots [y_{-\beta_k},h]\cdots ]]]\right)   \,. \label{mkunipotentb}
\end{eqnarray}
We transform the inner summand of this expression:
\begin{eqnarray*}
 && \ti{\psi}_x\left( [y_{-\beta_1},[\cdots [y_{-\beta_k},h]\cdots ]]\right)
   =  \iB{x}{ [y_{-\beta_1},[\cdots [y_{-\beta_k},h]\cdots ]] } 
   =  \iB{ [y_{-\beta_k},[\cdots [y_{-\beta_1},x] \cdots ]]}{h} \\
  && = \iB{ [y_{-\beta_k},[\cdots [y_{-\beta_1},x] \cdots ]]_0}{h} \,. 
\end{eqnarray*}
Here we used $\iB{\g_\al}{\h}=\{0\}$ for all $\al\in\Delta$ and $[y_{-\beta_k},[\cdots [y_{-\beta_1},x] \cdots ]]_0$ denotes the 0-homogeneous component of $[y_{-\beta_k},[\cdots [y_{-\beta_1},x] \cdots ]]$. Now we insert $h=\nu^{-1}(\La)$ and transform further:
\begin{eqnarray*}
  &&  \La( [y_{-\beta_k},[\cdots [y_{-\beta_1},x] \cdots ]]_0 ) 
       =  \delta_\La( [y_{-\beta_k},[\cdots [y_{-\beta_1},x]\cdots ]]_0 v_\La )
       =  \delta_\La( [y_{-\beta_k},[\cdots [y_{-\beta_1},x]\cdots ]] v_\La )   \\
  && = \left(\pi([y_{-\beta_k},[\cdots [y_{-\beta_1},x]\cdots ]]) \delta_\La\right)(  v_\La )
     = (-1)^k ([\pi(y_{-\beta_k}),[\cdots [\pi(y_{-\beta_1}),\pi(x)]\cdots ]] \delta_\La)(  v_\La )\,.
\end{eqnarray*}
Since $\pi(y_{-\beta_1})$, \ldots, $\pi(y_{-\beta_k})$ kills $\delta_\La$ it only remains
\begin{eqnarray*}
    (-1)^k \big( \pi(y_{-\beta_k})\cdots\pi(y_{-\beta_1})\pi(x) \delta_\La \big)(  v_\La ) 
  = (-1)^k (\pi(x) \delta_\La) (y_{-\beta_1}\cdots y_{-\beta_k}v_\La).
\end{eqnarray*}
Inserting in (\ref{mkunipotentb}) we get
\begin{eqnarray*}
  \psi_x(\Xi_h(exp(y))) 
         = (\pi(x) \delta_\La)\Big(\sum_{\gamma\in(\Delta^+)^J}\sum_{k\in\N}\frac{1}{k!}  \sum_{\beta_1,\,\beta_2,\,\ldots,\,\beta_k\in(\Delta^+)^J\atop \beta_1+\beta_2+\cdots+\beta_k =\gamma} y_{-\beta_1}y_{-\beta_2}\cdots y_{-\beta_k}v_\La\Big) \,.
\end{eqnarray*}
Here the sum over $\gamma \in (\Delta^+)^J$ can be replaced by a sum over $\gamma\in Q^+$, because of $\pi(x)\delta_\La\in \bigoplus_{\la\in\La-(\Delta^+)^J}L(\La)_\la^*$. If we use in addition $(\pi(x)\delta_\La)(v_\La)=\delta_\La(xv_\La)=0$ we find
\begin{eqnarray*}
 \psi_x(\Xi_h(exp(y))) = (\pi(x) \delta_\La)\Big(\sum_{k\in\N}\frac{1}{k!} y^k v_\La\Big) = (\pi(x) \delta_\La) (exp(y) v_\La)\,.
\end{eqnarray*}
\qed

\begin{rem} (a) V. Kac and D. Peterson described in Lemma 4.3 of \cite{KP2} the coordinate ring $\FK{G}\res{U}$ as a symmetric algebra in the matrix coefficients $\delta_\La (\ins x v_\La)$, $x\in\n^-$, where $\La\in P^+\cap C$ is arbitrary, the key steps as follows: Let $\n$ be equipped with the coordinate ring $\FK{\n}=\mb{Sym}(\bigoplus_{\al\in\Delta^+}\g_\al^*)$. Let $h\in\h$ such that $\al(h)\neq 0$ for all $\al\in\Delta$. V. Kac and D. Peterson
showed that the map $\psi:U\to\n$ obtained by $\psi(u)=Ad(u)h-h$, $u\in U$, induces a bijective comorphism $\psi^*:\FK{\n}\to \FK{G}\res{U}$.
This is tricky in the context of the minimal Kac-Moody group, in particular: The group $U$ is only defined by generators. The map $\psi$ is in general not surjective, and it is not at all obvious that its image is dense. 

The author tried in \cite{M1}, Theorem 5.6, to adapt the proof of V .Kac and D. Peterson to describe for $J\subseteq I$ the coordinate ring $\FK{G}\res{U^J}$, but only reached that $\FK{G}\res{U^J}$ is generated by the matrix coefficients $\delta_\La(\ins x v_\La)$, $x\in (\n^-)^J$, where $\La\in P^+\cap F_J$ is arbitrary. 

As shown above, the descriptions of these coordinate rings can be reached more directly in the context of formal Kac-Moody groups.

(b) The proofs of Theorem \ref{UniRing} and Theorem \ref{UniIso} work quite generally: It is easy to modify Theorem \ref{UniRing} such that it holds in the case where $(\Delta^-)^J$ is replaced by a subset $\Omega\subseteq\Delta^-$. 
It is easy to modify Theorem \ref{UniIso} such that it holds in the case where $(\Delta^-)^J$ is replaced by a bracket closed subset $\Omega\subseteq\Delta^-$. (Bracket closed means: If $\al,\,\beta\in\Omega$ and $\al+\beta\in\Delta^-$, then also $\al+\beta\in\Omega$.)
\end{rem}

The following easy Proposition will be useful later.
\begin{prop}\label{SpMk} Let $J\subseteq I$. Let $\phi\in L_J(\La)^{(*)}$, $v\in L_J(\La)_f$, $\La\in P^+$. We have
\begin{eqnarray*}
  \phi(uv)=\phi(v) \quad\mb{ for all }\quad u\in (U_f^-)^J.
\end{eqnarray*}
\end{prop}
\Proof This holds trivially for $J=I$. Let $J\neq I$. Then $\phi\in \bigoplus_{\la\in \La-(Q_J)_0^+}L(\La)^*_\la$ and $v \in \prod_{\la\in \La-(Q_J)_0^+}L(\La)_\la$. Furthermore, $u\in (U_f^-)^J$ acts on $v$ as $exp(x)$ for some $x\in\prod_{\al\in (\Delta^-)^J}\g_\al$. We find
\begin{eqnarray*}
   uv - v \in \prod_{\la\in \La-(Q_J)_0^+ - Q^+_0\setminus (Q_J)^+_0}L(\La)_\la .
\end{eqnarray*}
Therefore, $\phi(uv-v)=0$.
\qed

For $J\subseteq I$ set
\begin{eqnarray*}
    CA_J:=\bigoplus_{\La\in P^+_J}L_J(\La)^{(*)},
\end{eqnarray*}
which is a $P^+_J$-graded subalgebra of $CA$ by Proposition \ref{smdu} (a). In particular,  $CA_\emptyset=\F\,1$.
\begin{prop} For $\emptyset\neq J\subseteq I$ the algebra $CA_J$ is the Cartan algebra associated to $G_J$, the highest weight vectors chosen as $v_\La\in L_J(\La)_\La$, $\La\in P_J^+$.
\end{prop}
\Proof For every $\La\in P^+$ the $G$-orbit $G v_\La$ spans $L(\La)$, because $L(\La)$ is an irreducible $G$-module. Now let $\La,\, N\in P^+$ and $\phi\in L(\La)^{(*)}$, $\psi\in L(N)^{(*)}$. Then the Cartan product $\phi \prdca\psi\in L(\La+N)^{(*)}$ is determined by
\begin{eqnarray}\label{CPor}
  (\phi \prdca \psi)(g v_{\La+N})= (\phi\otimes\psi)(g(v_\La\otimes v_N))= (\phi\otimes\psi)(g v_\La\otimes g v_N)=\phi(g v_\La)\psi(g v_N)\,, \quad g\in G .
\end{eqnarray} 

The subgroup $G_J$ of $G$ identifies with the derived Kac-Moody group $G(A_J)'$ to the generalized Cartan submatrix $A_J$ of $A$. The sublattice $P_J$ of $P$ identifies with the weight lattice, and $P_J^+$ identifies with the set of dominant weights of the weight lattice of $G(A_J)'$. For $\La\in P^+_J$ the $G_J$-module $L_J(\La)$ identifies by Theorem \ref{sm} with the irreducible highest weight module of $G(A_J)'$, whose highest weight corresponds to $\La$.

Let $\La,\, N\in P^+_J$ and $\phi\in L_J(\La)^{(*)}$, $\psi\in L_J(N)^{(*)}$ Then $\phi \prdca \psi\in L_J(\La+N)^{(*)}$ by Proposition \ref{smdu} (a). Restricting (\ref{CPor}) to $g\in G_J$ shows: The product $\prdca$ is also the Cartan product of $CA_J$, where the highest weight vectors have been chosen as $v_\La\in L_J(\La)_\La$, $\La\in P_J^+$.
\qed

For $J\subseteq I$ let 
\begin{eqnarray*}
     \FK{P\cap(\overline{F_J}-\overline{F_J})} =\bigoplus_{\La\in P\cap(\overline{F_J}-\overline{F_J})}\F e_{\La}
\end{eqnarray*} 
be the group algebra of the lattice $P\cap(\overline{F_J}-\overline{F_J})$, equipped with its natural $P\cap(\overline{F_J}-\overline{F_J})$-gradation.
\begin{thm}\label{Sp4} Let $J\subseteq I$. 
\begin{itemize}
\item[(a)] There exists an isomorphism of graded algebras
\begin{eqnarray*}
 \Gamma: (\widetilde{P^+\cap\overline{F_J}})^{-1}CA \to \FK{(U^-_f)^J} \otimes CA_J \otimes \FK{P\cap(\overline{F_J}-\overline{F_J})}
\end{eqnarray*}
such that
\begin{eqnarray}
 \Gamma(\frac{\delta_N}{\delta_\La}) = 1 \otimes 1 \otimes   e_{N-\La}  \;\: &\;\mb{ for all }\;&  N,\,\La\in P^+\cap\overline{F_J},\label{mg1}\\
 \Gamma(\phi) = 1 \otimes\phi \otimes   1  \qquad\quad &\;\mb{ for all }\;&    \phi\in L_J(N)^{(*)},\;N\in P_J^+,\label{mg2}\\
 \Gamma(\frac{\phi}{\delta_N}) = \phi( \ins v_N) \otimes 1 \otimes 1 &\;\mb{ for all }\;&   \phi\in L(N)^{(*)},\;N\in P^+\cap\overline{F_J}. \label{mg3}
\end{eqnarray}
\item[(b)] The map
\begin{eqnarray*}
  \alpha:\,\mb{Proj\,}(\F[(U_f^-)^J]\otimes CA_J) &\to     & \qquad \qquad\qquad D(J)\\  
       Q \qquad\qquad\;   &\mapsto &  \Gamma^{-1}\big( Q\otimes \FK{P\cap (\overline{F_J}-\overline{F_J})} \big)\cap CA
\end{eqnarray*}
is a homeomorphism, mapping $\mb{Proj\,}(\F[(U_f^-)^J]\otimes CA_J)(\F)$ bijectively to $D(J)(\F)$.
\end{itemize}
\end{thm}

\Proof Part (b) of the theorem follows from part (a) and from Proposition \ref{Sp3} and Proposition \ref{prodlattice}. To show part (a) we first interpret and describe in ($\alpha $), ($\beta $), ($\gamma $), ($\delta $) the algebras involved as coordinate rings of functions. In ($\epsilon $) we realize $\Gamma$ as a comorphism. By this approach, $\Gamma$ automatically is a well defined map. It also shows how $\Gamma$ can be found.

($\alpha $)  The functions of $\F[G]^U$ take constant values on the elements of every coset $gU$, $g\in G$. Therefore, they induce a coordinate ring of functions on $G/U$ isomorphic to $\F[G]^U$. We denote this coordinate ring on $G/U$ also by $\F[G]^U$.
For $\La\in P^+$ set 
\begin{eqnarray*} 
  \theta_\La:=\delta_\La(\ins v_\La)\in \F[G]^U.
\end{eqnarray*} 
Define
\begin{eqnarray*}
 D_{G/U}(J):=  \Mklz{ gU\in G/U }{ \theta_\La(gU)\neq 0 \mb{ for all }\La\in P^+\cap\overline{F_J} },
\end{eqnarray*}
and equip this set with the coordinate ring of functions
\begin{eqnarray*}
    \FK{ D_{G/U}(J) } := \Mklz{\frac{f}{\theta_\La}\res{D_{G/U}(J)} }{f\in \F[G]^U,\;\La\in P^+\cap\overline{F_J}}.
\end{eqnarray*}

The Borel-Weil isomorphism of algebras
\begin{eqnarray*}
  BW:\; CA=\bigoplus_{\La\in P^+}L(\La)^{(*)}\to \F[G]^U
\end{eqnarray*}
is given by $BW(\phi):=\phi(\ins v_\La)$ for all $\phi\in L(\La)^{(*)}$, $\La\in P^+$, compare Section \ref{Preli} and Subsection \ref{subsectionCA}.
It is easy to check that it induces a well defined surjective homomorphism of algebras
\begin{eqnarray}\label{BW2}
   BW^J:\;\widetilde{(P^+\cap\overline{F_J})}^{-1}CA  \to    \FK{ D_{G/U}(J) }
\end{eqnarray}
by $BW^J( \frac{\phi}{\delta_\La} ):=\frac{BW(\phi)}{\theta_\La}\res{D_{G/U}(J)}$ for all $\phi\in CA$, $\La\in P^+\cap\overline{F_J}$. This map is also injective, because $D_{G/U}(J)$ is dense in $G/U$, which can be seen as follows: Let $N\in P\cap C$. By the Birkhoff decomposition of $G$ we find
\begin{eqnarray*}
 \Mklz{ gU\in G/U }{ \theta_N(gU)\neq 0 } = U^-TU/U \subseteq D_{G/U}(J).
\end{eqnarray*}
Now $\F[G]^U$ is an integral domain, compare Section 1. Therefore, the principal open set on the left and also $D_{G/U}(J)$ is dense in $G/U$.

($\beta $) We denote by $\FK{G_J}$ the restriction of $\FK{G}$ to $G_J$. We have $G_\emptyset=\{1\}$ by definition and $\FK{G_\emptyset}=\F 1$. For $J\neq\emptyset$ the subgroup $G_J$ of $G$ identifies with the derived Kac-Moody group $G(A_J)'$ to the generalized Cartan submatrix $A_J$ of $A$. It is well-known that $\FK{G_J}$ identifies with the algebra of strongly regular functions on $G(A_J)'$ of $\cite{KP2}$. (To check this use Theorem \ref{sm}.) 

We denote by $\F[G_J]^{U_J}$ the functions of $\FK{G_J}$ which take constant values on the elements of every coset $gU_J$, $g\in G_J$. If $J\neq\emptyset$ then by  Corollary 2.2 and Remark 2.1 of \cite{KP2} there holds a Borel-Weil theorem : An isomorphism of algebras
\begin{eqnarray}\label{BW1}
  BW_J:\; CA_J=\bigoplus_{\La\in P_J^+}L_J(\La)^{(*)}\to\F[G_J]^{U_J}
\end{eqnarray}
is given by $BW_J(\phi):=\phi(\ins v_\La)$ for all $\phi\in L_J(\La)^{(*)}$, $\La\in P_J^+$. By our definitions this also holds for $J=\emptyset$. We identify $\F[G_J]^{U_J}$ with its corresponding coordinate ring of functions on $G_J/U_J$.

($\gamma $) It holds 
\begin{eqnarray*}
    P\cap(\overline{F_J}-\overline{F_J}) 
            = \Mklz{\la\in P}{ \la(h_i)=0 \mb{ for all } i\in J } 
            = \bigoplus_{i\in\{1,\,\ldots,\,2n-l\}\setminus J} \mathbb{Z}\La_i\, .
\end{eqnarray*}
The torus $T^J$ has been defined as the subtorus of $T$ generated by $t_{h_i}(\F^\times)$, $i\in \{1,\,2,\,\ldots,\,2n-l\}\setminus J$. Therefore, the group algebra $\FK{P\cap(\overline{F_J}-\overline{F_J})}$ identifies with the classical coordinate ring of the torus $T^J$.
Explicitly, the generator $e_\la\in \FK{P\cap(\overline{F_J}-\overline{F_J})}$, $\la\in P\cap(\overline{F_J}-\overline{F_J})$, identifies with the function
\begin{eqnarray*}
    e_\la\Big(\prod_{i\in\{1,\,\ldots,\,2n-l\}\setminus J} h_i(s_i)\Big) = \prod_{i\in\{1,\,\ldots,\,2n-l\}\setminus J} (s_i)^{\la(h_i)}, \mb{ where }s_i\in\F^\times, \;i\in\{1,\,\ldots,\,2n-l\}\setminus J\,.
\end{eqnarray*}

($\delta $) By Corollary \ref{denseU} the restriction map $\FK{(U^-_f)^J} \to \FK{(U^-)^J} $ is an isomorphism of algebras. Hence we get
\begin{eqnarray*}
  \FK{(U^-)^J} = \mb{span}\Mklz{\phi(\ins v_N)}{\phi\in L(N)^{(*)},\;N\in P^+\cap\overline{F_J} }
\end{eqnarray*}
by Theorem \ref{UniRing}. By Corollary \ref{UniRes} we also get $\FK{G}\res{(U^-)^J}=\FK{(U^-)^J}$.

($\epsilon $) Now we show: The map
\begin{eqnarray*}
  m: (U^-)^J \times G_J/U_J \times T^J &\to& D_{G/U}(J)\\
        (u^-,\,gU_J,\,t) &\mapsto & u^- g t U
\end{eqnarray*} 
is well defined and induces a comorphism
\begin{eqnarray*}
m^*:\, \FK{D_{G/U}(J)}\to \FK{(U^-)^J} \otimes \F[G_J]^{U_J}\otimes \FK{P\cap(\overline{F_J}-\overline{F_J})},
\end{eqnarray*}
which is an isomorphism of algebras. Furthermore, 
\begin{eqnarray}
   m^*\left(\frac{\theta_N}{\theta_\La}\res{D_{G/U}(J)}\right) = 1\otimes 1\otimes e_{N-\La} &\mb{ for }&  N,\,\La\in P^+\cap\overline{F_J}, \label{mgl1}\\
   m^*\left(\phi(\ins v_N)\res{D_{G/U}(J)}\right)= 1\otimes \phi(\ins v_N)\res{G_J/U_J}\otimes 1  &\mb{ for }&  \phi\in L_J(N)^{(*)},\;N\in P_J^+.  \label{mgl2}\\
    m^*\left(\frac{\phi(\ins v_N)}{\theta_N}\res{D_{G/U}(J)}\right) = \phi(\ins v_N)\res{(U^-)^J}\otimes 1\otimes 1   &\mb{ for }&   \phi\in L(N)^{(*)},\;N\in P^+\cap\overline{F_J}. \label{mgl3}
\end{eqnarray}

We first show that $m$ is well defined as a map to $G/U$. Let $u^-\in (U^-)^J$, $gU_J=g'U_J\in G_J/U_J$, and $t\in T^J$. Then $u^-gtU=u^-g'tU$ if and only if $gtUt^{-1}=g'tUt^{-1}$ if and only if $gU=g'U$, which follows from $gU_J=g'U_J$.

Next we show that $m$ is a surjective map onto $D_{G/U}(J)$. Let $gU\in G/U$. Write $gU$ in the form $gU=u^-n_\sigma U$ with $u^-\in U^-$, $n_\sigma\in N$ projecting to $\sigma\in\We$. Then $gU\in D_{G/U}(J)$ if and only if
\begin{eqnarray*}
  \delta_\La(u^- n_\sigma v_\La) \neq 0   \quad\mb{ for all }\quad  \La\in P^+\cap\overline{F_J} , 
\end{eqnarray*}
if and only if $\sigma \La =\La $ for all $\La\in P^+\cap\overline{F_J}$, if and only if $\sigma\in\We_J$. Therefore, 
\begin{eqnarray*}
  D_{G/U}(J)=U^- N_J T^J/U=(U^-)^J U_J^- N_J T^J/U = m((U^-)^J,\,G_J/U_J,\,T^J). 
\end{eqnarray*}

Now let $\phi\in L(N)^*_\la $, $\la\in P(N)$, and $N\in P^+$. Choose dual bases $\phi_{\mu\,j}\in L(N)_\mu^*$, $v_{\mu\,j}\in L(N)_\mu$, $\mu\in P(N)$, $j=1,\,\ldots,\,\mb{dim}(L(N)_\mu)$. Decompose $N = N_J + N^J$ with $N_J\in P^+_J$ and $N^J\in P^+\cap\overline{F_J}$. 
Let $\La\in P^+\cap\overline{F_J}$. Note that $G_J\subseteq Z_G (L(\La)_\La)$. 
For all $u^-\in (U^-)^J$, $g\in G_J$, and $t\in T^J$ we find
\begin{eqnarray}
 && m^*\left(\frac{\phi(\ins v_N)}{\theta_\La}\res{D_{G/U}(J)}\right)(u^-,\,gU,\,t)=\frac{\phi(u^-gt v_N)}{\delta_\La(u^-gtv_\La)} 
  =\phi(u^-g v_N)e_{N^J-\La}(t)\qquad\qquad\qquad \label{gggl1}  \\
 &&\qquad\qquad\qquad  =\sum_{\mu,\,j}\phi(u^-v_{\mu\,j})\phi_{\mu\,j}(gv_N) e_{N^J-\La}(t) \nonumber\\
 &&\qquad\qquad\qquad = \bigg(\sum_{\mu,\,j\atop \mu\geq \la} \phi(\ins v_{\mu\,j})\otimes\phi_{\mu\,j}(\ins v_N)\otimes e_{N^J-\La}\bigg)(u^-,\,gU,\,t).  \label{gggl2}
\end{eqnarray}
Since the sum in (\ref{gggl2}) is finite, the comorphism $m^*$ exists. 
Again by $G_J\subseteq Z_G (L(N)_N)$ for all $N\in P^+\cap\overline{F_J}$ we get (\ref{mgl1}) and (\ref{mgl3}) from (\ref{gggl1}). By Proposition \ref{SpMk} we get  (\ref{mgl2}) from (\ref{gggl1}). The elements in (\ref{mgl1}) on the right generate $\F 1\otimes \F 1\otimes \FK{P\cap(\overline{F_J}-\overline{F_J})}$, in (\ref{mgl2}) they generate $\F 1 \otimes \F[G_J]^{U_J}\otimes \F 1$, and in (\ref{mgl3}) they generate $\FK{(U^-)^J} \otimes \F 1\otimes \F 1$. It follows that the morphism of algebras $m^*$ is surjective.
Furthermore, the comorphism $m^*$ is injective, because $m$ is surjective.

Now, by the isomorphisms and identifications of ($\alpha $), ($\beta $), ($\gamma $), ($\delta $) the isomorphism $m^*$ coincides with $\Gamma$. 

The algebra $\widetilde{(P^+\cap\overline{F_J})}^{-1}CA$ is graded by 
\begin{eqnarray*}
    P^+ -(P^+\cap\overline{F_J})\,.
\end{eqnarray*} 
The algebra $\FK{(U_f^-)^J}\otimes CA_J\otimes \FK{P\cap (\overline{F_J}-\overline{F_J})}$ is graded by 
\begin{eqnarray*}
   \{0\}\oplus P_J^+ \oplus (P\cap (\overline{F_J}-\overline{F_J})) \,.
\end{eqnarray*}
Both gradings identify with the submonoid $P\cap (\overline{C}-\overline{F_J})$ of $P$. The isomorphism $\Gamma$ preserves the degree of the generators in (\ref{mg1}), (\ref{mg2}), (\ref{mg3}). Hence $\Gamma$ is an isomorphism of graded algebras, the induced map on the gradings $P\cap (\overline{C}-\overline{F_J})$ given by the identity map.
\qed

As a first consequence of the previous theorem we obtain a description of the standard big cells:
\begin{cor}\label{BigCell2} Let $J\subseteq I$. 
\begin{itemize}
\item[(a)] There exists an isomorphism of graded algebras 
\begin{eqnarray*}
 \Gamma: (\widetilde{P^+\cap\overline{F_J}})^{-1}\bigoplus_{\La\in P^+\cap\overline{F_J}} L(\La)^{(*)}\to \FK{(U^-_f)^J} \otimes \FK{P\cap(\overline{F_J}-\overline{F_J})}
\end{eqnarray*}
such that
\begin{eqnarray*}
\begin{array}{lclcl}
 \Gamma(\frac{\delta_N}{\delta_\La}) &=& 1  \otimes   e_{N-\La} &\;\mb{ for all }\;& N,\,\La\in  P^+\cap\overline{F_J},\\
 \Gamma(\frac{\phi}{\delta_N}) &=& \phi( \ins v_N) \otimes 1 &\;\mb{ for all }\; & \phi\in L(N)^{(*)},\;N\in P^+\cap\overline{F_J}.
\end{array}
\end{eqnarray*}
\item[(b)] The map
\begin{eqnarray*}
  \beta:\,\mb{Spec\,}(\F[(U_f^-)^J]) &\to     & \qquad\qquad\qquad\qquad\qquad \qquad\qquad BC(J)\\  
       Q \quad\qquad\;   &\mapsto & \Big(\Gamma^{-1}\left( Q\otimes \FK{P\cap (\overline{F_J}-\overline{F_J})}\right)\cap \bigoplus_{\La\in P^+\cap\overline{F_J}}L(\La)^{(*)}\Big)\oplus \bigoplus_{\La\in P^+\setminus\overline{F_J}} L(\La)^{(*)}
\end{eqnarray*}
is a homeomorphism, mapping $\mb{Spec\,}(\F[(U_f^-)^J])(\F)$ bijectively to $BC(J)(\F)$.
\end{itemize}
\end{cor}
\Proof  The isomorphism of graded algebras of part (a) is obtained by restricting the graded isomorphism of Theorem \ref{Sp4}. Part (b) follows from (a) and from Proposition \ref{BigCell2a} and Proposition \ref{prodlattice}.
\qed

For $J\subseteq I$ identify the topological space $\overline{Or(J)}$ with $\mb{Proj\,}(\bigoplus_{\La\in P^+\cap\overline{F_J}}L(\La)^{(*)})$ by Proposition \ref{OrbitAbschluss2}. By this identification $\overline{Or(J)}$ gets the structure of a locally ringed space of $\F$-algebras. Because $Or(J)$ is covered by big cells we obtain:

\begin{cor} Let $J\subseteq I$. The open set $Or(J)$ of $\overline{Or(J)}$ equipped with its locally ringed substructure is a scheme. 
\end{cor}
\Proof This follows from Theorem \ref{BigCell1}, from Corollary \ref{BigCell2a},  Corollary \ref{principalscheme}, and from Corollary \ref{BigCell2}, Corollary \ref{prodlatticescheme}.
\qed

With Corollary \ref{BigCell2} we can also determine the $\F$-valued points contained in a standard big cell:
\begin{thm}\label{BigCell3} It is $BC(J)(\F)=(U_f^-)^J P(J)$ for all $J\subseteq I$.
\end{thm}
\Proof By Proposition \ref{UniSpecF} the vanishing ideals $I(u)$, $u\in (U_f^-)^J $, give the $\F$-valued points of $\FK{(U_f^-)^J} $. By Corollary  \ref{BigCell2} (b) its images under the map $\beta $ give the $\F$-valued points of $BC(J)$. 
For $u\in (U_f^-)^J$ we now compute $\beta(I(u))$ explicitly. Let $\phi\in L(\La)^{(*)}$, $\La\in P^+\cap\overline{F_J}$. Then 
\begin{eqnarray*}
  \frac{\phi}{\delta_0}\in \Gamma^{-1}\left( I(u)\otimes \FK{P\cap (\overline{F_J}-\overline{F_J})}\right)
\end{eqnarray*}
if and only if
\begin{eqnarray*}
\Gamma(\frac{\phi}{\delta_0})  = \phi(\ins v_\La)\otimes e_\La\in I(u)\otimes \FK{P\cap (\overline{F_J}-\overline{F_J})}
\end{eqnarray*}
if and only if
\begin{eqnarray*}
   0=\phi(uv_\La)=(\pi(u)\phi)(v_\La)
\end{eqnarray*}
if and only if $\pi(u)\phi\in L(\La)^{(*)}_{\neq\La}$. Since $\Gamma^{-1}\left( I(u)\otimes \FK{P\cap (\overline{F_J}-\overline{F_J})}\right)$ is a homogeneous ideal we have shown
\begin{eqnarray*}
\beta(I(u))= \bigoplus_{\La\in P^+\cap\overline{F_J}}\pi(u)^{-1}L(\La)^{(*)}_{\neq \La}\oplus \bigoplus_{\La\in P^+\setminus\overline{F_J}} L(\La)^{(*)} =uP(J).
\end{eqnarray*}
\qed 

Now we use again the covering of the Or-strata by big cells and finally reach:
\begin{cor}\label{surjofOmega} It is
\begin{eqnarray*}
  \PCAF \;=\;\dot{\bigcup_{J\subseteq I}}\,Or(J)(\F) \quad\mb{ with }\quad Or(J)(\F) \,=\, G_{fn}\, P(J) \quad\mb{ for all }\; J\subseteq I\,.
\end{eqnarray*}
In particular, the map $\omega$ of Theorem \ref{CAB1} is also surjective.
\end{cor}
\Proof The decomposition of $\PCAF$ is given by Proposition \ref{Orbit2} (c). From the previous theorem and Theorem \ref{BigCell1} (b) we get
\begin{eqnarray*}
   Or(J)(\F)= \bigcup_{g\in G_{fn}} g BC(J)(\F)=G_{fn}(U_f^-)^J P(J)=G_{fn}\, P(J).
\end{eqnarray*}
\qed

At last we investigate the transversal behaviour of the Or-stratification of $\PCA$. 
For $J\subseteq I$ we get a stratification of $\mb{Proj\,}(CA_J)$, similarly as for $\PCA$, as follows: For every $K\subseteq J$ it is 
$\bigoplus_{\La\in P^+_J\setminus P^+_{J\setminus K}}L_J(\La)^{(*)}$ a homogeneous prime ideal of $CA_J$. Set 
\begin{eqnarray*}
 \overline{Or_J(K)}:={\mathcal V}(\bigoplus_{\La\in P^+_J\setminus P^+_{J\setminus K}}L_J(\La)^{(*)}) \subseteq \mb{Proj\,}(CA_J).
\end{eqnarray*}
For $K\subseteq J$ define
\begin{eqnarray*}
 Or_J(K):=\overline{Or_J(K)}\setminus\bigcup_{K\subsetneqq T\subseteq J}\overline{Or_J(T)}\,.
\end{eqnarray*}
Then $Or_J(K)\neq \emptyset$, $Or_J(K)$ is open and dense in $\overline{Or_J(K)}$, and $\mb{Proj\,}(CA_J)=\dot{\bigcup}_{K\subseteq J}Or_J(K)$.

As a second consequence of Theorem \ref{Sp4} we obtain: 
\begin{cor}\label{stsp} Let $J\subseteq I$. 
\begin{itemize}
\item[(a)] There exists an isomorphism of graded algebras
\begin{eqnarray*}
 \Gamma: (\widetilde{P^+\cap\overline{F_J}})^{-1} \bigoplus_{\La\in P^+} L_J(\La)^{(*)} \to  CA_J \otimes \FK{P\cap(\overline{F_J}-\overline{F_J})}
\end{eqnarray*}
such that
\begin{eqnarray*}
 \Gamma(\frac{\delta_N}{\delta_\La}) =  1 \otimes   e_{N-\La}  \;\: &\;\mb{ for all }\;&  N,\,\La\in P^+\cap\overline{F_J},\\
 \Gamma(\phi) = \phi \otimes   1  \qquad\quad &\;\mb{ for all }\;&    \phi\in L_J(N)^{(*)},\;N\in P_J^+. 
\end{eqnarray*}
\item[(b)] The map
\begin{eqnarray*}
  \gamma:\,\mb{Proj\,}(CA_J) &\to     & \qquad\qquad\qquad\qquad \qquad\qquad S(J)\\  
       Q \qquad\;   &\mapsto & \Big(\Gamma^{-1}\left( Q\otimes \FK{P\cap (\overline{F_J}-\overline{F_J})}\right)\cap \bigoplus_{\La\in P^+ }L_J(\La)^{(*)}\Big)\oplus \bigoplus_{\La\in P^+} R_J(\La)^{(*)}
\end{eqnarray*}
is a homeomorphism, mapping $\mb{Proj\,}(CA_J)(\F)$ bijectively to $S(J)(\F)$. Furthermore, 
\begin{eqnarray}\label{Orrel}
  \gamma(Or_J(K))= S(J)\cap Or(K) \quad\mb{ for all }\quad K\subseteq J.
\end{eqnarray}
\end{itemize}
\end{cor}
\Proof  The isomorphism of graded algebras of part (a) is obtained by restricting the graded isomorphism of Theorem \ref{Sp4}. Part (b), except (\ref{Orrel}), follows from (a) and from Proposition \ref{Slice2} and Proposition \ref{prodlattice}. To prove (\ref{Orrel}) we first show 
\begin{eqnarray}\label{Orrel2}
  \gamma(\overline{Or_J(K)})= S(J) \cap \overline{Or(K)} \quad\mb{ for all }\quad K\subseteq J.
\end{eqnarray}
Because of
\begin{eqnarray}\label{orvan}
\gamma(\overline{Or_J(K)}) = \gamma( {\mathcal V}(\bigoplus_{\La\in P^+_J\setminus P^+_{J\setminus K}}L_J(\La)^{(*)}) ) 
                = {\mathcal V}(\gamma(\bigoplus_{\La\in P^+_J\setminus P^+_{J\setminus K}}L_J(\La)^{(*)}))\cap S(J) ,
\end{eqnarray}
we first determine $\gamma(\bigoplus_{\La\in P^+_J\setminus P^+_{J\setminus K}}L_J(\La)^{(*)})$. It is $P^+=P^+_J\oplus (P^+\cap\overline{F_J})$. Let $\La_1\in P_J^+$ and $\La_2\in P^+\cap\overline{F_J}$. From the definition of the map $\gamma$ it follows easily that
\begin{eqnarray*}
 \gamma(\bigoplus_{\La\in P^+_J\setminus P^+_{J\setminus K}}L_J(\La)^{(*)})_{\La_1+\La_2} = 
    \left\{\begin{array}{cll}
             L_J(\La_1)^{(*)} \prdca \delta_{\La_2} + R_J(\La_1+\La_2)^{(*)} & \mb{ if } & \La_1\in P^+_J\setminus P^+_{J\setminus K}\\
                    R_J(\La_1+\La_2)^{(*)}  & \mb{ if } & \La_1\in P^+_{J\setminus K}
           \end{array}\right. \,.
\end{eqnarray*}
Here, $L_J(\La_1)^{(*)} \prdca \delta_{\La_2}\subseteq L_J(\La_1+\La_2)^{(*)}$ by Proposition \ref{smdu} (a). To show equality note that the $\g_J$-module $L_J(\La_2)$ is trivial. It follows that $(G_J)^{op}$ acts trivially on $L_J(\La_2)^{(*)}=\F\delta_{\La_2}$.  We find
\begin{eqnarray*}
  \pi(g)\delta_{\La_1+\La_2} = \pi(g)(\delta_{\La_1} \prdca \delta_{\La_2}) =(\pi(g)\delta_{\La_1}) \prdca (\pi(g)\delta_{\La_2}) =(\pi(g)\delta_{\La_1}) \prdca \delta_{\La_2}\in L_J(\La_1)^{(*)} \prdca \delta_{\La_2}
\end{eqnarray*}
for all $g\in G_J$. Since the $(G_J)^{op}$-orbit of $\delta_{\La_1+\La_2}$ spans the  $(G_J)^{op}$-module $L_J(\La_1+\La_2)^{(*)}$ we obtain $L_J(\La_1)^{(*)} \prdca \delta_{\La_2} = L_J(\La_1+\La_2)^{(*)}$.
We have shown
\begin{eqnarray*}
 \gamma(\bigoplus_{\La\in P^+_J\setminus P^+_{J\setminus K}}L_J(\La)^{(*)})_{\La_1+\La_2} = 
    \left\{\begin{array}{cll}
             L(\La_1+\La_2)^{(*)}  & \mb{ if } & \La_1\in P^+_J\setminus P^+_{J\setminus K}\\
                    R_J(\La_1+\La_2)^{(*)}  & \mb{ if } & \La_1\in P^+_{J\setminus K}
           \end{array}\right. \,.
\end{eqnarray*}
With $P^+_{J\setminus K}+ (P^+\cap\overline{F_J})=P^+\cap\overline{F_K}$ we get
\begin{eqnarray*}
  \gamma(\bigoplus_{\La\in P^+_J\setminus P^+_{J\setminus K}}L_J(\La)^{(*)})
   = \bigoplus_{\La\in P^+\cap\overline{F_K}}R_J(\La)^{(*)} \oplus \bigoplus_{\La\in P^+\setminus\overline{F_K}}L(\La)^{(*)} \,.
\end{eqnarray*}
Now we insert in (\ref{orvan}) and we use that $S(J)\subseteq {\mathcal V }(\bigoplus_{\La\in P^+}R_J(\La)^{(*)})$  by definition: 
\begin{eqnarray*}
  && \gamma(\overline{Or_J(K)}) 
        \,=\, \Big\{Q\in S(J) \,\Big|\, Q\supseteq  \bigoplus_{\La\in P^+\cap\overline{F_K}}R_J(\La)^{(*)} \oplus \bigoplus_{\La\in P^+\setminus\overline{F_K}}L(\La)^{(*)} \Big\}\\
   &&    = \Big\{Q\in S(J) \,\Big|\, Q\supseteq  \bigoplus_{\La\in P^+\setminus\overline{F_K}}L(\La)^{(*)}\Big\} \,=\, S(J)\cap \overline{Or(K)}\,.
\end{eqnarray*}

From the definition of $Or_J(K)$ and (\ref{Orrel2}) we get
\begin{eqnarray}\label{OrSRel}
 \gamma(Or_J(K)) = (S(J)\cap \overline{Or(K)})\setminus \bigcup_{K\subsetneqq T\subseteq J}(S(J)\cap \overline{Or(T)})\,.
\end{eqnarray}
Suppose that $T\subseteq I$, $T\not\subseteq J$. Then by Proposition \ref{Orbit2} (c) and Proposition \ref{stratsl} we find
\begin{eqnarray*}
   S(J)\cap\overline{Or(T)} = \bigcup_{ I\supseteq L \supseteq T}(S(J)\cap Or(L))=\emptyset\, .
\end{eqnarray*} 
Therefore, the union in (\ref{OrSRel}) can be taken over all $T$ with $K\subsetneqq T\subseteq I$. This shows  (\ref{Orrel}).
\qed

Next we determine the $\F$-valued points of the standard transversal slices $S(J)$, $J\subseteq I$. We do not work with the previous corollary but with the description of the $\F$-valued points of $\PCA$ instead. In this way, we also find the $\F$-valued points of the principal open sets $D(J)$, $J\subseteq I$.
\begin{thm}\label{posF} For $J\subseteq I$ we have
\begin{eqnarray*}
   D(J)(\F) = \dot{\bigcup_{K\subseteq J}}(U^-_f)^J (G_{fn})_J P(K) = \dot{\bigcup_{K\subseteq J}}(G_{fn})_J (U^-_f)^J P(K)   \quad \mb{ and } \quad S(J)(\F) = \dot{\bigcup_{K\subseteq J}} (G_{fn})_J P(K)\,.
\end{eqnarray*}
\end{thm}
\Proof (a) By Corollary \ref{surjofOmega}, by the Birkhoff decomposition of $G_{fn}$, and by $U\subseteq (P_{fn})_K = \mb{Stab}_{G_{fn}}(P(K))$ for all $K\subseteq I$, compare Theorem \ref{CAB1}, we obtain
\begin{eqnarray*}
    \PCAF = \dot{\bigcup_{K\subseteq I}}G_{fn}\, P(K) = \dot{\bigcup_{K\subseteq I}}U^-_f  N \,P(K)\,.  
\end{eqnarray*}
We use Remark \ref{PrincipalOpen1} to describe $D(J)(\F)$. Fix an arbitrary $\La\in P^+\cap F_J$. For $K\subseteq I$, $u\in U_f^-$, and $n_\sigma\in N$ projecting to $\sigma\in\We$ we find
\begin{eqnarray*}
 && \delta_\La \;\not\in \; u n_\sigma P(K) = \bigoplus_{N\in P^+\cap\overline{F_K}}\pi(u n_\sigma)^{-1} L(N)_{\neq \La}^{(*)}\oplus 
        \bigoplus_{N\in P^+\setminus\overline{F_K}}L(N)^{(*)} \\
 && \iff\quad \La\in P^+\cap\overline{F_K} \;\mb{ and }\; \delta_\La\not\in \pi(u n_\sigma)^{-1} L(\La)^{(*)}_{\neq \La}   
                  \quad\iff\quad J\supseteq K  \;\mb{ and }\;\delta_\La(u n_\sigma v_\La) \neq  0 \\
 && \iff\quad J\supseteq K  \;\mb{ and }\; \sigma\La=\La \quad\iff\quad   J\supseteq K  \;\mb{ and }\; \sigma\in\We_J.
\end{eqnarray*}
Because of $T,\, U_J\subseteq (P_{fn})_K = \mb{Stab}_{G_{fn}}(P(K))$ for all $K\subseteq J \subseteq I $ we get
\begin{eqnarray*}
   U_f^- (\We_J T) P(K) =  (U_f^-)^J (U_f^-)_J N_J U_J P(K) = (U_f^-)^J (G_{fn})_J P(K).
\end{eqnarray*}
From the Levi decomposition of $(P_{fn}^-)_J$ it follows that $(G_{fn})_J \subseteq (L_{fn})_J$ normalizes $(U_f^-)^J$. Therefore, $(U^-_f)^J (G_{fn})_J = (G_{fn})_J (U^-_f)^J $.

(b) By definition, $S(J)= {\mathcal V}(\bigoplus_{\La\in P^+} R_J(\La)^{(*)})\cap D(J)$. We have seen in (a) that the elements of
$D(J)(\F)$ are of the form $g u P(K)$ with $g\in (G_{fn})_J$, $u\in (U_f^-)^J$, and $K\subseteq J $. Because $R_J(\La)^{(*)}$ is invariant under the action of $(G_{fn})_J^{op}$ for all $\La\in P^+$ we find
\begin{eqnarray*}
 && guP(K)   \;\supseteq\; \bigoplus_{\La\in P^+} R_J(\La)^{(*)} \\
  &&\iff  u P(K) = \bigoplus_{\La\in P^+\cap\overline{F_K}}\pi(u)^{-1}L(\La)_{\neq \La}^{(*)}\oplus 
        \bigoplus_{\La\in P^+\setminus\overline{F_K}}L(\La)^{(*)}
        \;\supseteq\; \bigoplus_{\La\in P^+} R_J(\La)^{(*)} \\
  && \iff \pi(u)^{-1}L(\La)_{\neq \La}^{(*)} = \Mklz{ \phi\in L(\La)^{(*)} }{ \phi (u v_\La) = 0 }\supseteq R_J(\La)^{(*)} \;\;\mb{ for all }\;\; \La\in P^+\cap \overline{F_K}\\
  &&\iff   u v_\La \in  L_J(\La)_f = \prod_{\la\,\in\, \La - (Q^+_0)_J} L(\La)_\la  \;\;\mb{ for all }\;\; \La\in P^+\cap \overline{F_K}\, .
\end{eqnarray*}
Because of $u\in (U_f^-)^J = exp((\n_f^-)^J)$ with $(\n_f^-)^J=\bigoplus_{\al\in (\Delta^-)^J}\g_\al$ this is equivalent to
\begin{eqnarray*}
  u L(\La)_\La \subseteq \prod_{\la\,\in\, \La - (\{0\}\cup Q^+_0\setminus(Q^+_0)_J) }L(\La)_\la  \cap  \prod_{\la\,\in\, \La - (Q^+_0)_J} L(\La)_\la \,=\, L(\La)_\La
\end{eqnarray*}
for all $\La\in P^+\cap \overline{F_K}$. These inclusions are equalities, because the highest weight spaces are one-di\-men\-sion\-al. By Theorem \ref{ngfnllala} (b), Theorem \ref{spfpp} (b), the definition of $(P_{fn})_K$, and  Lemma \ref{prodint} this in turn is equivalent to 
\begin{eqnarray*}
 u\in (U_f^-)^J \cap \bigcap_{T\supseteq K} (P_{fn})_T = (U_f^-)^J \cap (P_{fn})_K = (U_f^-)^J \cap (U_f^-)_K(\We_K T)U = (U_f^-)^J\cap (U_f^-)_K.
\end{eqnarray*} 
Since $K\subseteq J$ we obtain $u\in (U_f^-)^J\cap (U_f^-)_K \subseteq  (U_f^-)^J\cap (U_f^-)_J =\{1\}$.
\qed

Let $J\subseteq I$. The canonical injections
\begin{eqnarray*}
  \FK{(U_f^-)^J}\to \FK{(U_f^-)^J}\otimes\, CA_J \quad \mb{ and }\quad CA_J  \to \FK{(U_f^-)^J}\otimes\, CA_J
\end{eqnarray*}
induce a pair of continuous maps 
\begin{eqnarray*}
   \mb{Proj\,}(\FK{(U_f^-)^J}\otimes\, CA_J)\to \mb{Spec\,}(\FK{(U_f^-)^J} )\times \mb{Proj\,}(CA_J)\,.
\end{eqnarray*}
The projective spectrum $\mb{Proj\,}(\FK{(U_f^-)^J}\otimes\, CA_J)$ identifies by Theorem \ref{Sp4} with the principal open set $D(J)$ .  The spectrum $\mb{Spec\,}(\FK{(U_f^-)^J})$  identifies by Corollary \ref{BigCell2} with the standard big cell $BC(J)$.  The projective spectrum $\mb{Proj\,}(CA_J)$ identifies by Corollary \ref{stsp} with the standard transversal slice $S(J)$. Therefore, we get a pair of continuous maps
\begin{eqnarray*}
  p=(p_1,p_2):\,D(J)\to BC(J)\times S(J),
\end{eqnarray*}
whose properties are investigated next. Recall also: By definition, $BC(J)$ and $S(J)$ are closed subsets of $D(J)$. By Theorem \ref{Slice1}, $BC(J)\cap S(J)=\{P(J)\}$.

We first describe  the maps $p_1$ and $p_2$ explicitly:
\begin{prop}\label{desp1p2} Let $J\subseteq I$. The map $p_1:\, D(J)\to BC(J)$ is given by 
\begin{eqnarray*}
    p_1(Q)\;=\; \bigoplus_{\La\in P^+\cap\overline{F_J}}( Q\cap L(\La)^{(*)} ) \oplus \bigoplus_{\La\in P^+\setminus\overline{F_J}}L(\La)^{(*)}\, ,\quad Q\in D(J)\,.
\end{eqnarray*}
The map $p_2:\, D(J)\to S(J)$ is given by
\begin{eqnarray*}
    p_2(Q)\;=\; \bigoplus_{\La\in P^+} \Big( ( Q\cap L_J(\La)^{(*)} ) \oplus R_J(\La)^{(*)} \Big) \, ,\quad Q\in D(J)\,.
\end{eqnarray*}
\end{prop}
\Proof We only show the statement for $p_1$, the statement for $p_2$ follows similarly. Consider the following commutative diagram:
\begin{eqnarray*}
\xymatrix{
\FK{(U^-_f)^J}\otimes CA_J \ar[r] &  \FK{(U^-_f)^J}\otimes CA_J\otimes \FK{P\cap (\overline{F_J}-\overline{F_J})} \ar[r] &  (\widetilde{P^+\cap\overline{F_J}})^{-1} CA\\
  \FK{(U^-_f)^J} \ar[r]\ar[u] & \FK{(U^-_f)^J} \otimes \FK{P\cap (\overline{F_J}-\overline{F_J})} \ar[r]\ar[u]  & (\widetilde{P^+\cap\overline{F_J}})^{-1}\bigoplus_{\La\in P^+\cap \overline{F_J}}L(\La)^{(*)}\ar[u] 
}
\end{eqnarray*} 
Here the maps of the square on the left are the canonical injections of the tensor products. The remaining horizontal maps on the right are the inverse isomorphisms of graded algebras of Theorem \ref{Sp4} (a) and Corollary \ref{BigCell2} (a). The remaining vertical map on the right is the inclusion map, which is a morphism of graded algebras.

This diagram induces a commutative diagram of continuous maps between the corresponding spectra. The two horizontal maps on the left induce homeomorphisms between the corresponding projective spectra by Proposition \ref{prodlattice}. Trivially, the two horizontal maps on the right induce homeomorphisms between the corresponding projective spectra. 

By these homeomorphisms, the map induced by $\FK{(U_f^-)^J} \to \FK{(U_f^-)^J}\otimes\, CA_J$ identifies with the map 
\begin{eqnarray*}
 \mb{Proj\,}((\widetilde{P^+\cap\overline{F_J}})^{-1} CA) & \to & \mb{Proj\,}((\widetilde{P^+\cap\overline{F_J}})^{-1}\bigoplus_{\La\in P^+\cap \overline{F_J}}L(\La)^{(*)}) \\
 \ti{Q}\qquad\qquad\quad &\mapsto & \; \ti{Q}\;\cap \;(\widetilde{P^+\cap\overline{F_J}})^{-1}\bigoplus_{\La\in P^+\cap \overline{F_J}}L(\La)^{(*)}
\end{eqnarray*}
which in turn identifies by Propositions \ref{Sp3} and \ref{BigCell2a} with the map $p_1: D(J)\to BC(J)$. For $Q\in D(J)$ denote by $Q^e$ the ideal in  $(\widetilde{P^+\cap\overline{F_J}})^{-1} CA$ generated by $Q$. Since $Q$ is prime it is $Q^e\cap CA=Q$ . We find
\begin{eqnarray*}
 && p_1(Q) = \Big( \big( Q^e \cap (\widetilde{P^+\cap\overline{F_J}})^{-1}\bigoplus_{\La\in P^+\cap \overline{F_J}}L(\La)^{(*)} \big) \cap \bigoplus_{\La\in P^+\cap \overline{F_J}}L(\La)^{(*)} \Big) \oplus \bigoplus_{\La\in P^+\setminus \overline{F_J}}L(\La)^{(*)}\\
 &&        =  \Big( Q \cap \bigoplus_{\La\in P^+\cap \overline{F_J}}L(\La)^{(*)} \Big) \oplus \bigoplus_{\La\in P^+\setminus \overline{F_J}}L(\La)^{(*)} 
          = \bigoplus_{\La\in P^+\cap \overline{F_J}}(Q\cap L(\La)^{(*)} ) \oplus \bigoplus_{\La\in P^+\setminus \overline{F_J}}L(\La)^{(*)}\;.
\end{eqnarray*}
\qed

The map $p=(p_1,p_2):\,D(J)\to BC(J)\times S(J)$ has the following properties:
\begin{cor} Let $J\subseteq I$. Then:
\begin{itemize}
\item[(a)] $p(Q)= (Q,\,P(J))$ for all $Q\in BC(J)$.
\item[]         $p(Q)= (P(J),\,Q)$ for all $Q\in S(J)$.
\item[(b)] $p(D(J)\cap \overline{Or(K)})\subseteq BC(J)\times (S(J)\cap \overline{Or(K)})$ for all $K\subseteq J$.
\end{itemize}
\end{cor}
\Proof To (a): We only show the second equation, the first follows similarly. Let $Q\in S(J)$. The definition $S(J) = {\mathcal V}(\bigoplus_{\La\in P^+}R_J(\La)^{(*)})\cap D(J)$ implies
\begin{eqnarray}\label{pts1}
    Q\supseteq \bigoplus_{\La\in P^+} R_J(\La)^{(*)}\,.
\end{eqnarray}

From (\ref{pts1}) and the description of $p_1$ in Proposition \ref{desp1p2} we get
\begin{eqnarray*}
  p_1(Q) \supseteq \bigoplus_{\La\in P^+}R_J(\La)^{(*)}\,.
\end{eqnarray*}
By $BC(J)\subseteq D(J)$, by the definition of $S(J)$, and by Theorem \ref{Slice1} we find
\begin{eqnarray*}
  p_1(Q) \in  BC(J)\cap {\mathcal V}( \bigoplus_{\La\in P^+}R_J(\La)^{(*)} )  = BC(J)\cap S(J) = \{ P(J) \}\,.
\end{eqnarray*}

From (\ref{pts1}) and the description of $p_2$ in Proposition \ref{desp1p2} we get
\begin{eqnarray*}
  Q = \bigoplus_{\La\in P^+}\left( (Q\cap L_J(\La)^{(*)}) \oplus R_J(\La)^{(*)} \right) = p_2(Q).
\end{eqnarray*}

To (b): Let $K\subseteq J$. Let $Q\in D(J)$ with $ Q\supseteq \bigoplus_{\La\in P^+\setminus \overline{F_K}}L(\La)^{(*)}$. By the description of $p_2$ in Proposition \ref{desp1p2} we obtain $p_2(Q)\supseteq \bigoplus_{\La\in P^+\setminus \overline{F_K}}L(\La)^{(*)}$.
\qed

On the $\F$-valued points there holds the following description:
\begin{thm} Let $J\subseteq I$. The map $p:\,D(J)\to BC(J)\times S(J)$ restricts to a bijective map $p:\,D(J)(\F)\to BC(J)(\F)\times S(J)(\F)$, which is given by
\begin{eqnarray*}
    p(u g P(K)) = (uP(J),\,gP(K)) \quad\mb{ with }\quad u\in (U_f^-)^J,\; g\in (G_{fn})_J,\; K\subseteq J.
\end{eqnarray*}
In  particular, we have $p\big((D(J)\cap Or(K))(\F)\big) = BC(J)(\F)\times (S(J)\cap Or(K))(\F)$ for all $K\subseteq J$.
\end{thm}
\Proof Let $u\in (U_f^-)^J$, $g\in (G_{fn})_J$, and $K\subseteq J$. If we apply $p_1$ to 
\begin{eqnarray*}
  u g P(K) = \bigoplus_{\La\in P^+\cap\overline{F_K}}\pi(ug)^{-1}L(\La)_{\neq \La}^{(*)}\oplus\bigoplus_{\La\in P^+\setminus \overline{F_K}}L(\La)^{(*)}
\end{eqnarray*}
and use $\overline{F_K}\supseteq \overline{F_J}$ we find
\begin{eqnarray*}
  p_1( u g P(K) ) = \bigoplus_{\La\in P^+\cap\overline{F_J}}\pi(ug)^{-1}L(\La)_{\neq \La}^{(*)}\oplus\bigoplus_{\La\in P^+\setminus \overline{F_J}}L(\La)^{(*)} = ug P(J).
\end{eqnarray*}
Because of $\mb{Stab}_{G_{fn}}(P(J))=(P_{fn})_J\supseteq (G_{fn})_J$, compare Theorem \ref{CAB1}, we get $ p_1( u g P(K) ) = u g P(J) = u P(J)$.

If we apply $p_2$ to $u g P(K)$ we obtain
\begin{eqnarray*}
  p_2(u g P(K)) = 
     \bigoplus_{\La\in P^+\cap\overline{F_K}}\Big( \big(\pi(ug)^{-1}L(\La)_{\neq \La}^{(*)} \cap L_J(\La)^{(*)}\big)\oplus R_J(\La)^{(*)}\Big) \oplus\bigoplus_{\La\in P^+\setminus \overline{F_K}}L(\La)^{(*)}\;.
\end{eqnarray*}
Let $\La\in P^+\cap\overline{F_K}$. Because of $gv_\La\in L_J(\La)_f $ and Proposition \ref{SpMk} we find
\begin{eqnarray*}
  \big(\pi(ug)^{-1}L(\La)_{\neq \La}^{(*)} \cap L_J(\La)^{(*)}\big)\oplus R_J(\La)^{(*)} 
         = \Mklz{ \phi\in L_J(\La)^{(*)} }{ \phi(ugv_\La) = 0 } \oplus R_J(\La)^{(*)}  \quad\quad\quad \\  
          = \Mklz{ \phi\in L_J(\La)^{(*)} }{ \phi(g v_\La) = 0 } \oplus R_J(\La)^{(*)}   
         = \Mklz{ \phi\in L(\La)^{(*)} }{ \phi(g v_\La) = 0 } 
         = \pi(g)^{-1}L(\La)_{\neq \La}^{(*)} \,.
\end{eqnarray*} 
Therefore, we get $ p_2( u g P(K) ) =  g P(K) $.

The bijectivity can be concluded from Theorem \ref{prodF}, but it is also easy to check it directly: The surjectivity follows from Theorem \ref{BigCell3} and Theorem \ref{posF}. To show the injectivity let $p(ugP(K))=p(u'g'P(K'))$  with $u,\,u'\in  (U_f^-)^J$ and $g,\, g'\in  (G_{fn})_J $,  $K,\,K'\subseteq I$. From $uP(J)=u'P(J)$, Theorem \ref{CAB1} and Lemma \ref{prodint} we get
\begin{eqnarray*} 
u^{-1}u'\in \mb{Stab}_{G_{fn}}(P(J))\cap (U_f^-)^J  = (P_{fn})_J \cap  (U_f^-)^J =  (U_f^-)_J \cap  (U_f^-)^J=\{1\}\,.
\end{eqnarray*}
Since $gP(K)=g'P(K')$ we find $ugP(K)=u'g'P(K')$.
\qed
%

%
%
%
%
\section{An action of the face monoid on its building \label{afmb}}
Now we use the results of the previous section to obtain an algebraic geometric model of an action of the face monoid $\GD$ on the building of its unit group, the symmetrizable Kac-Moody group $G$.

We take as building the set
\begin{eqnarray*}
    \Omega := \dot{\bigcup_{J\subseteq I}}G/P_J =\Mklz{g P_J}{g\in G,\;J\subseteq I}
\end{eqnarray*}
partially ordered by the reverse inclusion, i.e., for $g,\,g'\in G$ and $J,\,J'\subseteq I$, 
\begin{eqnarray*}
    g P_J\leq g' P_{J'} \;:\iff\; g P_J\supseteq g' P_{J'}. 
\end{eqnarray*}
The Kac-Moody group $G$ acts order preservingly on $\Omega$ by multiplication from the left. 
We denote by
\begin{eqnarray*}
   {\mathcal A}:=\Mklz{n P_J}{n\in N,\,J\subseteq I}
\end{eqnarray*} 
the standard apartment of $\Omega$. The building $\Omega$ is covered by the apartments $g{\mathcal A}$, $g\in G$.\\

The building $\Omega$ is a substructure of the completed building $\Omega_{fn}$:
\begin{prop}\label{bfb} We get an order preserving, G-equivariant embedding by
\begin{eqnarray*}
     \begin{array}{rccc}
          j:  &\Omega &\to     & \Omega_{fn}\\
              &gP_J   &\mapsto & g(P_{fn})_J
     \end{array}.
\end{eqnarray*}
The standard apartment $\mathcal A$ is mapped bijectively to the standard apartment ${\mathcal A}_{fn}$.
\end{prop}
\Proof This follows immediately from Theorem \ref{ngfnllala} (a) and Theorem \ref{spfpp} (a).\qed
Recall that for $J\subseteq I$ we set
\begin{eqnarray*}
   P(J):=\bigoplus_{\La\in P^+\cap \overline{F_J}} L(\La)^{(*)}_{\neq\La}\oplus \bigoplus_{\La\in P^+\setminus\overline{F_J}} L(\La)^{(*)}   \subseteq CA 
\end{eqnarray*}
with $L(\La)^{(*)}_{\neq\La}:=\bigoplus_{\la\in P(\La)\setminus\{\La\}}L(\La)_\la^*$, a sum over the empty set defined to be $\{0\}$. 
\begin{cor} We get a $G$-equivariant embedding of the $G$-set $\Omega$ into the $G$-set $\PCAF$ by
\begin{eqnarray*}
\begin{array}{rccc}
  \omega: & \Omega &  \to      & \PCAF \\
          & gP_J   &  \mapsto  & g P(J)
\end{array}.
\end{eqnarray*}
Its image is dense in $\PCAF$. Furthermore, for $gP_J,\,hP_K\in\Omega$ we have
\begin{eqnarray*}
  gP_J\leq hP_K   \quad\iff\quad   \overline{\{g P(J)\}}^{\,pts}\cap \omega(\Omega)\subseteq \overline{\{h P(K)\}}^{\,pts}\cap \omega(\Omega).
\end{eqnarray*}
\end{cor}
\Proof The existence of the $G$-equivariant embedding $\omega:\Omega\to\PCAF$ follows from Theorem \ref{CAB1} and Proposition \ref{bfb}. 

For $\omega(\Omega)$ to be dense in $\PCAF$ it is sufficient that the $G$-orbit of $P(\emptyset)$ is dense in $\PCA$. It holds
\begin{eqnarray*}
   \bigcap_{g\in G}g P(\emptyset) &=& \bigcap_{g\in G} (\bigoplus_{\La\in P^+}\pi(g)^{-1}L(\La)_{\neq \La}^{(*)}) = \bigoplus_{\La\in P^+}(\bigcap_{g\in G}\pi(g)^{-1}L(\La)_{\neq \La}^{(*)})\\
  &=& \bigoplus_{\La\in P^+}\Mklz{\phi\in L(\La)^{(*)}} {\pi(g)\phi\in L(\La)^{(*)}_{\neq \La} \mb{ for all }g\in G}\\
  &=& \bigoplus_{\La\in P^+}\Mklz{\phi\in L(\La)^{(*)}} {\phi(gv_\La)=0 \mb{ for all }g\in G} .
\end{eqnarray*}
For every $\La\in P^+$ the orbit $Gv_\La$ spans $L(\La)$, because $L(\La)$ is an irreducible $G$-module. We find $\bigcap_{g\in G}g P(\emptyset)  =\{0\}$. Therefore, $\overline{G P(\emptyset)}={\mathcal V}(\{0\})=\PCA$.

By Theorem \ref{CAB1} and Proposition \ref{bfb} we have $ gP_J\leq hP_K $ if and only if
\begin{eqnarray}\label{GOmegacl1}
    \overline{\{g P(J)\}}^{\,pts}\subseteq \overline{\{h P(K)\}}^{\,pts},
\end{eqnarray}
from which we get
\begin{eqnarray}\label{GOmegacl2}
\overline{\{g P(J)\}}^{\,pts}\cap \omega(\Omega)\subseteq \overline{\{h P(K)\}}^{\,pts}\cap \omega(\Omega).
\end{eqnarray}
Now suppose that (\ref{GOmegacl2}) holds. Because of 
\begin{eqnarray*}
       gP(J)\in \overline{\{g P(J)\}}^{\,pts}\cap \omega(\Omega)\subseteq \overline{\{h P(K)\}}^{\,pts}\cap \omega(\Omega)\subseteq  \overline{\{h P(K)\}}^{\,pts}
\end{eqnarray*}
we get (\ref{GOmegacl1}).
\qed

The action of face monoid $\GD$ on $\PCAF$ induces an action of $\GD$ on the building $\Omega$:

\begin{thm}\label{CAB4} The image $\omega(\Omega)$ is a $\GD$-invariant subset of $\PCAF$, the action of $\GD$ on $\omega(\Omega)$ obtained as follows:
Let $g_1 e(R(\Th)) g_2\in \GD$ and $h P(J)\in \omega(\Omega)$. Decompose $g_2 h$ in the form
\begin{eqnarray*}
  g_2 h= a n_y   c
\end{eqnarray*}
with $a \in P_{\Th \cup \Th^\bot}^-$, $c\in P_J$, and $n_y \in N$ projecting to $y\in\mb{}^{\Th\cup\Th^\bot}\We^J$. Then 
\begin{eqnarray}\label{ga1}
  \left(g_1 e(R(\Th)) g_2\right)  h P(J) = g_1 p_\Th^-(a) P(\Th\cup J\cup \red{y}).
\end{eqnarray}
If we identify the building $\Omega$ with its image $\omega(\Omega)$, then this action coincides with good action 1 of $\GD$ on $\Omega$ given in Corollary 50 of \cite{M5}.
\end{thm}
\Proof It remains to show formula (\ref{ga1}). By Theorem \ref{CAB1} and by Theorem 9 (b) of \cite{M5} we get
\begin{eqnarray*}
   (g_1 e(R(\Th))g_2 ) h P(J)=g_1 e(R(\Th))a n_y c P(J) = g_1 p_\Th^-(a) e(R(\Th)) n_y P(J) .
\end{eqnarray*}
Therefore, we only have to show 
\begin{eqnarray}\label{g1red}
   e(R(\Th)) n_y P(J)=P(\Th\cup J\cup \red{y}).
\end{eqnarray} 
Let $\La\in P^+\setminus \overline{F_J}$. Trivially, 
\begin{eqnarray}\label{gl1}
  \pi(e(R(\Th))n_y)^{-1} P(J)_\La = \pi(e(R(\Th))n_y)^{-1} L(\La)^{(*)} = L(\La)^{(*)}.
\end{eqnarray}
Let $\La\in P^+\cap\overline{F_J}$. Then 
\begin{eqnarray*}
   \lefteqn{\pi(e(R(\Th))n_y)^{-1} P(J)_\La =   \pi(e(R(\Th))n_y)^{-1} L(\La)^{(*)}_{\neq \La}}&&\\
   &&   =\Mklz{\phi\in L(\La)^{(*)}}{\pi(e(R(\Th))n_y)\phi\in L(\La)^{(*)}_{\neq \La}}\\
   &&   =\Mklz{ \phi\in L(\La)^{(*)} }{ \phi(e(R(\Th))n_y v_\La)=0  }\\
   &&   =\left\{\begin{array}{ccl}
             \Mklz{ \phi\in L(\La)^{(*)} }{ \phi(n_y v_\La)=0  } & \mb{if} & y\La\in R(\Th)\\
              \Mklz{ \phi\in L(\La)^{(*)} }{ \phi(0)=0  } & \mb{if} & y\La\notin R(\Th)
            \end{array}\right.  \\
   &&    =\left\{\begin{array}{ccl}
             \Mklz{ \phi\in L(\La)^{(*)} }{ \phi(n_y v_\La)=0  } & \mb{if} & y\La\in R(\Th)\\
              L(\La)^{(*)}  &\mb{if} & y\La\notin R(\Th)
            \end{array}\right. .
\end{eqnarray*}
By Theorem 7 of \cite{M5} we obtain $R(\Th)\cap  y\overline{F_J}=\overline{F_{\Th\cup J\cup \red{y}}}$. Note also that $y$ fixes the points of $\overline{F_{\Th\cup J\cup \red{y}}}$. We find
\begin{eqnarray*}
  y\La\in R(\Th) \iff y\La \in \overline{F_{\Th\cup J\cup \red{y}}} \iff \La \in \overline{F_{\Th\cup J\cup \red{y}}} .
\end{eqnarray*}
Furthermore, for $\La \in \overline{F_{\Th\cup J\cup \red{y}}}$ we have $0\neq n_y v_\La\in L(\La)_{y\La}=L(\La)_\La$.
Therefore, for $\La\in P^+\cap\overline{F_J}$, we have shown
\begin{eqnarray}\label{gl2}
   \pi(e(R(\Th))n_y)^{-1} P(J)_\La = \left\{\begin{array}{ccl}
              L(\La)^{(*)}_{\neq\La} & \mb{if} & \La\in \overline{F_{\Th\cup J\cup \red{y}}}\\
              L(\La)^{(*)}  &\mb{if} &  \La\notin \overline{F_{\Th\cup J\cup \red{y}}}
            \end{array}\right. .
\end{eqnarray}
Because of $\overline{F_{\Th\cup J\cup \red{y}}}\subseteq \overline{F_J}$,  (\ref{gl1}) and (\ref{gl2}) imply (\ref{g1red}).
\qed

\section{Some open questions\label{quest}}
\begin{itemize}
\item[(a)] In the literature, there are many constructions of flag schemes, flag varieties, and flag manifolds for Kac-Moody groups. Compare the flag schemes $Or(J)$, $J\subseteq I$, with these constructions. In particular:
\begin{itemize}
\item[$\bullet$] Compare the scheme $Or(\emptyset)$ with the thick flag scheme of M. Kashiwara in \cite{Kas}.
\item[$\bullet$] Let $\La\in P^+\cap F_J$, $J\subseteq I$. V. Kac and D. Peterson constructed in \cite{KP1}, \cite{KP2} a thin flag variety as follows. They showed that the Kostant cone ${\mathcal V}_\La:=G (L(\La)_\La)$ is a closed subvariety of the linear space $L(\La)$,
its $\Nn$-graded coordinate ring $\C\,[{\mathcal V}_\La]$ isomorphic to $\bigoplus_{k\in\Nn}L(\La)^{(*)}\subseteq CA$. This realizes the thin flag variety $G/P_J$ as a closed subvariety of the projective space ${\mathbb P}( L(\La) )$. Compare the scheme $\mb{Proj\,}(\C\,[{\mathcal V}_\La])\setminus \{ \mb{ the irrelevant ideal }\}$ with the scheme $Or(J)$.
\end{itemize}
\item[(b)] In the article we investigated the projective spectrum of the Cartan algebra $\mb{Proj\,}(CA)$, on which the formal Kac-Moody group $G_{fn}$ acts by morphisms. In particular, we showed that the $G_{fn}$-space $\mb{Proj\,}(CA)(\F)$ is the union of all {\bf thick flag spaces}, i.e., $\mb{Proj\,}(CA)(\F)\cong\dot{\bigcup}_{J\subseteq I} G_{fn}/(P_{fn})_J$. In general, the subgroups $(P_{fn})_J$, $J\subseteq I$, of $G_{fn}$ are nonparabolic.

To define the Cartan algebra $CA$ the restricted duals $L(\La)^{(*)}$, $\La\in P^+$, have been used. Similarly, it is possible to obtain a Cartan algebra $CA_f$, the restricted duals replaced by the full duals $L(\La)^*$, $\La\in P^+$. 

Investigate the projective spectrum of the Cartan algebra $\mb{Proj\,}(CA_f)$, on which the formal Kac-Moody group $G_{fp}$ acts by morphisms. Is the $G_{fp}$-space $\mb{Proj\,}(CA_f)(\F)$ the union of all {\bf thin flag spaces}, i.e., $\mb{Proj\,}(CA_f)(\F)\cong\dot{\bigcup}_{J\subseteq I} G_{fp}/(P_{fp})_J$ where  $(P_{fp})_J$, $J\subseteq I$, are the standard parabolic subgroups of $G_{fp}$?
\end{itemize}
%
%
%
%
%

%
\end{document}